\documentclass[11pt]{article}

\usepackage{cite}             
\usepackage{amssymb}          
\usepackage{latexsym}         
\usepackage{amsmath}          
\usepackage[latin1]{inputenc} 
\usepackage{eucal}            
\usepackage{bbm}              
\usepackage{theorem}          
\usepackage[ps,dvips]{xypic}  
\usepackage{diagxy}           
\usepackage{paralist}         

\title{The $H$-Covariant Strong Picard Groupoid}

\author{\textbf{Stefan Jansen}\thanks{E-mail:
    Stefan.Jansen@physik.uni-freiburg.de},
  \addtocounter{footnote}{5}
  \textbf{Stefan Waldmann}\thanks{E-mail:
    Stefan.Waldmann@physik.uni-freiburg.de}
  \\[0.1cm]
  Fakult{\"a}t f{\"u}r Mathematik und Physik\\
  Albert-Ludwigs-Universit{\"a}t Freiburg\\
  Physikalisches Institut\\
  Hermann Herder Stra{\ss}e 3\\
  D 79104 Freiburg\\
  Germany}

\date{Revised version April 2005\\[0.5cm] FR-THEP 2004/16}

\renewcommand{\mathbb}[1]{\mathbbm{#1}} 

\newcommand{\textdef}[1] {\textbf{#1}}

\newcommand{\cc}[1]      {\overline{{#1}}}              
\newcommand{\id}         {\operatorname{\mathsf{id}}}   
\newcommand{\image}      {\operatorname{{\mathrm{im}}}} 
\newcommand{\Lie}        {\operatorname{\mathcal{L}}}

\newcommand{\Ad}         {\operatorname{\mathrm{Ad}}}    
\newcommand{\Hom}        {\operatorname{\mathsf{Hom}}}   
\newcommand{\End}        {\operatorname{\mathsf{End}}}   
\newcommand{\SP}[1]      {\left\langle{#1}\right\rangle} 
\newcommand{\ring}[1]    {\mathsf{#1}}                 
\newcommand{\Unit}       {\mathbb{1}}                  
\newcommand{\Aut}        {\operatorname{\mathsf{Aut}}} 
\newcommand{\cl}         {\mathrm{cl}}                    
\newcommand{\I}          {\mathrm{i}}

\newcommand{\group}[1]   {\mathrm{#1}}
\newcommand{\lie}[1]     {\mathfrak{#1}}
\newcommand{\acts}       {\mathbin{\triangleright}}
\newcommand{\actsp}      {\mathbin{\triangleright'}}
\newcommand{\actsT}[1][{}] {\mathbin{\triangleright^{#1}}}
\newcommand{\actsS}[1][{}] {\mathbin{\triangleright_{#1}}}
\newcommand{\ccacts}     {\mathbin{\cc{\triangleright}}}
\newcommand{\sweedler}[1] {{\scriptscriptstyle{(#1)}}}
\newcommand{\op}         {\mathrm{op}}
\newcommand{\twist}[1]   {\mathsf{#1}}
\newcommand{\zentrum}    {\mathcal{Z}}

\newcommand{\Bimod}[5] {\sideset{^{\scriptscriptstyle{#1}}_{\scriptscriptstyle{#2}}}{^{\scriptscriptstyle{#4}}_{\scriptscriptstyle{#5}}}{\operatorname{#3}}}

\newcommand{\EA}   {\Bimod{}{}{\mathcal{E}}{}{\mathcal{A}}}
\newcommand{\EpA}  {\Bimod{}{}{\mathcal{E}}{\prime}{\mathcal{A}}}
\newcommand{\BEA}  {\Bimod{}{\mathcal{B}}{\mathcal{E}}{}{\mathcal{A}}}
\newcommand{\BEpA} {\Bimod{}{\mathcal{B}}{\mathcal{E}}{\prime}{\mathcal{A}}}
\newcommand{\HD}   {\Bimod{}{}{\mathcal{H}}{}{\mathcal{D}}}
\newcommand{\HDbot}{\Bimod{}{}{\mathcal{H}}{\bot}{\mathcal{D}}}
\newcommand{\AHD}  {\Bimod{}{\mathcal{A}}{\mathcal{H}}{}{\mathcal{D}}}
\newcommand{\HpD}  {\Bimod{}{}{\mathcal{H}}{\prime}{\mathcal{D}}}
\newcommand{\CFB}  {\Bimod{}{\mathcal{C}}{\mathcal{F}}{}{\mathcal{B}}}
\newcommand{\CFPhiA}  {\Bimod{}{\mathcal{C}}{\mathcal{F}}{\Phi}{\mathcal{A}}}
\newcommand{\FB}   {\Bimod{}{}{\mathcal{F}}{}{\mathcal{B}}}
\newcommand{\GC}   {\Bimod{}{}{\mathcal{G}}{}{\mathcal{C}}}
\newcommand{\DGDp} {\Bimod{}{\mathcal{D}}{\mathcal{G}}{}{\mathcal{D}^\prime}}
\newcommand{\AccEB}{\Bimod{}{\mathcal{A}}{\cc{\mathcal{E}}}{}{\mathcal{B}}}
\newcommand{\AAA}  {\Bimod{}{\mathcal{A}}{\mathcal{A}}{}{\mathcal{A}}}
\newcommand{\BBB}  {\Bimod{}{\mathcal{B}}{\mathcal{B}}{}{\mathcal{B}}}
\newcommand{\DDD}  {\Bimod{}{\mathcal{D}}{\mathcal{D}}{}{\mathcal{D}}}
\newcommand{\BBPhiA}{\Bimod{}{\mathcal{B}}{\mathcal{B}}{\Phi}{\mathcal{A}}}
\newcommand{\CPsiEA}{\Bimod{\Psi}{\mathcal{C}}{\mathcal{E}}{}{\mathcal{A}}}
\newcommand{\BPhiAA}{\Bimod{\Phi}{\mathcal{B}}{\mathcal{A}}{}{\mathcal{A}}}
\newcommand{\CCPsiB}{\Bimod{}{\mathcal{C}}{\mathcal{C}}{\Psi}{\mathcal{B}}}
\newcommand{\CCPsiPhiA}{\Bimod{}{\mathcal{C}}{\mathcal{C}}{\Psi \circ \Phi}{\mathcal{A}}}
\newcommand{\PhiBEA}  {\Bimod{\Phi}{\mathcal{B}}{\mathcal{E}}{}{\mathcal{A}}}
\newcommand{\BEAH}  {\Bimod{}{\mathcal{B}\!\rtimes\!\mathit{H}}{\mathcal{E}\!\rtimes\!\mathit{H}}{}{\mathcal{A}\!\rtimes\!\mathit{H}}}
\newcommand{\BEpAH}  {\Bimod{}{\mathcal{B}\!\rtimes\!\mathit{H}}{\mathcal{E}^\prime\!\rtimes\!\mathit{H}}{}{\mathcal{A}\!\rtimes\!\mathit{H}}}
\newcommand{\AHAHAH}  {\Bimod{}{\mathcal{A}\!\rtimes\!\mathit{H}}{\mathcal{A}\!\rtimes\!\mathit{H}}{}{\mathcal{A}\!\rtimes\!\mathit{H}}}
\newcommand{\PA}   {\Bimod{}{}{\mathcal{P}}{}{\mathcal{A}}}
\newcommand{\PB}   {\Bimod{}{}{\mathcal{P}}{}{\mathcal{B}}}

\newcommand{\TMult}[1]  {\mathbin{\cdot_{\scriptscriptstyle{#1}}}}
\newcommand{\dotPhi} {\TMult{\Phi}}
\newcommand{\dotPsi} {\TMult{\Psi}}

\newcommand{\IP}[4]{{\,}_{\scriptscriptstyle{#2}\!\!}\left\langle{{#1}}\right\rangle^{\scriptscriptstyle{#3}}_{\scriptscriptstyle{#4}}}

\newcommand{\SPA}[1]     {\IP{{#1}}{}{}{\mathcal{A}}}
\newcommand{\SPD}[1]     {\IP{{#1}}{}{}{\mathcal{D}}}
\newcommand{\SPEA}[1]    {\IP{{#1}}{}{\mathcal{E}}{\mathcal{A}}}
\newcommand{\BSPE}[1]    {\IP{{#1}}{\mathcal{B}}{\mathcal{E}}{}}
\newcommand{\ASP}[1]     {\IP{{#1}}{\mathcal{A}}{}{}}
\newcommand{\SPFB}[1]    {\IP{{#1}}{}{\mathcal{F}}{\mathcal{B}}}

\newcommand{\SPFEA}[1]   {\IP{{#1}}{}{\mathcal{F}\otimes\mathcal{E}}{\mathcal{A}}}
\newcommand{\ASPccE}[1]  {\IP{{#1}}{\mathcal{A}}{\cc{\mathcal{E}}}{}}

\newcommand{\MnASP}[1]   {\IP{{#1}}{M_n(\mathcal{A})}{}{}}
\newcommand{\SPBPhiA}[1] {\IP{{#1}}{}{\mathcal{B}^\Phi}{\mathcal{A}}}
\newcommand{\SPFPhiA}[1] {\IP{{#1}}{}{\mathcal{F}^\Phi}{\mathcal{A}}}
\newcommand{\CSPPsiE}[1] {\IP{{#1}}{\mathcal{C}}{{\!}^{\Psi}\mathcal{E}}{}}
\newcommand{\BSPPhiE}[1] {\IP{{#1}}{\mathcal{B}}{{\!}^{\Phi}\mathcal{E}}{}}
\newcommand{\SPEoHAH}[1] {\IP{{#1}}{}{\mathcal{E}\!\otimes\!\mathit{H}}{\mathcal{A}\!\rtimes\!\mathit{H}}}
\newcommand{\SPEHAH}[1]  {\IP{{#1}}{}{\mathcal{E}\!\rtimes\!\mathit{H}}{\mathcal{A}\!\rtimes\!\mathit{H}}}
\newcommand{\SPFEHAH}[1] {\IP{{#1}}{}{(\mathcal{F} \tensM \mathcal{E})\!\rtimes\!\mathit{H}}{\mathcal{A}\!\rtimes\!\mathit{H}}}
\newcommand{\SPFHEHAH}[1] {\IP{{#1}}{}{(\mathcal{F}\!\rtimes\!\mathit{H}) \tensM (\mathcal{E}\!\rtimes\!\mathit{H})}{\mathcal{A}\!\rtimes\!\mathit{H}}}
\newcommand{\SPFHBH}[1] {\IP{{#1}}{}{\mathcal{F}\!\rtimes\!\mathit{H}}{\mathcal{B}\!\rtimes\!\mathit{H}}}
\newcommand{\AHSPccEH}[1]{\IP{{#1}}{\mathcal{A}\!\rtimes\!\mathit{H}}{\cc{\mathcal{E}}\!\rtimes\!\mathit{H}}{}}
\newcommand{\AHSPCCEH}[1]{\IP{{#1}}{\mathcal{A}\!\rtimes\!\mathit{H}}{\cc{\mathcal{E}\!\rtimes\!\mathit{H}}}{}}
\newcommand{\BHSPEoH}[1]{\IP{{#1}}{\mathcal{B}\!\rtimes\!\mathit{H}}{\mathcal{E}\!\otimes\!\mathit{H}}{}}
\newcommand{\BHSPEH}[1]{\IP{{#1}}{\mathcal{B}\!\rtimes\!\mathit{H}}{\mathcal{E}\!\rtimes\!\mathit{H}}{}}
\newcommand{\SPBHPhiIdAH}[1]{\IP{{#1}}{}{{\mathcal{B}\!\rtimes\!\mathit{H}}^{\Phi\otimes\id}}{\mathcal{A}\rtimes\!\mathit{H}}}
\newcommand{\SPBPhiHAH}[1]{\IP{{#1}}{}{\mathcal{B}^\Phi \rtimes\!\mathit{H}}{\mathcal{A}\rtimes\!\mathit{H}}}
\newcommand{\SPBHBH}[1]{\IP{{#1}}{}{\mathcal{B}\!\rtimes\!\mathit{H}}{\mathcal{B}\!\rtimes\!\mathit{H}}}
\newcommand{\SPPA}[1]    {\IP{{#1}}{}{\mathcal{P}}{\mathcal{A}}}

\newcommand{\tensor}[1][{}]{\mathbin{\otimes_{\scriptscriptstyle{#1}}}}
\newcommand{\tensM}[1][{}] {\mathbin{\widehat{\otimes}_{\scriptscriptstyle{#1}}}}
\newcommand{\tensB}[1][{}] {\mathbin{\widetilde{\otimes}_{\scriptscriptstyle{#1}}}}

\newcommand{\rep}[1][{}]  {\sideset{^*}{_{#1}}{\operatorname{\textrm{-}\mathsf{rep}}}}
\newcommand{\Rep}[1][{}]  {\sideset{^*}{_{#1}}{\operatorname{\textrm{-}\mathsf{Rep}}}}
\newcommand{\smod}[1][{}] {\sideset{^*}{_{#1}}{\operatorname{\textrm{-}\mathsf{mod}}}}
\newcommand{\sMod}[1][{}] {\sideset{^*}{_{#1}}{\operatorname{\textrm{-}\mathsf{Mod}}}}

\newcommand{\Lattice}[1][{}] {\sideset{}{_{#1}}{\operatorname{\mathfrak{L}}}}

\newcommand{\Pic}      {\operatorname{\mathsf{Pic}}}
\newcommand{\PicH}     {\sideset{}{_{H}}{\operatorname{\mathsf{Pic}}}}
\newcommand{\StrPic}   {\sideset{}{^{\mathrm{str}}}{\operatorname{\mathsf{Pic}}}}
\newcommand{\starPic}  {\sideset{}{^*}{\operatorname{\mathsf{Pic}}}}
\newcommand{\StrPicH}  {\sideset{}{^{\mathrm{str}}_{H}}{\operatorname{\mathsf{Pic}}}}
\newcommand{\starPicH} {\sideset{}{^*_{H}}{\operatorname{\mathsf{Pic}}}}

\newcommand{\starIso}  {\sideset{}{^*}{\operatorname{\mathsf{Iso}}}}
\newcommand{\starAut}  {\sideset{}{^*}{\operatorname{\mathsf{Aut}}}}
\newcommand{\starIsoH} {\sideset{}{^*_{H}}{\operatorname{\mathsf{Iso}}}}
\newcommand{\starAutH} {\sideset{}{^*_{H}}{\operatorname{\mathsf{Aut}}}}
\newcommand{\starInnAutH} {\sideset{}{^*_{H}}{\operatorname{\mathsf{InnAut}}}}
\newcommand{\IsoH}     {\sideset{}{_{H}}{\operatorname{\mathsf{Iso}}}}
\newcommand{\AutH}     {\sideset{}{_{H}}{\operatorname{\mathsf{Aut}}}}

\newcommand{\StrProjH}      {\sideset{}{^{\mathrm{str}}_{H}}{\operatorname{\underline{\mathsf{Proj}}}}}
\newcommand{\StrProjHClass} {\sideset{}{^{\mathrm{str}}_{H}}{\operatorname{\mathsf{Proj}}}}
\newcommand{\starProjH}     {\sideset{}{^*_{H}}{\operatorname{\underline{\mathsf{Proj}}}}}
\newcommand{\starProjHClass}{\sideset{}{^*_{H}}{\operatorname{\mathsf{Proj}}}}
\newcommand{\StrKH}         {\sideset{}{^{\mathrm{str}}_{0,H}}{\operatorname{\mathsf{K}}}}
\newcommand{\starKH}        {\sideset{}{^*_{0,H}}{\operatorname{\mathsf{K}}}}

\newtheorem{lemma} {Lemma} [section]
\newtheorem{proposition} [lemma] {Proposition}

\newtheorem{theorem} [lemma] {Theorem}
\newtheorem{corollary} [lemma] {Corollary}
\newtheorem{definition}[lemma] {Definition}

\theorembodyfont{\rmfamily}
\newtheorem{example}[lemma]{Example}
\newtheorem{remark}[lemma]{Remark}

\newenvironment{proof}[1][{}]{
  \par\noindent
  \textsc{Proof{#1}:}
}
{
  \hspace*{\fill} $\blacksquare$\newline
}

\numberwithin{equation}{section}

\setdefaultenum{\itshape i.)}{a.)}{I.)}{A.)}

\begin{document}

\maketitle

\begin{abstract}
    The notion of $H$-covariant strong Morita equivalence is
    introduced for $^*$-algebras over $\ring{C} = \ring{R}(\I)$ with
    an ordered ring $\ring{R}$ which are equipped with a $^*$-action
    of a Hopf $^*$-algebra $H$. This defines a corresponding
    $H$-covariant strong Picard groupoid which encodes the entire
    Morita theory.  Dropping the positivity conditions one obtains
    $H$-covariant $^*$-Morita equivalence with its $H$-covariant
    $^*$-Picard groupoid. We discuss various groupoid morphisms
    between the corresponding notions of the Picard groupoids.
    Moreover, we realize several Morita invariants in this context as
    arising from actions of the $H$-covariant strong Picard groupoid.
    Crossed products and their Morita theory are investigated using a
    groupoid morphism from the $H$-covariant strong Picard groupoid
    into the strong Picard groupoid of the crossed products.
\end{abstract}

\newpage

%
%

\tableofcontents

%
%

\section{Introduction}
\label{sec:Intro}

Morita equivalence is by now in many areas of mathematics an important
tool to compare and relate objects beyond the notion of isomorphism:
the general approach is to enhance a given category by allowing more
general morphisms while keeping the objects. This way, more objects
can become isomorphic in this enhanced category. The idea is that this
`Morita equivalence' of objects implies that the objects have an
equivalent `representation theory'.  Each such enhanced category
specifies its groupoid of invertible morphisms, which usually is
called the corresponding Picard groupoid in this context. This (large)
groupoid encodes then the entire Morita theory.

The list of examples is long, starting with Morita's original version
\cite{morita:1958a} where one considers the category of (unital) rings
with certain bimodules between them as generalized morphisms, see
e.g.\cite{lam:1999a, benabou:1967a, bass:1968a}. Beside various
algebraic refinements of the ring-theoretic notion, notions of Morita
equivalence have been developed also in completely different contexts,
notably for $C^*$-algebras by Rieffel~\cite{rieffel:1974a,
  rieffel:1974b} coining the notion of strong Morita equivalence which
is now one of the crucial ingredients for Connes' noncommutative
geometry \cite{connes:1994a}, for Poisson manifolds by
Xu~\cite{xu:1991a}, see also \cite{bursztyn.weinstein:2004a}, and for
Lie groupoids, see e.g.\cite{moerdijk.mrcun:2003a} and references
therein. We refer to \cite{landsman:2001c} for a comparison of the
later three concepts.

Among the algebraic notions one has $^*$-Morita equivalence for
involutive algebras by Ara~\cite{ara:1999b, ara:1999a} and strong
Morita equivalence of $^*$-algebras over a ring $\ring{C} =
\ring{R}(\I)$, where $\ring{R}$ is an ordered ring and $\I^2 = -1$, by
Bursztyn and Waldmann\cite{bursztyn.waldmann:2003a:pre,
  bursztyn.waldmann:2001a}. These notions provide a bridge between the
ring-theoretic notion and the $C^*$-algebraic framework and
incorporate already many ideas of the later like positivity in a
purely algebraic context. In particular, strong Morita equivalence of
$^*$-algebras was used to study formal star products in deformation
quantization, see \cite{bursztyn:2002a, bursztyn.waldmann:2002a} and
\cite{waldmann:2005a} for a review. Here the $^*$-algebra in
question is a formal associative deformation in the sense of
Gerstenhaber \cite{gerstenhaber.schack:1988a} of the Poisson algebra
of smooth functions $C^\infty(M)$ on a Poisson manifold $M$, which
plays the role of the phase space of a classical mechanical system,
see \cite{bayen.et.al:1978a} and \cite{gutt:2000a,
  dito.sternheimer:2002a} for recent reviews. The deformed algebra is
then interpreted as the observable algebra of the corresponding
quantum system. Since the understanding of the representation theory
is crucial for physical applications one is naturally interested in an
adapted Morita theory in this context. Moreover, it is believed that
Morita equivalence of star products is in some sense the quantum version
of Morita equivalence of the underlying Poisson structures in the
sense of Xu, while on the other hand, formal star products are seen to
be a step into the direction of a $C^*$-algebraic description of the
quantum observables.  Thus one expects relations between all three
types of Morita theory, see e.g.~the discussions in
\cite{bursztyn.weinstein:2004a:pre, landsman:2001c, landsman:1998a,
  landsman:2001b, cannasdasilva.weinstein:1999a}.

As symmetries play a fundamental role in the understanding of
classical and quantum mechanics, it is natural to ask for concepts of
Rieffel induction and Morita equivalence which are compatible with
given symmetries. In the $C^*$-algebraic framework such notions are
well-established for $C^*$-dynamical systems, see
e.g.~\cite[Chap.~7]{raeburn.williams:1998a} and reference therein as
well as \cite{vaes:2004a:pre, kustermans:2002a} for more general
constructions using locally compact quantum groups.

It is the purpose of this work to transfer these ideas from
$C^*$-algebra theory to the more general and algebraic framework of
$^*$-algebras over rings of the form $\ring{C} = \ring{R}(\I)$. This
framework still allows for the crucial notions of positivity but is
wide enough to treat $C^*$-algebras and formal star products on equal
footing. The notion of symmetry we are using is rather general, we
consider $^*$-actions of a Hopf $^*$-algebra $H$ on $^*$-algebras.
After establishing an adapted notion of Morita equivalence, one main
focus of this work is on the resulting notion of the Picard groupoid
and the Morita invariants which are seen to arise from actions of the
Picard groupoid.

The paper is organized as follows: In Section~\ref{sec:Preliminary} we
recall some well-known results on $^*$-algebras, their
$^*$-representation theory on pre-Hilbert modules, and Hopf
$^*$-algebras and their $^*$-actions. This allows us to set up the
basic notions of an $H$-covariant $^*$-representation theory. In
Section~\ref{sec:CovSME} we define the tensor product of bimodules
equipped with inner products and $H$-actions and discuss the
definition of $H$-covariant strong Morita equivalence. We show that it
is indeed an equivalence relation. Moreover, $H$-covariantly strongly
Morita equivalent $^*$-algebras are shown to have equivalent
$H$-covariant $^*$-representation theories on $H$-covariant
pre-Hilbert modules. Section~\ref{sec:CovPicardInvariants} is devoted
to the definition of the $H$-covariant strong Picard groupoid
$\StrPicH$. We discuss several natural groupoid morphisms in this
context, in particular the `forgetful' groupoid morphism $\StrPicH
\longrightarrow \StrPic$ into the strong Picard groupoid. We obtain a
full characterization of its kernel.
Section~\ref{subsec:ActionsInvariants} illustrates the principle that
Morita invariants are obtained from actions of the Picard groupoid. We
discuss the representation theories, the $H$-invariant central
elements, the $H$-equivariant strong $K_0$-groups, the lattices of
$(\mathcal{D}, H)$-closed ideals as well as the groups used to
classify the (inequivalent) $H$-actions on equivalence bimodules. In
Section~\ref{sec:CrossedProducts} we investigate crossed products and
establish a groupoid morphism from the $H$-covariant strong Picard
groupoid into the strong Picard groupoid of the crossed products,
proving thereby in particular that crossed products are strongly
Morita equivalent if the underlying $^*$-algebras are $H$-covariantly
strongly Morita equivalent, a theorem well-known in $C^*$-algebra
theory. Finally, Appendix~\ref{sec:FunnyGroup} contains the
construction of the groups used in the characterization of $H$-actions
on equivalence bimodules.

\medskip

\noindent
\textbf{Acknowledgments:} We would like to thank Martin Bordemann,
Henrique Bursztyn and Giuseppe Dito for valuable discussions and
remarks. Moreover, S. W.  would like to thank the University of Dijon
for the warm hospitality during his stay where parts of this work was
done. Finally, we would like to thank the referee for his detailed
comments and helpful suggestions.

%
%

\section{Preliminary results}
\label{sec:Preliminary}

In this section we recall some basic definitions and results from
representation theory of $^*$-algebras over ordered rings and Hopf
algebra theory in order to make this work self-contained and to set up
our notation, see \cite{bursztyn.waldmann:2003a:pre,
  bursztyn.waldmann:2001a, bursztyn.waldmann:2001b, waldmann:2005a}
for details on $^*$-algebras over ordered rings,
e.g.~\cite{kadison.ringrose:1997a, landsman:1998a, schmuedgen:1990a}
for the representation theory of $C^*$-algebras and $^*$-algebras over
$\mathbb{C}$ and \cite{majid:1995a, kassel:1995a,
  klimyk.schmuedgen:1997a, schmuedgen.wagner:2003a} for Hopf
$^*$-algebras.

%
%

\subsection{$^*$-Algebras over ordered rings}
\label{subsec:AlgebrasOrderedRings}

Let $\ring{R}$ be an \textdef{ordered ring} and let $\ring{C} =
\ring{R}(\I)$ with $\I^2 = -1$. Motivated by deformation quantization,
the main examples we have in mind are $\ring{R} = \mathbb{R}$ and
$\ring{R} = \mathbb{R}[[\lambda]]$ with their natural ordering
structures. Then a \textdef{$^*$-algebra} $\mathcal{A}$ over
$\ring{C}$ is an associative algebra over $\ring{C}$ with an
involutive $\ring{C}$-antilinear antiautomorphism, called the
\textdef{$^*$-involution}, which we shall denote by $a \mapsto a^*$
for $a \in \mathcal{A}$.

A linear functional $\omega: \mathcal{A} \longrightarrow \ring{C}$ is
called \textdef{positive} if $\omega(a^*a) \ge 0$ for all $a \in
\mathcal{A}$. This allows to define \textdef{positive algebra
  elements} $a \in \mathcal{A}$ by the requirement $\omega(a) \ge 0$
for all positive linear functionals. The set of positive algebra
elements is denoted by $\mathcal{A}^+$. Clearly, elements of the form
$\alpha_1 a_1^*a_1 + \cdots + \alpha_n a_n^*a_n$ are positive where $0
< \alpha_i \in \ring{R}$ and $a_i \in \mathcal{A}$. These elements
will be denoted by $\mathcal{A}^{++}$.  See \cite{schmuedgen:1990a}
for more general positive wedges and \cite{waldmann:2004a} for a
comparative discussion of these concepts of positive algebra elements.

A basic example of a $^*$-algebra is obtained as follows: a
\textdef{pre-Hilbert space} $\mathcal{H}$ is a $\ring{C}$-module with
a $\ring{C}$-valued sesquilinear inner product (linear in the second
argument) $\SP{\cdot,\cdot}: \mathcal{H} \times \mathcal{H}
\longrightarrow \ring{C}$ satisfying $\SP{\phi,\psi} =
\cc{\SP{\psi,\phi}}$ for $\phi,\psi \in \mathcal{H}$ and
$\SP{\phi,\phi} > 0$ for $\phi \ne 0$. Then a linear endomorphism $A
\in \End(\mathcal{H})$ is called \textdef{adjointable} if there is an
\textdef{adjoint} $A^*$ with $\SP{A \phi, \psi} = \SP{\phi, A^*\psi}$
for all $\phi, \psi \in \mathcal{H}$. It is easy to see that adjoints
are unique (if they exist at all) and the set of all adjointable
operators $\mathfrak{B}(\mathcal{H})$ becomes a unital $^*$-algebra in
the obvious way. Similarly, one defines the adjointable maps
$\mathfrak{B}(\mathcal{H},\mathcal{H}')$ from $\mathcal{H}$ to some
other pre-Hilbert space $\mathcal{H}'$.  By $\Theta_{\phi, \psi}$ we
denote the \textdef{rank one operator} $\Theta_{\phi,\psi}\chi = \phi
\SP{\psi, \chi}$ where $\phi, \psi, \chi \in \mathcal{H}$.  The span
of all rank one operators, i.e.~the \textdef{finite rank
  operators}, is denoted by $\mathfrak{F}(\mathcal{H})$. Clearly,
$\mathfrak{F}(\mathcal{H})$ is a $^*$-ideal in
$\mathfrak{B}(\mathcal{H})$.  Analogously, one defines
$\mathfrak{F}(\mathcal{H}, \mathcal{H}')$.

%
%

\subsection{Pre-Hilbert modules and $^*$-representation theory}
\label{subsec:PreHilbertModules}

Let $\mathcal{A}$ be a $^*$-algebra and $\EA$ a right
$\mathcal{A}$-module. We shall always assume that all occurring modules
have an underlying compatible $\ring{C}$-module structure. Then an
\textdef{$\mathcal{A}$-valued inner product} on $\EA$ is a
$\ring{C}$-sesquilinear map $\SP{\cdot,\cdot}: \EA \times \EA
\longrightarrow \mathcal{A}$ such that $\SP{x, y \cdot a} = \SP{x,y}a$
and $\SP{x,y} = \SP{y,x}^*$. Sometimes we indicate the dependence on
the module and the algebra explicitly by $\SPEA{\cdot,\cdot}$. We call
$\SP{\cdot,\cdot}$ \textdef{non-degenerate} if $\SP{x,y} = 0$ for all
$y$ implies $x = 0$, in which case $\EA$ is called a \textdef{inner
  product $\mathcal{A}$-module}. We also make use of left modules with
inner products defined analogously, only linear to the left in the
first argument. For an inner product $\mathcal{A}$-module one has the
$^*$-algebra $\mathfrak{B}(\EA)$ of adjointable (and necessarily right
$\mathcal{A}$-linear) endomorphisms of $\EA$, whence $\mathcal{E}$
becomes a $(\mathfrak{B}(\EA), \mathcal{A})$-bimodule. Similarly, one
defines $\mathfrak{B}(\EA, \EpA)$ as well as the finite rank operators
$\mathfrak{F}(\EA)$ and $\mathfrak{F}(\EA, \EpA)$.

An $\mathcal{A}$-valued inner product $\SP{\cdot,\cdot}$ on $\EA$ is
called \textdef{completely positive} if the matrix $(\SP{x_i,x_j}) \in
M_n(\mathcal{A})$ is positive for all $x_1, \ldots, x_n \in \EA$ and
$n \in \mathbb{N}$. Here $M_n(\mathcal{A})$ is endowed with its
canonical $^*$-algebra structure. An inner product
$\mathcal{A}$-module with completely positive inner product is called
a \textdef{pre-Hilbert $\mathcal{A}$-module}.

Let $\mathcal{D}$ be another $^*$-algebra.  A
\textdef{$^*$-representation} of a $^*$-algebra $\mathcal{A}$ on an
inner product $\mathcal{D}$-module $\HD$ is a $^*$-homomorphism $\pi:
\mathcal{A} \longrightarrow \mathfrak{B}(\HD)$, generalizing thereby
the usual notion of a $^*$-representation on a pre-Hilbert space,
where $\mathcal{D} = \ring{C}$. An \textdef{intertwiner} $T$ between
two $^*$-representations $(\HD, \pi)$ and $(\HpD, \pi')$ of
$\mathcal{A}$ is an adjointable map $T \in \mathfrak{B}(\HD, \HpD)$
with $T\pi(a) = \pi'(a)T$ for all $a \in \mathcal{A}$.  It is easy to
see that the $^*$-representations of $\mathcal{A}$ on inner product
$\mathcal{D}$-modules form a category, denoted by
$\smod[\mathcal{D}](\mathcal{A})$, where morphisms are intertwiners.
The subcategory of \textdef{strongly non-degenerate}
$^*$-representations, i.e.~those with $\pi(\mathcal{A})\HD = \HD$, is
denoted by $\sMod[\mathcal{D}](\mathcal{A})$ and the subcategories of
(strongly non-degenerate) $^*$-representations on pre-Hilbert
$\mathcal{D}$-modules are denoted by $\rep[\mathcal{D}](\mathcal{A})$
and $\Rep[\mathcal{D}](\mathcal{A})$, respectively.
\begin{remark}
    \label{remark:UnitalOrNot}
    In the following we shall mainly be interested in \emph{unital}
    $^*$-algebras where we shall adopt the convention that
    $^*$-homomorphisms preserve units and units act as identities on
    modules. Thus for unital $^*$-algebras we have $\smod = \sMod$ and
    $\rep = \Rep$ by convention. In the non-unital case we still need
    some replacement for the units in order to obtain a reasonably
    good behavior. The right choices are \textdef{idempotent} and
    \textdef{non-degenerate} $^*$-algebras, see the discussion in
    \cite{bursztyn.waldmann:2003a:pre}.
\end{remark}

From \cite[Eq.~(4.7)]{bursztyn.waldmann:2003a:pre} one has a
functorial tensor product of inner product modules
\begin{equation}
    \label{eq:tensM}
    \tensM[\mathcal{B}]: 
    \smod[\mathcal{B}](\mathcal{C}) \times
    \smod[\mathcal{A}](\mathcal{B}) 
    \longrightarrow
    \smod[\mathcal{A}](\mathcal{C}),
\end{equation}
for three $^*$-algebras $\mathcal{A}$, $\mathcal{B}$, $\mathcal{C}$,
which is obtained as follows: For $\CFB \in
\smod[\mathcal{B}](\mathcal{C})$ and $\BEA \in
\smod[\mathcal{A}](\mathcal{B})$ one endows the algebraic tensor
product $\CFB \tensor[\mathcal{B}] \BEA$ with the $\mathcal{A}$-valued
inner product defined by
\begin{equation}
    \label{eq:SPtensorSP}
    \SPFEA{y\tensor[\mathcal{B}] x, y' \tensor[\mathcal{B}] x'}
    =
    \SPEA{x, \SPFB{y,y'} \cdot x'},
\end{equation}
where $x, x' \in \mathcal{E}$ and $y, y' \in \mathcal{F}$. Then one
divides by the (possibly non-empty) degeneracy space $(\CFB
\tensor[\mathcal{B}] \BEA)^\bot$ to obtain a non-degenerate
$\mathcal{A}$-valued inner product on the quotient $\CFB
\tensM[\mathcal{B}] \BEA = (\CFB \tensor[\mathcal{B}] \BEA) \big/
(\CFB \tensor[\mathcal{B}] \BEA)^\bot$, which is then a
$^*$-representation of $\mathcal{C}$. This construction generalizes
Rieffel's internal tensor product \cite{rieffel:1974a, rieffel:1974b},
which is a fundamental tool in $C^*$-algebra and Hilbert $C^*$-module
theory, see e.g.~\cite{raeburn.williams:1998a, lance:1995a}.  The
tensor product $\tensM$ is associative up to the usual canonical
isomorphism.  Moreover, if the inner products where both completely
positive then the resulting inner product \eqref{eq:SPtensorSP} is
completely positive again, see
\cite[Thm.~4.7]{bursztyn.waldmann:2003a:pre}.  Thus $\tensM$ restricts
to a functor
\begin{equation}
    \label{eq:tensMRep}
    \tensM[\mathcal{B}]: 
    \rep[\mathcal{B}](\mathcal{C}) \times 
    \rep[\mathcal{A}](\mathcal{B})
    \longrightarrow
    \rep[\mathcal{A}](\mathcal{C}).
\end{equation}

\begin{remark}[Complex conjugation of bimodules]
    \label{remark:LeftRightAndCC}
    Of course, we can also define $^*$-re\-pre\-sen\-tations from the
    right on inner product left modules. Then the analogous statements
    are still true. Furthermore, we can pass from one to the other by
    \textdef{complex conjugation} of the bimodule. For $\BEA \in
    \smod[\mathcal{A}](\mathcal{B})$ we define the $(\mathcal{A},
    \mathcal{B})$-bimodule $\AccEB$ by $\cc{\mathcal{E}} =
    \mathcal{E}$ as $\ring{R}$-module with $\ring{C}$-module structure
    given by $\alpha \cc{x} = \cc{\cc{\alpha} x}$ for $\alpha \in
    \ring{C}$ and $x \in \mathcal{E}$, where $\mathcal{E} \ni x
    \mapsto \cc{x} \in \cc{\mathcal{E}}$ denotes the identity map of
    the underlying $\ring{R}$-module. Then the $(\mathcal{A},
    \mathcal{B})$-bimodule structure is defined by
    \begin{equation}
        \label{eq:AccEAndccEB}
        a \cdot \cc{x} = \cc{x \cdot a^*}
        \quad
        \textrm{and}
        \quad
        \cc{x} \cdot b = \cc{b^* \cdot x}.
    \end{equation}
    The $\mathcal{A}$-left linear $\mathcal{A}$-valued inner product
    is defined by
    \begin{equation}
        \label{eq:ASPccE}
        \ASPccE{\cc{x}, \cc{y}} = \SPEA{x, y},
    \end{equation}
    which is clearly compatible with the right $\mathcal{B}$-module
    structure.  Then it is easily shown that $\ASPccE{\cdot,\cdot}$ is
    completely positive if and only if $\SPEA{\cdot,\cdot}$ is
    completely positive.
\end{remark}

%
%

\subsection{Hopf $^*$-algebras and $^*$-actions}
\label{subsec:HopfStarAlgebras}

Let $H$ be a Hopf algebra over $\ring{C}$ with comultiplication
$\Delta$, counit $\epsilon$ and antipode $S$. For $\Delta$ we use
Sweedler's notation, i.e.~$\Delta(g) = g_\sweedler{1} \otimes
g_\sweedler{2}$, etc. Now assume that $H$ is in addition a
$^*$-algebra. Then $H$ is called a \textdef{Hopf $^*$-algebra} if
$\Delta$ and $\epsilon$ are $^*$-homomorphisms and $S(S(g)^*)^* = g$,
see e.g.~\cite[Sect.~IV.8]{kassel:1995a}. In particular, $S$ is
invertible with inverse $S^{-1}(g) = S(g^*)^*$. The basic examples are
group algebras $\ring{C}[G]$ for a group $G$ and complexified
universal enveloping algebras $U_{\ring{C}}(\mathfrak{g}) =
U(\mathfrak{g}) \otimes_{\ring{R}} \ring{C}$ for Lie algebras
$\mathfrak{g}$ over $\ring{R}$, each endowed with the canonical Hopf
and $^*$-algebra structures. Both of them are \textdef{cocommutative},
i.e.~ $\Delta = \Delta^\op$, where $\Delta^\op(g) = g_\sweedler{2}
\otimes g_\sweedler{1}$ denotes the opposite comultiplication.

Let $\mathcal{A}$ be a $^*$-algebra over $\ring{C}$. A (left)
\textdef{$^*$-action} of $H$ on $\mathcal{A}$ is a (left) $H$-module
structure on $\mathcal{A}$, denoted by $(g, a) \mapsto g \acts a$ for
$g \in H$ and $a \in \mathcal{A}$, such that in addition
\begin{equation}
    \label{eq:gactsab}
    g \acts (ab) = (g_\sweedler{1} \acts a)(g_\sweedler{2} \acts b)
\end{equation}
\begin{equation}
    \label{eq:gactsaStar}
    (g \acts a)^* = S(g)^* \acts a^*.
\end{equation}
We shall use resulting formulas like $(g \acts a) b = g_{\sweedler{1}}
\acts(a S(g_\sweedler{2}) \acts b)$ for $a, b \in \mathcal{A}$ and $g
\in H$ as well as $a (g \acts b) = g_\sweedler{2}
\acts((S^{-1}(g_\sweedler{1}) \acts a) b)$ in the sequel, see
e.g.~\cite{majid:1995a, kassel:1995a, klimyk.schmuedgen:1997a} for
more details on the calculus with Hopf $^*$-algebras and $^*$-actions.

If $\mathcal{A}$ is unital we require $g \acts \Unit_{\mathcal{A}} =
\epsilon(g)\Unit_{\mathcal{A}}$ and we always assume that $\Unit_H
\acts a = a$. Recall that $a \in \mathcal{A}$ is called
\textdef{$H$-invariant} if $g \acts a = \epsilon(g) a$.

Similarly, one can define right $^*$-actions which we shall not need
in the sequel.  From now on, all $^*$-algebras are thought of being
equipped with a particular $^*$-action of $H$.
\begin{example}[$^*$-Actions of $H$]
    \label{example:StarActions}
    For later use we shall mention the following basic examples:
    \begin{compactenum}
    \item The \textdef{trivial $^*$-action} of $H$ on $\mathcal{A}$ is
        given by the counit $\epsilon$, i.e.
        \begin{equation}
            \label{eq:TrivialAction}
            g \acts a = \epsilon(g) a.
        \end{equation}
        It is easy to see that this is indeed a $^*$-action. The ring
        $\ring{C}$ and the matrices $M_n(\ring{C})$ are always assumed
        to carry the trivial $^*$-action.
    \item The \textdef{adjoint action}, see
        e.g.~\cite[Ex.~1.6.3]{majid:1995a}, of $H$ on itself is given
        by
        \begin{equation}
            \label{eq:AdjointAction}
            \Ad_g h = g_\sweedler{1} h S(g_\sweedler{2})
        \end{equation}
        and it turns out to be a $^*$-action as well.
    \item If $\mathcal{A}$ has a $^*$-action then the matrices
        $M_n(\mathcal{A})$ are endowed with a $^*$-action of $H$ as
        well by applying $g \in H$ componentwise.
    \end{compactenum}
\end{example}

Let $\mathcal{H}$ be a pre-Hilbert space over $\ring{C}$ and let
$\acts: H \longrightarrow \End_{\ring{C}}(\mathcal{H})$ be an action of
$H$ on $\mathcal{H}$. Then $\acts$ is called \textdef{unitary} if
\begin{equation}
    \label{eq:UnitaryOnPreHilbertSpace}
    \epsilon(g) \SP{x,y} 
    = \SP{S(g_\sweedler{1})^* \acts x, g_\sweedler{2} \acts y}
\end{equation}
for all $x,y \in \mathcal{H}$ and $g \in H$. Clearly, this gives
unitary representations of groups and (anti-) Hermitian
representations of real Lie algebras when applied to $\ring{C}[G]$ and
$U_{\ring{C}}(\mathfrak{g})$, respectively.

We generalize \eqref{eq:UnitaryOnPreHilbertSpace} as follows: Let
$\HD$ be a right $\mathcal{D}$-module, where $\mathcal{D}$ is an
auxiliary $^*$-algebra over $\ring{C}$ endowed with a $^*$-action of
$H$. Moreover, let $\HD$ be endowed with an $H$-module structure and
with a $\mathcal{D}$-valued inner product. Then the $H$-module
structure is called \textdef{compatible with the right
  $\mathcal{D}$-module structure} if
\begin{equation}
    \label{eq:gactsxd}
    g \acts (x \cdot d) 
    = (g_\sweedler{1} \acts x) \cdot (g_\sweedler{2} \acts d)
\end{equation}
and it is called \textdef{compatible with the inner product} if
\begin{equation}
    \label{eq:gactsSP}
    g \acts \SP{x,y} 
    = \SP{S(g_\sweedler{1})^* \acts x, g_\sweedler{2} \acts y}
\end{equation}
for $x,y \in \HD$, $d \in \mathcal{D}$ and $g \in H$.
\begin{lemma}
    \label{lemma:NiceAxiomForSP}
    The covariance condition \eqref{eq:gactsSP} is equivalent to the
    condition
    \begin{equation}
        \label{eq:NiceCovariance}
        \SP{x, g \acts y} 
        = g_\sweedler{2} \acts \SP{g_\sweedler{1}^* \acts x, y}
    \end{equation}
    for $x,y \in \HD$ and $g \in H$. If the $\mathcal{D}$-valued inner
    product is non-degenerate then \eqref{eq:gactsSP} implies
    \eqref{eq:gactsxd}.
\end{lemma}
\begin{proof}
    The equivalence of the two conditions \eqref{eq:NiceCovariance}
    and \eqref{eq:gactsSP} is a simple computation. Moreover, applying
    \eqref{eq:NiceCovariance} twice one obtains
    \[
    \SP{x, g\acts(y\cdot d)} 
    = \SP{x, (g_\sweedler{1} \acts y) \cdot (g_\sweedler{2} \acts d)},
    \]
    so the non-degeneracy of $\SP{\cdot,\cdot}$ implies
    \eqref{eq:gactsxd}. This generalizes
    \cite[Rem.~7.3]{raeburn.williams:1998a}.
\end{proof}
\begin{definition}
    \label{definition:CovariantDModule}
    An inner product right $\mathcal{D}$-module is called
    $H$-covariant if the inner product satisfies \eqref{eq:gactsSP}
    and hence also \eqref{eq:NiceCovariance} and \eqref{eq:gactsxd}.
\end{definition}

In case of a pre-Hilbert space \eqref{eq:NiceCovariance} simply means
that the action with $g \in H$ is \emph{adjointable} and we have
\begin{equation}
    \label{eq:CovariancePreHilbertSpace}
    \SP{x, g \acts y} = \SP{g^* \acts x, y}.
\end{equation}
\begin{remark}
    \label{remark:gactsNotAdjointable}
    Note that in general the operator $x \mapsto g \acts x$ is
    \emph{not} adjointable as endomorphism of $\HD$ since the action
    `outside' the inner product can be non-trivial.
\end{remark}
\begin{remark}
    \label{remark:DegeneracySpace}
    For a possibly degenerate inner product the condition
    \eqref{eq:NiceCovariance} immediately ensures $H \acts \HDbot
    \subseteq \HDbot$. Thus if the inner product $\SP{\cdot,\cdot}$ on
    $\HD$ is degenerate the $H$-action passes to the quotient $\HD
    \big/ \HDbot$ which then becomes an $H$-covariant inner product
    $\mathcal{D}$-module.
\end{remark}
\begin{proposition}
    \label{proposition:InducedAdjointAction}
    Let $\HD$ be an $H$-covariant inner product
    $\mathcal{D}$-module. Then
    \begin{equation}
        \label{eq:AdjointActionForBH}
        (g \acts A) x 
        = g_\sweedler{1} \acts (A S(g_\sweedler{2}) \acts x),
    \end{equation}
    for $A \in \mathfrak{B}(\HD)$ and $x \in \HD$, defines a
    $^*$-action of $H$ on $\mathfrak{B}(\HD)$ uniquely determined by
    the property
    \begin{equation}
        \label{eq:BDHHDBimoduleCovariant}
        g \acts (Ax) 
        = (g_\sweedler{1} \acts A) (g_\sweedler{2} \acts x).
    \end{equation}
    Moreover, we have for the rank one operators
        \begin{equation}
        \label{eq:HactsOnFiniteRank}
        g \acts \Theta_{x,y} 
        = 
        \Theta_{g_\sweedler{1} \acts x, S(g_\sweedler{2})^* \acts y},
    \end{equation}
    whence we have $H \acts \mathfrak{F}(\HD) \subseteq
    \mathfrak{F}(\HD)$.
\end{proposition}
\begin{proof}
    Using the same kind of calculations as for the adjoint action of
    $H$ on itself one shows that \eqref{eq:AdjointActionForBH} defines
    an action of $H$ on \emph{all} endomorphisms
    $\End_{\ring{C}}(\HD)$, which is uniquely determined by the
    property \eqref{eq:BDHHDBimoduleCovariant}. This part is fairly
    standard and well-known. It remains to show that for $A \in
    \mathfrak{B}(\HD)$ the result $g \acts A$ is again adjointable
    with adjoint given by $S(g)^* \acts A^*$. Note that this is
    non-trivial according to Remark~\ref{remark:gactsNotAdjointable}.
    One computes
    \[
    \SP{x, g_\sweedler{1} \acts (A S(g_\sweedler{2}) \acts y)}
    = g_\sweedler{2} \acts S(g_\sweedler{3})_\sweedler{2} \acts 
    \SP{
      S(g_\sweedler{3})_\sweedler{1}^* \acts 
      (A^* g_\sweedler{1}^* \acts x), y}
    = \SP{S(g_\sweedler{2})^* \acts (A^* g_\sweedler{1}^* \acts x), y},
    \]
    using twice \eqref{eq:NiceCovariance} and the fact that $A$ is
    adjointable as well as $S \otimes S \circ \Delta^\op = \Delta
    \circ S$. This shows that $g\acts A$ is indeed adjointable with
    adjoint given by
    \[
    (g \acts A)^* x 
    = S(g_\sweedler{2})^* \acts A^* g_\sweedler{1} \acts x
    = (S(g)^* \acts A^*) x,
    \]
    whence the action \eqref{eq:AdjointActionForBH} is a $^*$-action.
    The fact that $g \acts A$ is again right $\mathcal{D}$-linear
    follows from the existence of an adjoint.  The statement about the
    rank one operator follows analogously.
\end{proof}

For obvious reasons we call the action on $\mathfrak{B}(\HD)$ the
\textdef{adjoint action} induced by the action on $\HD$.

Let us now mention one of our motivating examples from geometry:
\begin{example}[Lie algebra action on a manifold]
    \label{example:LieAlgebraAction}
    Let $\lie{g}$ be a real finite-dimensional Lie algebra and $M$ a
    smooth manifold and let $H = U_{\mathbb{C}}(\lie{g})$ be the
    complexified universal enveloping algebra of $\lie{g}$, viewed as
    Hopf $^*$-algebra. Then a $^*$-action of $H$ on the complex-valued
    smooth functions $C^\infty(M)$ is equivalent to a Lie algebra
    action of $\lie{g}$ on $M$, i.e.  an Lie algebra homomorphism
    $\lie{g} \longrightarrow \Gamma^\infty(TM)$. This follows from the
    fact that the condition~\eqref{eq:gactsab} implies that $\xi \in
    \lie{g}$ acts as derivation on $C^\infty(M)$ together with the
    fact that any derivation of $C^\infty(M)$ is given by a vector
    field.
\end{example}
\begin{example}[Lie group action on a manifold]
    \label{example:LieGroupAction}
    If $G$ is a Lie group and $\Phi: G \times M \longrightarrow M$ a
    smooth Lie group action on a manifold $M$, then the pull-backs
    $\Phi_g^*$ act on $C^\infty(M)$ by $^*$-automorphisms. This yields
    a $^*$-action of the Hopf $^*$-algebra $\mathbb{C}[G]$. Note
    however, that in our (algebraic) definition of the group algebra
    $H = \mathbb{C}[G]$ no topological information about $G$ is
    contained.  Thus not every $^*$-action of $H$ on $C^\infty(M)$
    comes from a \emph{smooth} action of $G$ on $M$. Here one has to
    impose additional conditions which go beyond our purely algebraic
    treatment.
\end{example}
\begin{remark}
    \label{remark:StarProducts}
    The above two examples provide the framework for symmetries in
    classical mechanics. In deformation quantization such symmetries
    are encoded in the notions of invariant star products, see e.g.
    \cite{arnal.cortet.molin.pinczon:1983a, gutt.rawnsley:2003a,
      mueller-bahns.neumaier:2004a} and references therein. Here we
    have to pass from $\mathbb{R}$ and $\mathbb{C}$ to the ordered
    ring $\mathbb{R}[[\lambda]]$ and $\mathbb{C}[[\lambda]]$.
\end{remark}

%
%

\subsection{$H$-covariant representation theory}
\label{subsec:CovRepTheory}

Let $\mathcal{A}$ be a $^*$-algebra and let $\mathcal{D}$ be an
auxiliary $^*$-algebra as above, both endowed with a fixed
$^*$-action of $H$. If $\HD$ is an $H$-covariant inner product right
$\mathcal{D}$-module then a $^*$-representation $\pi$ of $\mathcal{A}$
on $\HD$ is called \textdef{$H$-covariant} if
\begin{equation}
    \label{eq:RepCovariant}
    \pi(g \acts a) x 
    = g_\sweedler{1} \acts (\pi(a) S(g_\sweedler{2}) \acts x)
\end{equation}
holds for all $a \in \mathcal{A}$ and $x \in \HD$. Again, applied to
pre-Hilbert spaces and for group algebras or complexified universal
enveloping algebras one recovers the usual notion of a covariant
$^*$-representation. Another way to view \eqref{eq:RepCovariant} is
that the map $\pi: \mathcal{A} \longrightarrow \mathfrak{B}(\HD)$ is
$H$-equivariant with respect to the adjoint action on
$\mathfrak{B}(\HD)$ induced by the action on $\HD$.

An intertwiner $T: \HD \longrightarrow \HpD$ between two $H$-covariant
$^*$-representations $(\HD, \pi)$ and $(\HpD, \pi')$ of $\mathcal{A}$
is called \textdef{$H$-covariant} if $T$ also intertwines the
$H$-module structure, i.e.
\begin{equation}
    \label{eq:TIntertwinergacts}
    T(g \acts x) = g \acts T(x).
\end{equation}
Then one obtains the \textdef{category of $H$-covariant
  $^*$-representations} of $\mathcal{A}$ on $H$-covariant inner
product right $\mathcal{D}$-modules, denoted by $\smod[\mathcal{D},
H](\mathcal{A})$, where $H$-covariant intertwiners are used as
morphisms. Analogously, one defines the sub-categories
$\sMod[\mathcal{D}, H](\mathcal{A})$ as well as $\rep[\mathcal{D},
H](\mathcal{A})$ and $\Rep[\mathcal{D}, H](\mathcal{A})$.
\begin{remark}
    \label{remark:ccExtendsHcov}
    Also in this framework we can pass from left to right
    $^*$-representations. For a left $\mathcal{B}$-representation on a
    $H$-covariant inner product right $\mathcal{A}$-module $\BEA$ we
    define the $H$-action $\ccacts$ on $\AccEB$ by
    \begin{equation}
        \label{eq:ccacts}
        g \ccacts \cc{x} = \cc{S(g)^* \acts x},
    \end{equation}
    which can be shown to be an $H$-action compatible with the complex
    conjugated bimodule structure as well as with the complex
    conjugated inner product $\ASPccE{\cdot,\cdot}$. This is a
    straightforward computation.  Moreover, $\cc{\cc{\mathcal{E}}} =
    \mathcal{E}$, including all its structures.
\end{remark}

The prototype of an $H$-covariant $^*$-representation is obtained by
the GNS representation with respect to an $H$-invariant positive linear
functional on $\mathcal{A}$:
\begin{example}[$H$-invariant GNS construction]
    \label{example:InvariantGNS}
    The usual GNS construction of a $^*$-re\-pre\-sen\-tation out of a
    positive linear functional can be generalized immediately to the
    $H$-covariant framework. Let $\omega: \mathcal{A} \longrightarrow
    \ring{C}$ be a \textdef{$H$-invariant} positive linear functional,
    i.e.~we have $\omega(a^*a) \ge 0$ and $\omega( g \acts a) =
    \epsilon(g) \omega(a)$ for all $a \in \mathcal{A}$ and $g \in H$.
    Then we consider the inner product
    \begin{equation}
        \label{eq:GNSInnerProductOnA}
        \SP{a, b}_\omega = \omega(a^*b)
    \end{equation}
    on $\mathcal{A}$, viewed as $(\mathcal{A}, \ring{C})$-bimodule.
    We have
    \begin{equation}
        \label{eq:SPGNSCovariant}
        \SP{g^* \acts a, b}_\omega = \SP{a, g \acts b}_\omega
    \end{equation}
    by a straightforward computation using the invariance of $\omega$,
    whence $\SP{\cdot,\cdot}_\omega$ is compatible with the
    $H$-action. Thus we can apply Remark~\ref{remark:DegeneracySpace}
    and divide by the (possibly non-empty) degeneracy space of
    $\SP{\cdot,\cdot}_\omega$ to obtain a pre-Hilbert module
    $\mathcal{H}_\omega = \mathcal{A} \big/ \mathcal{A}^\bot$ over
    $\ring{C}$, i.e.~a pre-Hilbert space. Note that
    \begin{equation}
        \label{eq:DegSpaceIsGelfandIdeal}
        \mathcal{A}^\bot = \mathcal{J}_\omega 
        = \{a \in \mathcal{A} \; | \; \omega(a^*a) = 0 \}
    \end{equation}
    is just the Gel'fand ideal of $\omega$. Thus we recover the usual
    GNS representation $\pi_\omega$ of $\mathcal{A}$ on
    $\mathcal{H}_\omega$ together with an $H$-action making the GNS
    representation $H$-covariant.

    If in addition $\mathcal{A}$ is unital then the class of
    $\Unit_{\mathcal{A}}$ in $\mathcal{H}_\omega$ is a
    \textdef{cyclic} $H$-invariant vector, the \textdef{vacuum
      vector}. Every other $H$-covariant $^*$-representation
    $(\mathcal{H}, \pi)$ of $\mathcal{A}$ with $H$-invariant cyclic
    vector $\Omega \in \mathcal{H}$ such that $\omega(a) = \SP{\Omega,
      \pi(a)\Omega}$ is unitarily equivalent to the GNS representation
    via the usual $H$-covariant intertwiner. Needless to say, this
    example is of fundamental importance for the understanding of the
    $H$-covariant $^*$-representation theory of $\mathcal{A}$.
\end{example}

%
%

\subsection{The lattice of $(\mathcal{D}, H)$-closed $^*$-ideals}
\label{subsec:LatticeDef}

For a $C^*$-algebra a $^*$-ideal is topologically closed if and only
if it is the kernel of a $^*$-representation, see e.g.~\cite[Chap.~I,
Thm.~1.3.10]{landsman:1998a}. This fact was the motivation to define a
$^*$-ideal in a $^*$-algebra to be \textdef{closed} if it is the kernel
of a $^*$-representation of $\mathcal{A}$ on a pre-Hilbert space, see
\cite{bursztyn.waldmann:2001b}. We extend this definition now in two
directions, allowing for $^*$-representations on pre-Hilbert
$\mathcal{D}$-modules instead of pre-Hilbert spaces and incorporating
$H$-covariance.

For reasons which become clear in
Section~\ref{subsec:IdealMoritaInvariant} we have to restrict the
auxiliary $^*$-algebras $\mathcal{D}$. The problem is that for a
pre-Hilbert $\mathcal{D}$-module $\HD$ the inner product
$\SPD{\cdot,\cdot}$ is completely positive and non-degenerate
\emph{but} there may be elements $\phi \in \HD$ with $\SPD{\phi, \phi}
= 0$ and $\phi \ne 0$. The Grassmann algebra $\Lambda (\ring{C}^n)$
provides a simple example, see
\cite[Ex.~3.5]{bursztyn.waldmann:2003a:pre}. In order to avoid this we
state the following definition: We call $\mathcal{D}$
\textdef{admissible} if on any pre-Hilbert $\mathcal{D}$-module $\HD$
the inner product $\SPD{\cdot,\cdot}$ is in addition \textdef{positive
  definite}, i.e. $\SPD{\phi, \phi} = 0$ implies $\phi = 0$.  This is
the case if e.g.~$\mathcal{D}$ has sufficiently many positive linear
functionals in the sense that for any Hermitian element $d = d^* \ne
0$ we find a positive linear functional $\omega: \mathcal{D}
\longrightarrow \ring{C}$ with $\omega(d) \ne 0$ and if $d + d = 0$
implies $d = 0$ for all $d \in \mathcal{D}$, see
\cite[Ex.~3.6]{bursztyn.waldmann:2003a:pre}. Then we can state the
following definition:
\begin{definition}
    \label{definition:DHclosedIdeal}
    Let $\mathcal{D}$ be admissible. Then $\mathcal{J} \subseteq
    \mathcal{A}$ is called a $(\mathcal{D}, H)$-closed ideal if
    $\mathcal{J} = \ker \pi$ for some $^*$-representation $(\HD, \pi)
    \in \rep[\mathcal{D}, H](\mathcal{A})$. The set of all
    $(\mathcal{D}, H)$-closed ideals is denoted by
    $\Lattice[\mathcal{D}, H](\mathcal{A})$.
\end{definition}

Clearly, if $\mathcal{D} = \ring{C}$ the non-covariant version gives
the lattice of closed ideals $\Lattice(\mathcal{A})$ as in
\cite[Sect.~4]{bursztyn.waldmann:2001b}.  We collect a few first
properties of $\Lattice[\mathcal{D}, H](\mathcal{A})$ which can be
obtained completely analogously as for $\Lattice(\mathcal{A})$.
\begin{lemma}
    \label{lemma:Lattice}
    Let $\mathcal{D}$ be admissible and let $\mathcal{A}$ be
    idempotent.
    \begin{compactenum}
    \item If $\mathcal{J} \in \Lattice[\mathcal{D}, H](\mathcal{A})$
        then $\mathcal{J}$ is an $H$-invariant $^*$-ideal of
        $\mathcal{A}$.
    \item If $\mathcal{J} \in \Lattice[\mathcal{D}, H](\mathcal{A})$
        then there exists a strongly non-degenerate
        $^*$-representation $(\HD, \pi) \in \Rep[\mathcal{D},
        H](\mathcal{A})$ with $\mathcal{J} = \ker \pi$.
    \item Let $I$ be a set and $\mathcal{J}_\alpha \in
        \Lattice[\mathcal{D}, H](\mathcal{A})$ for $\alpha \in I$.
        Then $\bigcap_{\alpha \in I} \mathcal{J}_\alpha \in
        \Lattice[\mathcal{D}, H](\mathcal{A})$.
    \item For an arbitrary subset $\mathcal{J} \subseteq \mathcal{A}$
        let
        \begin{equation}
            \label{eq:closure}
            \mathcal{J}^\cl = 
            \bigcap_{\mathcal{J} \subseteq \mathcal{I}, 
              \mathcal{I} \in \Lattice[\mathcal{D}, H](\mathcal{A})} 
            \mathcal{I}.
        \end{equation}
        Then $\mathcal{J}^\cl \in \Lattice[\mathcal{D},
        H](\mathcal{A})$ is the smallest $(\mathcal{D}, H)$-closed
        ideal containing $\mathcal{J}$ and $(\mathcal{J}^\cl)^\cl =
        \mathcal{J}^\cl$.
    \item The operations $\mathcal{I} \wedge \mathcal{J} = \mathcal{I}
        \cap \mathcal{J}$ and $\mathcal{I} \vee \mathcal{J} =
        (\mathcal{I} \cup \mathcal{J})^\cl$ define the structure of a
        lattice on $\Lattice[\mathcal{D}, H](\mathcal{A})$ such that
        $\mathcal{I} \le \mathcal{J}$ if and only if $\mathcal{I}
        \subseteq \mathcal{J}$.
    \end{compactenum}
\end{lemma}
\begin{proof}
    The proof is completely analogous to the corresponding ones in
    \cite{bursztyn.waldmann:2001b} since inner products are always
    positive definite.
\end{proof}

We call $\Lattice[\mathcal{D}, H](\mathcal{A})$ the \textdef{lattice
  of $(\mathcal{D}, H)$-closed ideals} of the $^*$-algebra
$\mathcal{A}$. Note that only for the second part of the lemma one
needs that $\mathcal{D}$ is admissible.

%
%

\section{$H$-Covariant strong Morita equivalence}
\label{sec:CovSME}

In this section we adapt the tensor product $\tensM$ to the
$H$-covariant situation and obtain this way an $H$-covariant version of
Rieffel induction. This tensor product will allow a definition of
$H$-covariant strong Morita equivalence which implies the usual strong
Morita equivalence, see e.g.~\cite{raeburn.williams:1998a} for the
analogous construction for $G$-covariant strong Morita equivalence of
$C^*$-algebras.

%
%

\subsection{$H$-covariant tensor products}
\label{subsec:TensorProducts}

First we show how the functor $\tensM$ from \eqref{eq:tensMRep}
restricts to $H$-covariant $^*$-representations. For a right
$\mathcal{B}$-module $\FB$ with $H$-covariant $\mathcal{B}$-valued
inner product $\SPFB{\cdot,\cdot}$ and an $H$-covariant $(\mathcal{B},
\mathcal{A})$-bimodule $\BEA$ with $H$-covariant $\mathcal{A}$-valued
inner product $\SPEA{\cdot,\cdot}$ compatible with the
$\mathcal{B}$-action in the sense that $\SPEA{b \cdot x, y} = \SPEA{x,
  b^* \cdot y}$ we have the inner product $\SPFEA{\cdot,\cdot}$ on
$\FB \tensor[\mathcal{B}] \BEA$ given by \eqref{eq:SPtensorSP}.  On
the tensor product we also have canonically an action of $H$ defined
as usual by
\begin{equation}
    \label{eq:ActionTensorProduct}
    g \acts (x \tensor[\mathcal{B}] y)
    =
    g_\sweedler{1} \acts x 
    \tensor[\mathcal{B}] g_\sweedler{2} \acts y,
\end{equation}
which is indeed easily shown to be well-defined over
$\tensor[\mathcal{B}]$ and an action of $H$.
\begin{lemma}
    \label{lemma:InternalTensor}
    The canonical $H$-action on $\FB \tensor[\mathcal{B}] \BEA$ given
    by the tensor product of the action on $\FB$ and $\BEA$ makes
    $\SPFEA{\cdot,\cdot}$ an $H$-covariant inner product. Moreover, the
    $H$-action passes to the quotient $\FB \tensor[\mathcal{B}] \BEA
    \big/ (\FB \tensor[\mathcal{B}] \BEA)^\bot$ which becomes a
    $H$-covariant inner product $\mathcal{A}$-module.
\end{lemma}
\begin{proof}
    Let $x, x' \in \mathcal{F}$ and $y, y' \in \mathcal{E}$ as well as
    $g \in H$. From the $H$-covariance of the inner products
    $\SPEA{\cdot,\cdot}$ and $\SPFB{\cdot,\cdot}$ we conclude that
    \[
    \begin{split}
        g \acts \SPFEA{x \otimes y, x' \otimes y'}
        &= g \acts \SPEA{y, \SPFB{x, x'} \cdot y'} \\
        &= \SPEA{S(g_\sweedler{1})^* \acts y, 
          g_\sweedler{2} \acts \left(\SPFB{x, x'} \cdot y'\right)} \\
        &= \SPEA{S(g_\sweedler{1})^* \acts y, 
          \SPFB{S(g_\sweedler{2})^* \acts x, g_\sweedler{3} \acts x'}
          \cdot (g_\sweedler{4} \acts y')} \\
        &= \SPFEA{
          S(g_\sweedler{2})^* \acts x \otimes
          S(g_\sweedler{1})^* \acts y, 
          g_\sweedler{3} \acts x' \otimes
          g_\sweedler{4} \acts y'} \\
        &= \SPFEA{
          S(g_\sweedler{1})^* \acts (x \otimes y), 
          g_\sweedler{2} \acts(x' \otimes y')},
    \end{split}
    \]
    using $S \otimes S \circ \Delta^\op = \Delta \circ S$ in the last
    step.  This proves the compatibility of the inner product with the
    $H$-action.  The passage to the quotient follows immediately from
    Remark~\ref{remark:DegeneracySpace}.
\end{proof}

Thus we can define the \textdef{$H$-covariant internal tensor product}
of $\FB$ and $\BEA$ to be the right $\mathcal{A}$-module $\FB
\tensM[\mathcal{B}] \BEA = \FB \tensor[\mathcal{B}] \BEA \big/ (\FB
\tensor[\mathcal{B}] \BEA)^\bot$ endowed with its $H$-action and its
$H$-covariant $\mathcal{A}$-valued inner product. If $\FB$ carries in
addition an $H$-covariant $^*$-representation of some $^*$-algebra
$\mathcal{C}$ then the induced $^*$-representation of $\mathcal{C}$ on
$\FB \tensM[\mathcal{B}] \BEA$ is again $H$-covariant.  The
functoriality of the tensor product of $H$-actions, i.e.~tensor
products of intertwiners give intertwiners, together with the
functoriality of the internal tensor product of inner products as in
\cite[Lem.~4.16]{bursztyn.waldmann:2003a:pre} finally gives a functor
\begin{equation}
    \label{eq:CovariantInternalTensorProduct}
    \tensM[\mathcal{B}]: 
    \smod[\mathcal{B}, H](\mathcal{C}) \times
    \smod[\mathcal{A}, H](\mathcal{B}) 
    \longrightarrow
    \smod[\mathcal{A}, H](\mathcal{C}).
\end{equation}
It is easy to see that the usual associativity of the tensor product
gives associativity of $\tensM$ up to the usual canonical isomorphism,
i.e.
\begin{equation}
    \label{eq:GFEAsso}
    (\GC \tensM[\mathcal{C}] \CFB) \tensM[\mathcal{B}] \BEA
    \cong
    \GC \tensM[\mathcal{C}] (\CFB \tensM[\mathcal{B}] \BEA),
\end{equation}
see \cite[Lem.~4.5]{bursztyn.waldmann:2003a:pre} for the non-covariant
case. Since $\tensM$ is compatible with complete positivity of inner
products \cite[Thm.~4.7]{bursztyn.waldmann:2003a:pre} the functor
\eqref{eq:CovariantInternalTensorProduct} restricts to a functor
\begin{equation}
    \label{eq:repCovInternalTensor}
    \tensM[\mathcal{B}]:
    \rep[\mathcal{B}, H](\mathcal{C}) \times
    \rep[\mathcal{A}, H](\mathcal{B}) 
    \longrightarrow
    \rep[\mathcal{A}, H](\mathcal{C}).
\end{equation}
Fixing one of the two arguments of $\tensM$ we get the $H$-covariant
versions of Rieffel induction and the change of the base ring as in
\cite[Ex.~4.9 \& 4.10]{bursztyn.waldmann:2003a:pre}. The
\textdef{$H$-covariant Rieffel induction} with some $\BEA \in
\rep[\mathcal{A}, H](\mathcal{B})$ is denoted by
\begin{equation}
    \label{eq:RieffelInduction}
    \mathsf{R}_{\mathcal{E}} = \BEA \tensM \cdot :
    \rep[\mathcal{D}, H](\mathcal{A}) 
    \longrightarrow
    \rep[\mathcal{D}, H](\mathcal{B})
\end{equation}
and the \textdef{$H$-covariant change of the base ring} with some
$\DGDp \in \rep[\mathcal{D}', H](\mathcal{D})$ is denoted by
\begin{equation}
    \label{eq:ChangeOfBase}
    \mathsf{S}_{\mathcal{G}} = \cdot \tensM \DGDp :
    \rep[\mathcal{D}, H](\mathcal{A})
    \longrightarrow
    \rep[\mathcal{D}', H](\mathcal{A}).
\end{equation}
The functors $\mathsf{R}_{\mathcal{E}}$ and $\mathsf{S}_{\mathcal{G}}$
commute up to the usual natural transformation induced by
\eqref{eq:GFEAsso}.

%
%

\subsection{$H$-covariant strong Morita equivalence}
\label{subsec:HcovMoritaEquivalence}

We are now able to adapt the notions of Ara's $^*$-Morita equivalence
\cite{ara:1999a} and strong Morita equivalence
\cite{bursztyn.waldmann:2003a:pre} to the $H$-covariant framework.
Recall that an inner product $\SPEA{\cdot,\cdot}$ is \textdef{full} if
the $\ring{C}$-span of the elements $\SPEA{x,y} \in \mathcal{A}$ gives
the whole $^*$-algebra $\mathcal{A}$ and analogously for
$\BSPE{\cdot,\cdot}$.
\begin{definition}
    \label{definition:HcovariantEquivBimodule}
    A $(\mathcal{B}, \mathcal{A})$-bimodule $\BEA$ with inner products
    $\BSPE{\cdot,\cdot}$ and $\SPEA{\cdot,\cdot}$ is called a
    $H$-covariant $^*$-equivalence bimodule if it is a
    $^*$-equivalence bimodule in the sense of
    \cite[Def.~5.1]{ara:1999a} together with an action of $H$ such
    that
    \begin{equation}
        \label{eq:HcovSPEA}
        g \acts \SPEA{x,y} 
        = \SPEA{S(g_\sweedler{1})^* \acts x,
          g_\sweedler{2} \acts y}
    \end{equation}
    and
    \begin{equation}
        \label{eq:HcovBSPE}
        g \acts \BSPE{x, y} =
        \BSPE{g_\sweedler{1} \acts x, S(g_\sweedler{2})^* \acts y}
    \end{equation}
    for all $x,y \in \BEA$ and $g \in H$. It is called an $H$-covariant
    strong equivalence bimodule if in addition the underlying
    $^*$-equivalence bimodule is a strong equivalence bimodule in the
    sense of \cite[Def.~5.1]{bursztyn.waldmann:2003a:pre}, i.e.~the
    inner products are both completely positive.
\end{definition}

Recall that $\BEA$ is a $^*$-equivalence bimodule in Ara's sense if
both inner products are non-degenerate, full and satisfy the
compatibility conditions
\begin{equation}
    \label{eq:SPCompatible}
    \SPEA{b \cdot x, y} = \SPEA{x, b^* \cdot y},
    \quad
    \BSPE{x \cdot a, y} = \BSPE{x, y \cdot a^*}
    \quad
    \textrm{and}
    \quad
    \BSPE{x, y} \cdot z = x \cdot \SPEA{y,z}
\end{equation}
and if $\mathcal{B} \cdot \mathcal{E} = \mathcal{E} = \mathcal{E}
\cdot \mathcal{A}$, see \cite[Sect.~5.1]{ara:1999a} for details.  The
compatibility \eqref{eq:SPCompatible} can also be interpreted as
$\BSPE{x, y} \cdot z = \Theta_{x, y}(z)$.  Moreover,
\eqref{eq:HcovSPEA} and \eqref{eq:HcovBSPE} imply by
Lemma~\ref{lemma:NiceAxiomForSP} and the non-degeneracy of the inner
products that the bimodule structure is compatible with the
$H$-action, i.e.
\begin{equation}
    \label{eq:EquivBimComHAction}
    g \acts( b \cdot x) 
    = (g_\sweedler{1} \acts b) \cdot (g_\sweedler{2} \acts x)
    \quad
    \textrm{and}
    \quad
    g \acts (x \cdot a)
    = (g_\sweedler{1} \acts x) \cdot (g_\sweedler{2} \acts a).
\end{equation}
Thus $(\BEA, \SPEA{\cdot,\cdot})$ is an $H$-covariant
$^*$-representation of $\mathcal{B}$ on an $H$-covariant inner product
right $\mathcal{A}$-module and analogously for exchanged roles of
$\mathcal{A}$ and $\mathcal{B}$.
\begin{definition}
    \label{definition:HcovMorita}
    Two $^*$-algebras $\mathcal{A}$ and $\mathcal{B}$ with $^*$-action
    of $H$ are called $H$-covariantly $^*$-Morita equivalent (resp.
    $H$-covariantly strongly Morita equivalent) if there exists a
    $H$-covariant $^*$-Morita (resp. strong Morita) equivalence
    bimodule for them.
\end{definition}

Clearly, $H$-covariant $^*$- or strong Morita equivalence implies
$^*$- or strong Morita equivalence, respectively, and $H$-covariant
strong Morita equivalence implies $H$-covariant $^*$-Morita
equivalence.  Moreover, as expected, $H$-covariant $^*$- as well as
strong Morita equivalence turn out to be equivalence relations when
applied to non-degenerate and idempotent $^*$-algebras. This
restriction is necessary according to
\cite{bursztyn.waldmann:2003a:pre, ara:1999a}.
\begin{theorem}
    \label{theorem:MoritaEquivalenceRelation}
    Within the class of idempotent and non-degenerate $^*$-algebras
    with $^*$-actions of $H$, $H$-covariant $^*$- or strong Morita
    equivalence are both equivalence relations. Moreover,
    $H$-equivariantly $^*$-isomorphic $^*$-algebras are
    $H$-covariantly strongly Morita equivalent and hence also
    $H$-covariantly $^*$-Morita equivalent.
\end{theorem}
\begin{proof}
    We already know that the underlying $^*$- or strong Morita
    equivalence is an equivalence relation where the bimodule $\AAA$
    with the canonical inner products
    \begin{equation}
        \label{eq:SPAASP}
        \SPA{a,b} = a^*b
        \quad
        \textrm{and}
        \quad
        \ASP{a,b} = a b^*
    \end{equation}
    gives reflexivity. The complex conjugate bimodule $\AccEB$, see
    Remark~\ref{remark:LeftRightAndCC}, gives symmetry. Finally the
    internal tensor product $\tensM$ gives transitivity, see
    \cite{ara:1999a, bursztyn.waldmann:2003a:pre}.  Thus it remains to
    show that the three constructions are compatible with the
    $H$-covariance.  Clearly, the bimodule structure on $\AAA$ is
    $H$-covariant and we have
    \[
    g \acts \SPA{a,b} 
    = g \acts (a^*b)
    = (g_\sweedler{1} \acts a^*) (g_\sweedler{2} \acts b)
    = \left(S(g_\sweedler{1})^* \acts a\right)^* (g_\sweedler{2} \acts b)
    = \SPA{S(g_\sweedler{1})^* \acts a, g_\sweedler{2} \acts b},
    \]
    and similarly for $\ASP{\cdot,\cdot}$. Thus the inner products on
    $\AAA$ are $H$-covariant which proves reflexivity.  
    On $\AccEB$ we have already constructed the candidate for the
    $H$-action in Remark~\ref{remark:ccExtendsHcov}. A simple
    computation shows that $\ccacts$ is compatible with the
    $\mathcal{B}$-valued inner product as well. This follows
    immediately from the compatibility \eqref{eq:SPCompatible} or from
    a straightforward direct computation.
    Finally, transitivity follows from our considerations in
    Lemma~\ref{lemma:InternalTensor} where we have already shown that
    $\tensM$ is compatible with $H$-actions. Note however, that now we
    have to check the compatibility with two inner products, which can
    be done in a completely analogous way as for one. Thus
    $H$-covariant $^*$- or strong Morita equivalence is an equivalence
    relation.  Let us finally consider a $^*$-isomorphism $\Phi:
    \mathcal{A} \longrightarrow \mathcal{B}$ such that $\Phi$ is
    $H$-equivariant, i.e.~$\Phi(g \acts a) = g \acts \Phi(a)$ for all
    $a \in \mathcal{A}$ and $g \in H$. Then we claim that $\BBPhiA$ is
    an $H$-covariant strong Morita equivalence bimodule, where the
    right $\mathcal{A}$-module structure on $\mathcal{B}$ is defined
    by $b \dotPhi a = b\Phi(a)$ and the $\mathcal{A}$-valued inner
    product is $\SPBPhiA{b_1, b_2} = \Phi^{-1}(b_1^* b_2)$. Again, we
    only have to check the $H$-covariance which is a simple
    computation.
\end{proof}
\begin{remark}
    \label{remark:NonDegIdempotent}
    From now on we shall always assume that the $^*$-algebras in
    question are idempotent and non-degenerate since otherwise Morita
    theory becomes somewhat pathological as Morita equivalence no
    longer defines a reflexive relation.
\end{remark}

As we shall need the tensor product of equivalence bimodules
throughout this article, we introduce a new notation: For two
equivalence bimodules (either $H$-covariant $^*$- or strong Morita
equivalence) $\CFB$ and $\BEA$ we denote their internal tensor product
by $\CFB \tensB[\mathcal{B}] \BEA$ to stress that now \emph{two} inner
products are involved. From
\cite[Lem.~5.7]{bursztyn.waldmann:2003a:pre} we know that the
degeneracy spaces of the two inner products on the algebraic tensor
product $\CFB \tensor[\mathcal{B}] \BEA$ coincide if each of the
bimodules is an equivalence bimodule. This is a simple consequence of
\eqref{eq:SPCompatible}. Thus dividing by the degeneracy space is
non-ambiguous. It is clear that $\tensB$ enjoys analogous
functoriality properties as $\tensM$.

Let us now discuss some basic consequences of $H$-covariant $^*$- or
strong Morita equivalence:
\begin{proposition}
    \label{proposition:ActionOnBisTheOneOfF}
    Let $\mathcal{A}$, $\mathcal{B}$ be non-degenerate and idempotent
    $^*$-algebras over $\ring{C}$ and let $\BEA$ be an $H$-covariant
    $^*$-Morita equivalence bimodule. Then $\mathcal{B}$ is
    canonically $^*$-isomorphic to $\mathfrak{F}(\EA)$ via the action
    map
    \begin{equation}
        \label{eq:BisotoFEA}
        \mathcal{B} \ni b 
        \mapsto 
        (x \mapsto b \cdot x) \in \mathfrak{F}(\EA)
    \end{equation}
    and the $^*$-action of $H$ on $\mathcal{B}$ corresponds under
    \eqref{eq:BisotoFEA} to the adjoint action on $\mathfrak{F}(\EA)$
    induced by the action on $\mathcal{E}$.  In particular, if
    $\mathcal{B}$ is unital then $\mathcal{B} \cong \mathfrak{F}(\EA)
    = \mathfrak{B}(\EA)$.  Conversely, if $\EA$ is a right
    $\mathcal{A}$-module with $H$-action and compatible full inner
    product $\SPEA{\cdot,\cdot}$ such that $\EA = \EA \cdot
    \mathcal{A}$ then the $^*$-algebra $\mathfrak{F}(\EA)$, equipped
    with the adjoint action of $H$, is $H$-covariantly $^*$-Morita
    equivalent to $\mathcal{A}$ via
    $\Bimod{}{\mathfrak{F}(\mathcal{E}_{\mathcal{A}})}{\mathcal{E}}{}{\mathcal{A}}$.
\end{proposition}
\begin{proof}
    The non-covariant part of this proposition is well-known, see
    Ara's work \cite{ara:1999a} as well as the discussion in
    \cite{bursztyn.waldmann:2003a:pre}. Thus we only have to determine
    the $H$-action induced on $\mathfrak{F}(\EA)$ by the isomorphism
    \eqref{eq:BisotoFEA}. Since
    \[
    g \acts (b \cdot x) 
    = (g _\sweedler{1} \acts b) \cdot (g_\sweedler{2} \acts x)
    \]
    by compatibility, we see by
    Proposition~\ref{proposition:InducedAdjointAction} that this is
    precisely the defining property of the adjoint action. The other
    direction also follows directly from this observation.
\end{proof}
\begin{remark}
    \label{remark:InjectiveActions}
    It follows from the proposition that the maps $b \mapsto (x
    \mapsto b \cdot x)$ as well as $a \mapsto (x \mapsto x \cdot a)$
    are \emph{injective} for an equivalence bimodule.
\end{remark}
\begin{remark}
    \label{remark:SMEHaltMitPos}
    The case of $H$-covariant strong Morita equivalence is analogous
    with the only additional requirement that both inner products are
    completely positive.
\end{remark}
\begin{example}
    \label{example:StandardExampleHSME}
    As usual the standard example is the Morita equivalence of
    $\mathcal{A}$ and $M_n(\mathcal{A})$ via the bimodule
    $\mathcal{A}^n$ where $\mathcal{A}$ acts componentwisely from the
    right and $M_n(\mathcal{A})$ acts by matrix multiplication from
    the left. The canonical, completely positive, full and
    non-degenerate inner product is
    \begin{equation}
        \label{eq:CanonicalInnerProdAn}
        \SPA{x,y} = \sum_{i=1}^n x_i^*y_i,
    \end{equation}
    which determines $\MnASP{\cdot,\cdot}$ by
    compatibility~\eqref{eq:SPCompatible}. The $H$-action on
    $\mathcal{A}^n$ is componentwise and the induced $^*$-action on
    $M_n(\mathcal{A}) =
    \mathfrak{F}(\mathcal{A}^n_{\scriptscriptstyle{\mathcal{A}}})$ is
    just the one from Example~\ref{example:StarActions}, part
    \textit{iii.)}. Thus we get the $H$-covariant strong Morita
    equivalence of $\mathcal{A}$ and $M_n(\mathcal{A})$.
\end{example}

One of the original aims of Morita theory is to establish the
equivalence of representation theories. In our case this is based on
the following observation inspired by \cite[Lem.~5.13 \&
Lem.~5.14]{bursztyn.waldmann:2003a:pre}:
\begin{proposition}
    \label{proposition:TensorTildeNett}
    Let $\mathcal{A}$, $\mathcal{B}$, $\mathcal{C}$, $\mathcal{D}$ be
    idempotent and non-degenerate $^*$-algebras with $^*$-actions of
    $H$. Let $\CFB$ and $\BEA$ be $H$-covariant $^*$-equivalence
    bimodules and let $\AHD \in \sMod[\mathcal{D}, H](\mathcal{A})$ be
    a strongly non-degenerate $^*$-representation of $\mathcal{A}$
    such that in addition $\AHD \cdot \mathcal{D} = \AHD$.
    \begin{compactenum}
    \item One has
        \begin{equation}
            \label{eq:TenssTildeTensHatAssoc}
            \left(\CFB \tensB[\mathcal{B}] \BEA\right)
            \tensM[\mathcal{A}] \AHD
            \cong
            \CFB \tensM[\mathcal{B}] 
            \left(\BEA \tensM[\mathcal{A}] \AHD\right)
        \end{equation}
        via the usual natural $H$-covariant isometric isomorphism.
    \item One has
        \begin{equation}
            \label{eq:AAAUnit}
            \AAA \tensM[\mathcal{A}] \AHD 
            \cong \AHD \cong
            \AHD \tensM[\mathcal{D}] \DDD
        \end{equation}
        via the canonical $H$-covariant isometric isomorphisms $a
        \otimes x \mapsto a \cdot x$ and $x \otimes d \mapsto x \cdot
        d$, respectively.
    \item One has
        \begin{equation}
            \label{eq:AccEBInverserBim}
            \AccEB \tensB[\mathcal{B}] \BEA \cong \AAA
            \quad
            \textrm{and}
            \quad
            \BEA \tensB[\mathcal{A}] \AccEB
            \cong \BBB
        \end{equation}
        via the natural $H$-covariant isometric isomorphisms $\cc{x}
        \otimes y \mapsto \SPEA{x,y}$ and $y \otimes \cc{x} \mapsto
        \BSPE{y,x}$, respectively.
    \end{compactenum}
\end{proposition}
\begin{proof}
    The only thing to be checked is that the isomorphisms are
    compatible with the $H$-actions. The remaining properties where
    already shown in \cite[Lem~5.13,
    Lem~5.14]{bursztyn.waldmann:2003a:pre}. The compatibility for the
    first part is contained in \eqref{eq:GFEAsso}. The action on
    $\mathcal{A} \otimes \mathcal{H}$ is by definition $g \acts (a
    \otimes x) = g_\sweedler{1} \acts a \otimes g_\sweedler{2} \acts
    x$ whence $g \acts (a \otimes x)$ is mapped to $(g_\sweedler{1}
    \acts a) \cdot (g_\sweedler{2} \acts x) = g \acts (a \cdot x)$
    under the isomorphism~\eqref{eq:AAAUnit}. This shows the second
    part as the argument for $\DDD$ is analogous. For the third part
    recall that the action on the complex conjugate bimodule is $g
    \ccacts \cc{x} = \cc{S(g)^* \acts x}$ whence the action on
    $\cc{\mathcal{E}} \otimes \mathcal{E}$ is given by $g \acts
    (\cc{x} \otimes y) = \cc{S(g_\sweedler{1})^* \acts x} \otimes
    g_\sweedler{2} \acts y$.  Thus $g \acts (\cc{x} \otimes y)$ is
    mapped under \eqref{eq:AccEBInverserBim} to $\SPEA{S(g_{1})^*
      \acts x, g_\sweedler{2} \acts y} = g \acts \SPEA{x,y}$ by
    $H$-covariance of the inner product showing the $H$-covariance of
    the first isomorphism. The $H$-covariance of the second
    isomorphism in \eqref{eq:AccEBInverserBim} is analogous.
\end{proof}
\begin{corollary}
    \label{corollary:RieffelFunctorsTensorProducts}
    For equivalence bimodules $\CFB$ and $\BEA$ there is
    a natural equivalence
    \begin{equation}
        \label{eq:RFRERFtensorE}
        \mathsf{R}_{\mathcal{F}} \circ \mathsf{R}_{\mathcal{E}} \cong
        \mathsf{R}_{\mathcal{F} \tensB \mathcal{E}}
    \end{equation}
    for the $H$-covariant Rieffel induction functors.
    Furthermore, when restricted to $\sMod$ (or $\Rep$ in the
    completely positive case, respectively) there are natural
    equivalences
    \begin{equation}
        \label{eq:RAIdentitaet}
        \mathsf{R}_{\mathcal{A}} 
        \cong \id_{^*\textrm{-}\mathsf{Mod}(\mathcal{A})},
    \end{equation}
    \begin{equation}
        \label{eq:RieffelccEinvRieffelE}
        \mathsf{R}_{\mathcal{E}} \circ \mathsf{R}_{\cc{\mathcal{E}}}
        \cong \id_{^*\textrm{-}\mathsf{Mod}(\mathcal{B})}
        \quad
        \textrm{and}
        \quad
        \mathsf{R}_{\cc{\mathcal{E}}} \circ \mathsf{R}_{\mathcal{E}}
        \cong \id_{^*\textrm{-}\mathsf{Mod}(\mathcal{A})}
    \end{equation}
    for the $H$-covariant Rieffel induction functors. Analogous
    statements hold for the functor $\mathsf{S}_{\mathcal{E}}$.
\end{corollary}
\begin{corollary}
    \label{corollary:RieffelHSMEEquivalence}
    Let $\mathcal{A}$, $\mathcal{B}$ be $H$-covariantly $^*$-Morita
    equivalent via $\BEA$. Then
    \begin{equation}
        \label{eq:RieffelEEquivalencesMod}
        \mathsf{R}_{\mathcal{E}}: 
        \sMod[\mathcal{D}, H](\mathcal{A})
        \stackrel{\cong}{\longrightarrow}
        \sMod[\mathcal{D}, H](\mathcal{B})
    \end{equation}
    is an equivalence of categories with `inverse'
    $\mathsf{R}_{\cc{\mathcal{E}}}$. If in addition $\BEA$ is even a
    $H$-covariant strong Morita equivalence bimodule, then
    $\mathsf{R}_{\mathcal{E}}$
    restricts to an equivalence
    \begin{equation}
        \label{eq:RepAequivRepB}
        \mathsf{R}_{\mathcal{E}}:
        \Rep[\mathcal{D}, H](\mathcal{A})
        \stackrel{\cong}{\longrightarrow}
        \Rep[\mathcal{D}, H](\mathcal{B}).
    \end{equation}
\end{corollary}

This is the $H$-covariant version of
\cite[Cor.~5.15]{bursztyn.waldmann:2003a:pre} which itself is the
algebraic generalization of Rieffel's theorem on equivalent
$^*$-representation theories of $C^*$-algebras \cite{rieffel:1974b}.
In the case of $C^*$-algebras and strongly continuous group actions of
locally compact groups analogous statements are well-known, see e.g.
the discussion in \cite[Sect.~7.2]{raeburn.williams:1998a}. It is an
interesting problem whether and how one can use our purely algebraic
approach to obtain those results. We will address these questions in
future projects.

%
%

\section{The $H$-covariant Picard groupoid and Morita invariants}
\label{sec:CovPicardInvariants}

As already mentioned, Morita theory can be seen as resulting from an
extended notion of morphisms between algebras: one considers
isomorphism classes of bimodules as morphisms and obtains a new
category with the same underlying class of objects but bigger classes
of morphisms. Isomorphism in this category is then precisely Morita
equivalence. This point of view is classical for ring-theoretic Morita
equivalence, see e.g.~\cite{bass:1968a, benabou:1967a}, and it was
discussed in detail for the $^*$- and strong Morita equivalence of
$^*$-algebras in \cite[Sect.~6]{bursztyn.waldmann:2003a:pre}. See also
Landsman's work \cite{landsman:2001b, landsman:2001c} in the context
of $C^*$-algebras.  Alternatively, one could use a bicategorical
approach by not identifying the bimodules up to isomorphism in a first
step.

Important for us is that any such enlarged category defines its
groupoid of invertible arrows, the corresponding \textdef{Picard
  groupoid}.  Strictly speaking, this is not an honest groupoid for
two reasons: first the class of units (here the class of
$^*$-algebras) is not a set, so it can not be a \emph{small} category.
Second, the class of invertible arrows between two units is, a priori,
not known to be a set either.  This is more severe, but in the case of
unital $^*$-algebras one actually can show that the space of arrows
between two units in the Picard groupoid forms a set as it is given by
equivalence classes of certain \emph{finitely generated projective
  modules}.  Thus we shall ignore these subtleties in the following
and focus mainly on the unital case.  In any case, throughout this
section all algebras will be idempotent and non-degenerate.

%
%

\subsection{The $H$-covariant Picard groupoid}
\label{subsec:HcovPicard}

Instead of defining the Picard groupoid in the above described way, we
give a more direct definition using the equivalence bimodules
directly. Both approaches are completely equivalent which can easily
be obtained from an $H$-covariant version of
\cite[Thm.~6.1]{bursztyn.waldmann:2003a:pre} using
Proposition~\ref{proposition:ActionOnBisTheOneOfF}.
\begin{definition}
    \label{definition:PicStrH}
    Let $\mathcal{A}$, $\mathcal{B}$ be $^*$-algebras over $\ring{C}$
    and define $\starPicH(\mathcal{B}, \mathcal{A})$ to be the class
    of isomorphism classes of $H$-covariant $^*$-Morita equivalence
    bimodules $\BEA$ and set $\starPicH(\mathcal{A}) =
    \starPicH(\mathcal{A}, \mathcal{A})$. Similarly, we define
    $\StrPicH(\mathcal{B}, \mathcal{A})$ to be the class of
    isomorphism classes of $H$-covariant strong Morita equivalence
    bimodules $\BEA$ and set $\StrPicH(\mathcal{A}) =
    \StrPicH(\mathcal{A}, \mathcal{A})$.
\end{definition}
Here and in the following `isomorphism' of equivalence bimodules
includes all relevant structures, i.e.~the $H$-action, the bimodule
structure as well as the inner products.
\begin{theorem}
    \label{theorem:PicardGroupoidIsGroupoid}
    Viewing $\starPicH(\mathcal{B}, \mathcal{A})$ as space of arrows
    $\mathcal{A} \longrightarrow \mathcal{B}$ one obtains the
    $H$-covariant $^*$-Picard groupoid $\starPicH$, where the
    composition law is $\tensB$, the units are the $^*$-algebras
    themselves with the classes of the canonical bimodules
    $\left[\AAA\right]$ as unit arrows. The inverse arrows are the
    classes of the complex conjugated bimodules. Similarly, one obtains
    the $H$-covariant strong Picard groupoid $\StrPicH$.
\end{theorem}
The proof is obvious by use of
Proposition~\ref{proposition:TensorTildeNett} and the fact that
$\tensB$ is functorial and hence well-defined on isomorphism classes.
\begin{remark}
    \label{remark:UnitalPicCase}
    For unital $^*$-algebras $\starPicH(\mathcal{B}, \mathcal{A})$ as
    well as $\StrPicH(\mathcal{B}, \mathcal{A})$ are in bijection to
    certain finitely generated projective modules and hence they are
    sets. Thus $\starPicH$ as well as $\StrPicH$ become `large'
    groupoids in this case. For non-unital $^*$-algebras this is a
    priori not clear. Dropping the information about the inner
    products one obtains the ring-theoretic notions of the Picard
    groupoid which we denote by $\PicH$ and $\Pic$, respectively, see
    also \cite{bass:1968a, benabou:1967a}.
\end{remark}
The isotropy groups $\StrPicH(\mathcal{A})$ and
$\starPicH(\mathcal{A})$, respectively, of the Picard groupoids are
called the $H$-covariant \textdef{strong} (resp. \textdef{$^*$-})
\textdef{Picard groups} of $\mathcal{A}$.

By successively forgetting the additional structures one obtains
groupoid morphisms:
\begin{equation}
    \label{eq:NiceTrianglePicDiagram}
    \bfig
    \Vtriangle(0,350)<500,200>[\StrPicH`\starPicH`\Pic_H;``]
    \Vtriangle(0,0)/@{>}|\hole`>`>/<500,200>[\StrPic`\starPic,`\Pic;``]
    \morphism(0,525)<0,-250>[`;]
    \morphism(500,325)<0,-250>[`;]
    \morphism(1000,525)<0,-250>[`;]
    \efig
\end{equation}
Each of them induces the identity on the units of the groupoids.
Clearly, all combinations of possible compositions commute in this
diagram. Thus an interesting program will be to investigate under
which reasonable restrictions and conditions on the $^*$-algebras and
the $H$-actions one can say something on the images and kernels of
these groupoid morphisms. In the situation without $H$-action the
lower commuting triangle in \eqref{eq:NiceTrianglePicDiagram} has been
investigated in some detail in \cite{bursztyn.waldmann:2003a:pre} for
a large class of unital $^*$-algebras.

Before investigating \eqref{eq:NiceTrianglePicDiagram} we relate the
Picard groupoid to the \textdef{isomorphism groupoid} as we want to
interpret the elements in $\starPicH(\mathcal{B}, \mathcal{A})$ and
$\StrPicH(\mathcal{B}, \mathcal{A})$ as generalized isomorphisms of
$^*$-algebras. We denote by $\starIsoH(\mathcal{B}, \mathcal{A})$ the
\textdef{$H$-equivariant $^*$-isomorphisms} from $\mathcal{A}$ to
$\mathcal{B}$ and set $\starAutH(\mathcal{A}) = \starIsoH(\mathcal{A},
\mathcal{A})$ for the \textdef{$H$-equivariant $^*$-automorphisms} of
$\mathcal{A}$. The non-equivariant case is denoted by
$\starIso(\mathcal{B}, \mathcal{A})$ and $\starAut(\mathcal{A})$,
respectively.
\begin{remark}
    \label{remark:IsoGroupoid}
    Viewing $\starIsoH(\mathcal{B}, \mathcal{A})$ as space of arrows
    from $\mathcal{A}$ to $\mathcal{B}$ one obtains the usual (large)
    groupoid of $H$-equivariant $^*$-isomorphisms with the
    $H$-equivariant $^*$-automorphism groups as isotropy groups.
\end{remark}

Let us first recall and adapt some results and definitions from
\cite[Sect.~2]{bursztyn.waldmann:2004a}. Let $\Phi \in
\starIso(\mathcal{B},\mathcal{A})$ be given and let $\CFB$ be a
representative for a class $\left[\CFB\right] \in
\starPic(\mathcal{C},\mathcal{B})$ (or in
$\StrPic(\mathcal{C},\mathcal{B})$, respectively), Then we can twist
$\mathcal{F}$ into a right $\mathcal{A}$-module by setting
\begin{equation}
    \label{eq:PhiTwistModule}
    x \dotPhi a = x \cdot \Phi(a)
\end{equation}
for $x \in \mathcal{F}$ and $a \in \mathcal{A}$. Clearly, this is
gives a $(\mathcal{C}, \mathcal{A})$-bimodule, denoted by
$\CFPhiA$. Moreover, we define
\begin{equation}
    \label{eq:PhiTwistInnerProduct}
    \SPFPhiA{x, y} = \Phi^{-1} \left(\SPFB{x, y}\right).
\end{equation}
Since $\Phi$ is a $^*$-isomorphism, this gives a full and
non-degenerate $\mathcal{A}$-valued inner product on $\CFPhiA$
(completely positive in the case of a strong equivalence bimodule)
which is compatible with the $\mathcal{C}$-module structure and with
the $\mathcal{C}$-valued inner product on $\mathcal{F}$. Thus we
obtain a $^*$- respectively strong Morita equivalence bimodule
$\CFPhiA$. A last simple check ensures that the class
$\left[\CFPhiA\right]$ only depends on the class $\left[\CFB\right]$.
Similarly, we can twist equivalence bimodules $\BEA$ from the left
with some $\Psi \in \starIso(\mathcal{C}, \mathcal{B})$ by setting
\begin{equation}
    \label{eq:LeftPsiTwists}
    c \dotPsi x = \Psi^{-1}(c) \cdot x
    \quad
    \textrm{and}
    \quad
    \CSPPsiE{x, y} = \Psi\left(\BSPE{x, y}\right)
\end{equation}
and obtain an equivalence bimodule $\CPsiEA$. Again, this works either
for $^*$-equivalence or strong equivalence bimodules. The
$H$-covariant situation is as follows:
\begin{lemma}
    \label{lemma:PhiTwistHcov}
    Let $\Phi \in \starIso(\mathcal{B}, \mathcal{A})$ and let
    $\left[\CFB\right] \in \starPicH(\mathcal{C}, \mathcal{B})$ (or in
    $\StrPicH(\mathcal{B}, \mathcal{A})$, respectively). Then
    $\CFPhiA$ is an $H$-covariant $^*$- (or strong, respectively)
    equivalence bimodule if and only if $\Phi$ is $H$-equivariant. In
    this case $\left[\CFPhiA\right] \in \starPicH(\mathcal{C},
    \mathcal{A})$ (or in $\starPicH(\mathcal{C}, \mathcal{A})$,
    respectively) is well-defined.
\end{lemma}
\begin{proof}
    If $\Phi$ is $H$-equivariant, then it follows from a simple
    computation that the twisted bimodule is an $H$-covariant
    equivalence bimodule as well. For the converse direction, we note
    that if the twisted bimodule is still $H$-covariant then
    \[
    x \cdot \Phi(g \acts a)
    = x \dotPhi (g \acts a)
    = g_\sweedler{2} \acts 
    \left( (S^{-1}(g_\sweedler{2}) \acts x) \dotPhi a \right)
    = g_\sweedler{2} \acts 
    \left( (S^{-1}(g_\sweedler{2}) \acts x) \cdot \Phi (a) \right)
    = x \cdot (g \acts \Phi(a)).
    \]
    Thanks to Remark~\ref{remark:InjectiveActions} the
    $H$-equivariance of $\Phi$ follows. The well-definedness of the
    class is clear.
\end{proof}

The analogous statement holds if we twist with an $H$-equivariant
$^*$-isomorphism from the left.  In particular, we can twist the
`unit' bimodules (as we have already done implicitly in the proof of
Theorem~\ref{theorem:MoritaEquivalenceRelation}) continuing the
discussion of \cite{bursztyn.waldmann:2004a} as well as the well-known
ring-theoretic case.
\begin{proposition}
    \label{proposition:IsoToPicMorphism}
    Let $\Phi \in \starIsoH(\mathcal{B}, \mathcal{A})$ and denote by
    $\ell(\Phi) \in \StrPicH(\mathcal{B}, \mathcal{A})$ the class of
    the bimodule $\BBPhiA$.
    \begin{compactenum}
    \item We have $\big[\BBPhiA\big] = \big[\BPhiAA\big]$.
    \item The map
        \begin{equation}
            \label{eq:ellIsoToPic}
            \ell: \starIsoH \longrightarrow \StrPicH
        \end{equation}
        is a groupoid morphism inducing the identity on the units.
    \item For $\left[\CFB\right] \in \StrPicH(\mathcal{C},
        \mathcal{B})$ we have
        \begin{equation}
            \label{eq:CFBellPhiCFPhiA}
            \left[\CFB\right] \tensB \ell(\Phi) 
            =
            \left[\CFPhiA\right].
        \end{equation}
    \end{compactenum}
    We can replace strong by $^*$-Picard groupoids as well.
\end{proposition}
\begin{proof}
    The bimodule isomorphism for the first part is simply given by $b
    \mapsto \Phi^{-1}(b)$. Now let $\Phi \in \starIsoH(\mathcal{B},
    \mathcal{A})$ and $\Psi \in \starIsoH(\mathcal{C}, \mathcal{B})$
    be given. Then we consider $\ell(\Psi \circ \Phi)$ and compare it
    with $\ell(\Psi) \tensB \ell(\Phi)$. We consider the map defined
    by
    \[
    \CCPsiB \tensB[\mathcal{B}] \BBPhiA \ni
    c \otimes b \mapsto
    c \Psi(b) \in \CCPsiPhiA.
    \]
    On the level of $\tensor[\mathcal{B}]$ rather than
    $\tensB[\mathcal{B}]$ it is easy to see that this map is
    well-defined over $\tensor[\mathcal{B}]$. Moreover, it is
    surjective since $\mathcal{C}$ is idempotent. A straightforward
    check shows that it is a $(\mathcal{C}, \mathcal{A})$-bimodule
    morphism isometric with respect to both inner products. Thus the
    quotient by the degeneracy spaces yields an injective map,
    well-defined over $\tensB[\mathcal{B}]$. Hence we end up with a
    bimodule isomorphism. A last simple computation using the
    $H$-equivariance of $\Psi$ shows that it is even an $H$-equivariant
    isomorphism as wanted. This proves the second part as
    \eqref{eq:ellIsoToPic} clearly maps the unit $\id_{\mathcal{A}}$
    to the unit $\left[\AAA\right]$. For the last part we check that
    \[
    \CFB \tensB[\mathcal{B}] \BBPhiA \ni x \otimes b
    \mapsto x \cdot b
    \in \CFPhiA
    \]
    is the desired isomorphism. This is again a straightforward
    computation. In the proof the positivity of the inner products was
    not essential.
\end{proof}

At least for unital $^*$-algebras on can describe the \emph{kernel} of
the groupoid morphism \eqref{eq:ellIsoToPic} rather explicitly. We
slightly extend and specialize the arguments from
\cite[Prop.~2.3]{bursztyn.waldmann:2004a} for our purposes. First we
denote by $\starInnAutH(\mathcal{A})$ those \textdef{inner
  $^*$-automorphisms} $a \mapsto u a u^{-1}$ where $u^* = u^{-1}$ is
unitary and $H$-invariant $g \acts u = \epsilon(g) u$. Clearly,
$\starInnAutH(\mathcal{A}) \subseteq \starAutH(\mathcal{A})$ is a
normal subgroup.
\begin{proposition}
    \label{proposition:InnAutAutPicSequenz}
    Let $\mathcal{A}$, $\mathcal{B}$ be unital $^*$-algebras.
    \begin{compactenum}
    \item For $\Phi \in \starAutH(\mathcal{B})$ and $\left[\BEA\right]
        \in \StrPicH(\mathcal{B}, \mathcal{A})$ we have
        $\left[\PhiBEA\right] = \left[\BEA\right]$ if and only if
        $\Phi \in \starInnAutH(\mathcal{B})$.
    \item We have the exact sequence of groups
        \begin{equation}
            \label{eq:InnAutAutPicSequenz}
            1 
            \longrightarrow
            \starInnAutH(\mathcal{A})
            \longrightarrow
            \starAutH(\mathcal{A})
            \stackrel{\ell}{\longrightarrow}
            \StrPicH(\mathcal{A})
        \end{equation}
    \end{compactenum}
    Again, strong can be replaced by $^*$-Picard groupoids.
\end{proposition}
\begin{proof}
    Assume $U: \PhiBEA \longrightarrow \BEA$ is an isomorphism. Then
    $U(x \cdot a) = U(x) \cdot a$ implies that there exists an
    invertible element $u \in \mathcal{B}$ with $U(x) = u \cdot x$
    thanks to Proposition~\ref{proposition:ActionOnBisTheOneOfF} and
    since $\mathcal{B}$ is unital. Then $b \cdot u \cdot x = b \cdot
    U(x) = U(b \dotPhi x) = u \cdot \Phi^{-1}(b) \cdot x$ implies
    $\Phi(b) = u b u^{-1}$ thanks to
    Remark~\ref{remark:InjectiveActions}. Note that $\Phi$ being a
    $^*$-automorphism does not necessarily imply that $u$ is
    unitary. Nevertheless, we have by isometry of $U$
    \[
    u \BSPE{x, y} u^*
    =
    \BSPE{u \cdot x, u \cdot y}
    =
    \BSPE{U(x), U(y)}
    =
    \BSPPhiE{x, y}
    =
    \Phi(\BSPE{x, y}),
    \]
    whence by fullness $\Phi(b) = u b u^*$. Thus $u^* = u^{-1}$
    turns out to be unitary. Finally, we have from $g \acts U(x) = U(g
    \acts x)$ the relation
    \[
    (g \acts u) \cdot x 
    = g_\sweedler{1} \acts (u \cdot S(g_\sweedler{2}) \acts x)
    = g_\sweedler{1} \acts (U(S(g_\sweedler{2}) \acts x))
    = U(g_\sweedler{1} \acts (S(g_\sweedler{2}) \acts x))
    = \epsilon(g) u \cdot x.
    \]
    Again by Remark~\ref{remark:InjectiveActions} we see $g \acts u =
    \epsilon(g) u$. This proves the first statement as the converse
    is a trivial computation. The second part is then an easy
    consequence if we apply the first part to $\mathcal{B} =
    \mathcal{A}$ and $\BEA = \AAA$.
\end{proof}

From this proposition we see that $\StrPicH$ as well as $\starPicH$
indeed generalize $\starIsoH$ in a very precise way. Though the
kernel of $\ell$ in \eqref{eq:ellIsoToPic} can be described by this
proposition explicitly, the lack of surjectivity usually depends very
much on the example.

%
%

\subsection{The groupoid morphism $\StrPicH{} \longrightarrow \StrPic{}$}
\label{subsec:HPicToPicmorphism}

We shall now discuss the canonical groupoid morphism of `forgetting'
the $H$-covariance
\begin{equation}
    \label{eq:PicHtoPic}
    \StrPicH \longrightarrow \StrPic
\end{equation}
where we treat the case of the strong Picard groupoids. For the
$^*$-Picard groupoids and the ring-theoretic Picard groupoids the
results will be analogous.

In general, the question whether $\StrPicH(\mathcal{B}, \mathcal{A})
\longrightarrow \StrPic(\mathcal{B}, \mathcal{A})$ is surjective for
two given $^*$-algebras with $^*$-action of $H$ is very difficult and
depends very much on the example. The problem is to `lift' the action
of $H$ from the algebras to an equivalence bimodule $\BEA$. In general
there will be obstructions for this lifting. 
\begin{example}
    \label{example:Lifting}
    Coming back to the Examples~\ref{example:LieAlgebraAction} and
    \ref{example:LieGroupAction}, we notice that for $\mathcal{A} =
    C^\infty(M)$ any equivalence bimodule is of the form
    $\Gamma^\infty(E)$ where $E \longrightarrow M$ is a complex vector
    bundle and $C^\infty(M)$ acts by pointwise multiplication from the
    right. Then the Morita equivalent algebra $\mathcal{B}$ is
    isomorphic to $\Gamma^\infty(\End(E))$ where the action is by
    pointwise application of the endomorphism. Now if a Lie algebra
    action of $\lie{g}$ on $M$ is given, then the question is whether
    one can lift this action to an action on the sections of $E$. In
    general, there are obstructions: Consider a Lie group $G$ with
    action of its Lie algebra $\lie{g}$ by left invariant vector
    fields. Suppose $E \longrightarrow G$ is a complex vector bundle
    admitting a lifted action $\xi \mapsto \Lie_\xi$ of $\lie{g}$,
    i.e. $\Lie_\xi: \Gamma^\infty(E) \longrightarrow \Gamma^\infty(E)$
    is a representation of $\lie{g}$ and satisfies $\Lie_\xi(s f) =
    \Lie_\xi (s) f + s X_\xi(f)$, for all $\xi \in \lie{g}$, where
    $X_\xi$ is the corresponding left invariant vector field on $G$.
    Since for a basis $e_1, \ldots, e_n$ of $\lie{g}$ we obtain a
    module basis $X_{e_1}, \ldots, X_{e_n}$ of all vector fields
    $\Gamma^\infty(TG)$ on $G$, we define by
    \begin{equation}
        \label{eq:FlacheKovAbl}
        \nabla_Y = \sum_{i=1}^n Y^i \Lie_{X_{e_i}},
        \quad
        \textrm{with}
        \quad
        Y = \sum_{i=1}^n Y^i X_{e_i}
        \quad
        \textrm{and}
        \quad
        Y^i \in C^\infty(G)
    \end{equation}
    a covariant derivative, which is easily shown to be \emph{flat}.
    In general the existence of a flat covariant derivative is a
    cohomological obstruction on $E$, unless $E$ is a trivial vector
    bundle.
\end{example}

The question about injectivity is in how many ways such a lifting can
be done.  Surprisingly, there is a general answer to this question
which is even independent on the particular bimodule but universal for
all bimodules $\BEA$ as long as they allow for lifting at all.

In the following we fix a strong Morita equivalence bimodule $\BEA$
and assume that there is at least one $H$-action $\acts$ on $\BEA$
such that it becomes an $H$-covariant strong Morita equivalence
bimodule. If $\actsp$ is another such $H$-action then we define
\begin{equation}
    \label{eq:ugDef}
    u_g (x) 
    = g_\sweedler{1} \acts 
    \left( S(g_\sweedler{2}) \actsp x \right)
\end{equation}
to `measure' the difference between the two actions, where $x \in
\BEA$. Knowing $\acts$ and all the maps $g \mapsto u_g \in
\End_{\ring{C}}(\BEA)$ allows to reconstruct $\actsp$ by
\begin{equation}
    \label{eq:actspFromuacts}
    g \actsp x = 
    g_\sweedler{2} \acts 
    \left(u_{S^{-1}(g_\sweedler{1})} (x)\right)
\end{equation}
and conversely $\acts$ is determined by
\begin{equation}
    \label{eq:actsFromugAndactsp}
    g \acts x = 
    u_{g_\sweedler{1}} \left( g_\sweedler{2} \actsp x\right).
\end{equation}
Thus we have to investigate the maps $u_g$ and find conditions such
that for a given $H$-action, say $\actsp$, the formula
\eqref{eq:actsFromugAndactsp} defines again an $H$-action with the same
properties.
\begin{lemma}
    \label{lemma:ugRightAlinear}
    Let $(\BEA, \actsp)$ be an $H$-covariant strong Morita equivalence
    bimodule and let $H \ni g \mapsto u_g \in \End_{\ring{C}}(\BEA)$ be a
    linear map. Then $\acts$ defined by \eqref{eq:actsFromugAndactsp}
    satisfies $g \acts (x \cdot a) = (g_\sweedler{1} \acts x) \cdot
    (g_\sweedler{2} \acts a)$ if and only if $u_g$ is right
    $\mathcal{A}$-linear for all $g \in H$.
\end{lemma}
\begin{proof}
    First we assume that $u_g$ is right $\mathcal{A}$-linear. Then
    \[
    \begin{split}
        g \acts (x \cdot a)
        &= 
        u_{g_\sweedler{1}}\left(g_\sweedler{2} \actsp (x \cdot a)\right)
        = u_{g_\sweedler{1}}\left((g_\sweedler{2} \actsp x) \cdot
            (g_\sweedler{3} \acts a)\right) \\
        &= 
        \left(u_{g_\sweedler{1}} (g_\sweedler{2} \actsp x)\right)
        \cdot (g_\sweedler{3} \acts a)
        = (g_\sweedler{1} \acts x) \cdot (g_\sweedler{2} \acts a).
    \end{split}
    \]
    For the converse, note first that $g_\sweedler{1} \acts
    (S(g_\sweedler{2}) \actsp x) = u_g(x)$ since $\actsp$ is an action
    (whether $\acts$ is an action or not). If $\acts$ is an action,
    \[
    u_g(x \cdot a)
    = g_\sweedler{1} \acts 
    \left(
        (S(g_\sweedler{3}) \actsp x)
        \cdot
        (S(g_\sweedler{2}) \acts a)
    \right)
    = \left(g_\sweedler{1} \acts \left(S(g_\sweedler{3}) \actsp
            x\right)\right)
    \cdot (\epsilon(g_\sweedler{2}) a)
    = u_g(x) \cdot a.
    \]
\end{proof}

Thus we have to investigate right $\mathcal{A}$-linear endomorphisms
of $\BEA$. Now the crucial observation is that in the \emph{unital}
case any right $\mathcal{A}$-linear endomorphism is a left
multiplication by a unique element in $\mathcal{B}$. To make use of
this drastic simplification we shall assume that in this section all
$^*$-algebras are unital.

Thus we can rephrase Lemma~\ref{lemma:ugRightAlinear} in the following
way: If we want to pass from $\acts$ to $\actsp$ then it will be
necessary and sufficient to consider a map $u_g$ of the form
\begin{equation}
    \label{eq:ugIsLeftMultWithbg}
    u_g(x) = \twist{b}(g) \cdot x 
    = g_\sweedler{1} \acts (S(g_\sweedler{2}) \actsp x),
\end{equation}
if we want to keep the compatibility with the right
$\mathcal{A}$-module structure. Here $\twist{b} \in \Hom_{\ring{C}}(H,
\mathcal{B})$.

The following proposition clarifies under which conditions on
$\twist{b}$ we stay in the class of $H$-covariant strong Morita
equivalence bimodules.
\begin{proposition}
    \label{proposition:twistbIStwist}
    Let $(\BEA, \actsp)$ be an $H$-covariant strong Morita equivalence
    bimodule and let $\twist{b} \in \Hom_{\ring{C}}(H,
    \mathcal{B})$. Then for $\acts$ defined by 
    \begin{equation}
        \label{eq:actsByTwistbactsp}
        g \acts x = 
        \twist{b}(g_\sweedler{1}) \cdot (g_\sweedler{2} \actsp x)
    \end{equation}
    one has the following properties:
    \begin{compactenum}
    \item $\Unit_H \acts x = x$ if and only if $\twist{b}(\Unit_H) =
        \Unit_{\mathcal{B}}$.
    \item $\acts$ is an $H$-action if and only if for all $g, h \in H$
        \begin{equation}
            \label{eq:twistActionCondition}
            \twist{b}(gh) = 
            \twist{b}(g_\sweedler{1}) 
            (g_\sweedler{2} \acts \twist{b}(h)). 
        \end{equation}
    \item $\acts$ is compatible with the left $\mathcal{B}$-module
        structure if and only if for all $g \in H$ and $b \in
        \mathcal{B}$
          \begin{equation}
              \label{eq:twistModuleCondition}
              (g_\sweedler{1} \acts b) \twist{b}(g_\sweedler{2})
              =
              \twist{b}(g_\sweedler{1}) (g_\sweedler{2} \acts b).
          \end{equation}
      \item $\acts$ is compatible with the inner product $\BSPE{\cdot,
            \cdot}$ if and only if $\twist{b}$ fulfills
          \eqref{eq:twistModuleCondition} and for all $g \in H$
          \begin{equation}
              \label{eq:twistUnitaryCondition}
              \twist{b}(g_\sweedler{1})
              \left(\twist{b}(S(g_\sweedler{2})^*)\right)^*
              = \epsilon(g) \Unit_{\mathcal{B}}.
          \end{equation}
    \end{compactenum}
    If $\acts$ fulfills \textit{i.)}--\textit{iv.)} then $\acts$ is
    compatible with the inner product $\SPEA{\cdot,\cdot}$, too.
\end{proposition}
\begin{proof}
    The first part is trivial. For the second we compute under
    assumption of \eqref{eq:twistActionCondition}
    \[
    \begin{split}
        g \acts (h \acts x)
        &= 
        \twist{b}(g_\sweedler{1}) \cdot 
        \left(
            g_\sweedler{2} \actsp \left( 
                \twist{b}(h_\sweedler{1})
                \cdot (h_\sweedler{2} \actsp x)
            \right)
        \right) \\
        &= 
        \left( 
            \twist{b}(g_\sweedler{1}) 
            (g_\sweedler{2} \acts \twist{b}(h_\sweedler{1}))
        \right) 
        \cdot
        \left(
            g_\sweedler{3} \actsp (h_\sweedler{2} \actsp x)
        \right) \\
        &= 
        \twist{b}(g_\sweedler{1}h_\sweedler{1})
        \cdot 
        \left(
            (g_\sweedler{2} h_\sweedler{2}) \actsp x
        \right) \\
        &= 
        (gh) \acts x,
    \end{split}
    \]
    using the compatibility of $\actsp$ with the module structure as
    well as \eqref{eq:twistActionCondition} and that $\actsp$ is an
    action. Conversely, we have $\twist{b}(g) \cdot x = g_\sweedler{1}
    \acts (S(g_\sweedler{2}) \actsp x)$ whether $\acts$ is an action
    or not. Now, if $\acts$ is an action, too, then
    \[
    \begin{split}
        \twist{b}(gh) \cdot x
        &= 
        (g_\sweedler{1} h_\sweedler{1}) \acts 
        \left(
            S(g_\sweedler{2}h_\sweedler{2}) \actsp x
        \right)\\
        &=
        g_\sweedler{1} \acts \left(
            h_\sweedler{1} \acts \left(
                S(g_\sweedler{2}h_\sweedler{2}) \actsp x
            \right)
        \right) \\
        &=
        \twist{b}(g_\sweedler{1}) \cdot\left(
            g_\sweedler{2} \actsp \left(
                \twist{b}(h_\sweedler{1}) \cdot
                \left(
                    h_\sweedler{2} \actsp \left(
                        S(g_\sweedler{3}h_\sweedler{3}) \actsp x
                    \right)
                \right)
            \right)
        \right) \\
        &=
        \left(
            \twist{b}(g_\sweedler{1}) 
            (g_\sweedler{2} \acts \twist{b}(h_\sweedler{1}))
        \right)
        \cdot
        \left(
            \left(g_\sweedler{3}h_\sweedler{2}
                S(g_\sweedler{4}h_\sweedler{3})
            \right)
            \actsp x
        \right) \\
        &= 
        \left(
            \twist{b}(g_\sweedler{1}) 
            (g_\sweedler{2} \acts \twist{b}(h))
        \right)
        \cdot
        x,
    \end{split}                        
    \]
    whence by Remark~\ref{remark:InjectiveActions} the second part
    follows. For the third part we assume
    \eqref{eq:twistModuleCondition} and compute
    \[
    \begin{split}
        g \acts (b \cdot x)
        &=
        \twist{b}(g_\sweedler{1}) \cdot
        (g_\sweedler{2} \actsp (b \cdot x)) \\
        &=
        \left(
            \twist{b}(g_\sweedler{1}) (g_\sweedler{2} \acts b)
        \right) 
        \cdot
        (g_\sweedler{3} \actsp x) \\
        &=
        \left(
            (g_\sweedler{1} \acts b)\twist{b}(g_\sweedler{2})
        \right)
        \cdot
        (g_\sweedler{3} \actsp x) \\
        &=
        (g_\sweedler{1} \acts b)
        \cdot
        (g_\sweedler{2} \acts x).
    \end{split}
    \]
    Conversely, assuming $\acts$ is compatible with the module
    structure gives by a similar computation
    \[
    \left((g_\sweedler{1} \acts b) \twist{b}(g_\sweedler{2})\right)
    \cdot x
    =
    \left(\twist{b}(g_\sweedler{1}) (g_\sweedler{2} \acts b)\right) 
    \cdot x,
    \]
    whence again by Remark~\ref{remark:InjectiveActions} the third
    part follows. For the fourth part, note that
    \eqref{eq:twistModuleCondition} is necessary by
    Lemma~\ref{lemma:NiceAxiomForSP} anyway whence we assume
    \eqref{eq:twistModuleCondition}. Then we have
    \[
    \begin{split}
        \BSPE{g_\sweedler{1} \acts x, S(g_\sweedler{2})^* \acts y}
        &=
        \BSPE{
          \twist{b}(g_\sweedler{1}) \cdot(g_\sweedler{2} \actsp x),
          \twist{b}(S(g_\sweedler{4})^*) \cdot 
          (S(g_\sweedler{3})^* \actsp y)
        } \\
        &=
        \twist{b}(g_\sweedler{1})
        \BSPE{g_\sweedler{2} \actsp x, 
          S(g_\sweedler{3})^* \actsp y}
        \left(\twist{b}(S(g_\sweedler{4})^*)\right)^* \\
        &=
        \twist{b}(g_\sweedler{1})
        \left(
            g_\sweedler{2}
            \acts
            \BSPE{x, y}
        \right)
        \left(\twist{b}(S(g_\sweedler{3})^*)\right)^* \\
        &=
        \left(
            g_\sweedler{1} \acts 
            \BSPE{x, y}
        \right)
        \twist{b}(g_\sweedler{2})
        \left(\twist{b}(S(g_\sweedler{3})^*)\right)^*.
    \end{split}
    \]
    Now if \eqref{eq:twistUnitaryCondition} is fulfilled, then the
    last line gives $g \acts \BSPE{x, y}$ whence $\acts$ is compatible
    with the inner product. Conversely, if $\acts$ is compatible, then
    we obtain from this computation
    \[
    g \acts \BSPE{x, y}
    = 
    \left(
        g_\sweedler{1} \acts 
        \BSPE{x, y}
    \right)
    \twist{b}(g_\sweedler{2})
    \left(\twist{b}(S(g_\sweedler{3})^*)\right)^*.
    \]
    Since the inner product is full we can take linear combinations in
    $x$ and $y$ to get $g \acts \Unit_{\mathcal{B}} = \epsilon(g)
    \Unit_{\mathcal{B}}$ on the left hand side. Then the right hand
    side gives $\epsilon(g_\sweedler{1}) \twist{b}(g_\sweedler{2})
    \left(\twist{b}(S(g_\sweedler{3})^*)\right)^*$ whence
    \eqref{eq:twistUnitaryCondition} follows. From the compatibility
    of the two inner products as in \eqref{eq:SPCompatible} and the
    compatibility of one of them with the $H$-action $\acts$ the
    compatibility of the other with the $H$-action follows in general.
\end{proof}

This proposition has now the following easy interpretation in terms of
the group $\group{U}(H, \mathcal{B})$ as defined in
Definition~\ref{definition:GLHAandUHA}. Clearly, we can exchange the
roles of $\acts$ and $\actsp$ again (only for aesthetic reasons) as we
have a bijective correspondence.
\begin{corollary}
    \label{corollary:UHBISCOOL}
    Let $(\BEA, \acts)$ be an $H$-covariant strong Morita equivalence
    bimodule. Then all other compatible $H$-actions on $\mathcal{E}$
    are parametrized in a unique way by elements $\twist{b} \in
    \group{U}(H, \mathcal{B})$ by
    \begin{equation}
        \label{eq:actsbDef}
        g \actsT[\twist{b}] x = 
        \twist{b}(g_\sweedler{1}) \cdot (g_\sweedler{2} \acts x).
    \end{equation}
\end{corollary}
\begin{proof}
    The four conditions in Proposition~\ref{proposition:twistbIStwist}
    are precisely the defining relations for elements in $\group{U}(H,
    \mathcal{B})$, thereby explaining the names of the conditions in
    Definition~\ref{definition:GLHAandUHA}.
\end{proof}

We want to understand in which case two given actions give an
isomorphic bimodule and hence the same element in the Picard groupoid.
We use some notation from Appendix~\ref{sec:FunnyGroup}.
\begin{lemma}
    \label{lemma:AutBEAACtsOnActs}
    Let $\BEA$ be a strong Morita equivalence bimodule.
    \begin{compactenum}
    \item The group of isometric bimodule automorphisms $\Aut(\BEA)$
        of $\BEA$ is canonically isomorphic to
        $\group{U}(\zentrum(\mathcal{B}))$ via
        \begin{equation}
            \label{eq:UZBisAUTBEA}
            \group{U}(\zentrum(\mathcal{B}))
            \ni c
            \mapsto
            (\Phi_c: x \mapsto c \cdot x) \in \Aut(\BEA).
        \end{equation}
    \item Assume that $\BEA$ allows for compatible $H$-actions such
        that it becomes an $H$-covariant strong Morita equivalence
        bimodule. Then $\Aut(\BEA) = \group{U}(\zentrum(\mathcal{B}))$
        acts on the set of such compatible $H$-actions by
        \begin{equation}
            \label{eq:AutBEAActsOnActs}
            (\Phi, \acts) \mapsto
            \actsT[\Phi]
            \quad
            \textrm{where}
            \quad
            g \actsT[\Phi] x = \Phi(g \acts \Phi^{-1}(x)).
        \end{equation}
        Two $H$-actions $\acts$, $\actsp$ define isomorphic
        $H$-covariant strong Morita equivalence bimodules if and only
        if $\acts$ and $\actsp$ lie in the same $\Aut(\BEA)$-orbit.
    \item For $c \in \group{U}(\zentrum(\mathcal{B}))$ we have
        $\actsT[\Phi_c] = \acts$ if and only if $c \in
        \group{U}(\zentrum(\mathcal{B}))^H$.
    \item Let $\twist{b} \in \group{U}(H, \mathcal{B})$. Then $\acts$
        and $\actsT[\twist{b}]$ define isomorphic $H$-covariant strong
        Morita equivalence bimodules if and only if $\twist{b} =
        \hat{c}$ for some $c \in \group{U}(\zentrum(\mathcal{B}))$
        where $\hat{c}(g) = c (g \acts c^{-1})$.
    \end{compactenum}
\end{lemma}
\begin{proof}
    The first part is obvious as any bimodule homomorphism can be
    written as the left multiplication with a unique \emph{central}
    element of $\mathcal{B}$. Then the isometry condition implies
    immediately $c^* = c^{-1}$. For the second part a straightforward
    computation shows that $\actsT[\Phi]$ is indeed a compatible
    $H$-action again. The remaining statements are obvious.  The third
    part is clear. The fourth part is then a simple consequence as
    \[
    g \actsT[\Phi_c] x
    =
    c \cdot (g \acts (c^{-1} \cdot x))
    =
    (c (g_\sweedler{1} \acts c^{-1})) \cdot (g_\sweedler{2} \acts x)
    =
    \hat{c}(g_\sweedler{1}) \cdot (g_\sweedler{2} \acts x).
    \]
\end{proof}

The last ingredient we need to describe the kernel of the groupoid
morphism \eqref{eq:PicHtoPic} is the following statement:
\begin{lemma}
    \label{lemma:UHBActsOnActs}
    Let $\BEA$ be an $H$-covariant strong Morita equivalence bimodule.
    \begin{compactenum}
    \item The group $\group{U}(H, \mathcal{B})$ acts transitively and
        freely of the set of all $H$-actions which make $\BEA$ a
        $H$-covariant strong Morita equivalence bimodule by
        \begin{equation}
            \label{eq:UHAActionOnActs}
            (\twist{b}, \acts)
            \mapsto
            \actsT[\twist{b}].
        \end{equation}
    \item The action of $\group{U}(H, \mathcal{B})$ and
        $\group{U}(\zentrum(\mathcal{B}))$ are compatible with the
        group morphism $\group{U}(\zentrum(\mathcal{B}))
        \longrightarrow \group{U}(H, \mathcal{B})$, i.e.~we have
        \begin{equation}
            \label{eq:UHBActsUZBNice}
            \actsT[\hat{c}] = \actsT[\Phi_c]
        \end{equation}
        for all $c \in \group{U}(\zentrum(\mathcal{B}))$.
    \end{compactenum}
\end{lemma}
\begin{proof}
    For the first part we already know that $\group{U}(H,
    \mathcal{B})$ parametrizes the $H$-actions in a one-to-one
    correspondence. Thus we only have to show that
    \eqref{eq:UHAActionOnActs} with $\actsT[\twist{b}]$ as in
    \eqref{eq:actsbDef} defines a group action. Let $\twist{b},
    \tilde{\twist{b}} \in \group{U}(H, \mathcal{B})$ be given. Then a
    straightforward computation shows $\actsT[\twist{b} *
    \tilde{\twist{b}}] = (\actsT[\tilde{\twist{b}}])^{\twist{b}}$
    and $\actsT[\twist{e}] = \acts$, whence the first part follows.
    The second is obvious from the preceding lemma.
\end{proof}

According to Lemma~\ref{lemma:AutBEAACtsOnActs}, part~\textit{iv.)}
the interesting twists $\actsT[\twist{b}]$ of the action $\acts$ are
described by the quotient $\group{U}_0(H, \mathcal{B}) = \group{U}(H,
\mathcal{B}) \big/ \widehat{\group{U}(\mathcal{Z}(\mathcal{B}))}$, see
also \eqref{eq:UNullDef}. Thus we can determine the kernel of the
groupoid morphism completely.
\begin{theorem}
    \label{theorem:KernelPicHPic}
    For unital $^*$-algebras the kernel of the groupoid morphism
    \eqref{eq:PicHtoPic} can be described as follows: Let
    $\mathcal{A}$, $\mathcal{B}$ be unital $^*$-algebras with
    $^*$-action of $H$ such that $\StrPic(\mathcal{B}, \mathcal{A})
    \ne \emptyset$. Then we have the following alternatives:
    \begin{compactenum}
    \item $\StrPicH(\mathcal{B}, \mathcal{A}) = \emptyset$.
    \item $\StrPicH(\mathcal{B}, \mathcal{A}) \longrightarrow
        \image(\StrPicH(\mathcal{B}, \mathcal{A})) \subseteq
        \StrPic(\mathcal{B}, \mathcal{A})$ is a principal
        $\group{U}_0(H, \mathcal{B})$-bundle over the image
        $\image(\StrPicH(\mathcal{B}, \mathcal{A}))$, i.e.~the group
        $\group{U}_0(H, \mathcal{B})$ acts freely and transitively
        (from the left) on the fibers of the projection.
    \end{compactenum}
\end{theorem}

Thus it is of major importance to understand the group $\group{U}_0(H,
\mathcal{B})$ for a given $^*$-algebra $\mathcal{B}$. As shown in the
appendix, this group can be quite non-trivial. Note that the first
alternative in the theorem may well happen and note also that the
image in the second case may not exhaust the whole set
$\StrPic(\mathcal{B}, \mathcal{A})$, see
Example~\ref{example:Lifting}.

The symmetry of the relation `$H$-covariant strong Morita equivalence'
already suggests that if $\StrPicH(\mathcal{B}, \mathcal{A})$ is
non-empty then $\group{U}_0(H, \mathcal{B}) \cong \group{U}_0(H,
\mathcal{A})$. This is indeed the case and will be investigated more
systematically in Section~\ref{subsec:UHAUNullHAMoritaInvariant}.

Finally, note that the Theorem is literally the same for $\StrPicH$
and $\StrPic$ being replaced by $\starPicH$ and $\starPic$,
respectively, as we have never used the positivity of the inner
products. It is also valid in the ring-theoretic situation if one
replaces $\group{U}_0(H, \mathcal{B})$ by $\group{GL}_0(H,
\mathcal{B})$.

\begin{remark}
    \label{remark:PicHDeformed}
    One of our original motivations was to understand the covariant
    Morita theory for star products. Using the techniques developed in
    this section one would like to proceed analogously to
    \cite{bursztyn.waldmann:2004a} in order to understand the
    covariant strong Picard groupoid of deformed algebras in terms of
    their classical limits. We address these topics in a project
    together with Nikolai Neumaier
    \cite{jansen.neumaier.waldmann:2005a:pre}.
\end{remark}

%
%

\section{Morita invariants and actions of the Picard groupoid}
\label{subsec:ActionsInvariants}

We shall now use the Picard groupoid in the spirit of
\cite{waldmann:2004a} to obtain \textdef{Morita invariants} (most
of which are well-known) as arising from \textdef{actions} of
$\StrPicH$ (or $\starPicH$) on `something'. Here an $H$-covariant,
strong Morita invariant is a property $\mathcal{P}$ of $^*$-algebras
with a $^*$-action of $H$ such that if $\mathcal{A}$ has this
property then any algebra $\mathcal{B}$ which is $H$-covariantly and
strongly Morita equivalent to $\mathcal{A}$ has this property
$\mathcal{P}$ as well, see also \cite[Def.~18.4]{lam:1999a} for the
ring-theoretic definition.

From this point of view, the Picard groups are the most fundamental
Morita invariant as they arise from the Picard groupoid acting on
itself by multiplication.  Hence (as for any groupoid) the isotropy
groups are all isomorphic along an orbit.

%
%

\subsection{The representation theories}
\label{subsec:REpTheoriesInvariant}

The statements of
Corollary~\ref{corollary:RieffelFunctorsTensorProducts} can be
rephrased in the following way, specializing the discussion in
\cite{bursztyn.waldmann:2003a:pre, waldmann:2004a} to the
$H$-covariant situation. Up to natural unitary equivalence the Picard
groupoids `act' on the representation theories by Rieffel induction
\begin{equation}
    \label{eq:starPicActsMod}
    \textrm{``}\mathsf{R}: \starPicH(\mathcal{B}, \mathcal{A}) 
    \times
    \sMod[\mathcal{D}, H](\mathcal{A})
    \longrightarrow
    \sMod[\mathcal{D}, H](\mathcal{B})\textrm{''}
\end{equation}
and
\begin{equation}
    \label{eq:StrPicActsRep}
    \textrm{``}\mathsf{R}: \StrPicH(\mathcal{B}, \mathcal{A}) 
    \times
    \Rep[\mathcal{D}, H](\mathcal{A})
    \longrightarrow
    \Rep[\mathcal{D}, H](\mathcal{B})\textrm{''},
\end{equation}
where we have of course \emph{not} an honest action as the Rieffel
induction functor $\mathsf{R}_{\mathcal{E}}$ depends on $\mathcal{E}$
and not only on its class in $\starPicH$ (or $\StrPicH$, respectively)
and the action properties $\mathsf{R}_{\mathcal{F}} \circ
\mathsf{R}_{\mathcal{E}} \cong \mathsf{R}_{\mathcal{F} \tensB
  \mathcal{E}}$ and $\mathsf{R}_{\mathcal{A}} \cong \id$ are only
fulfilled up to a natural transformation.

Thus \eqref{eq:starPicActsMod} and \eqref{eq:StrPicActsRep} become
actions once we pass to unitary equivalence classes of $H$-covariant
$^*$-representations. Alternatively, one should view
\eqref{eq:starPicActsMod} and \eqref{eq:StrPicActsRep} as an action of
the Picard \emph{bigroupoids}, where we have not yet identified
isomorphic bimodules. We shall not give a precise definition of an
action of a bigroupoid on a collection of categories (though in
principle this could be done) but leave this as a suggestive picture.
In any case, this gives a conceptually clear picture why
$H$-covariantly strongly (or $^*$-) Morita equivalent $^*$-algebras
have equivalent $H$-covariant representation theories. Moreover, we
see that the Picard groups $\starPicH(\mathcal{A})$ and
$\StrPicH(\mathcal{A})$, respectively, act on the unitary equivalence
classes of $H$-covariant $^*$-representations. In more physical terms,
these are just the \textdef{super-selection rules} of the $^*$-algebra
$\mathcal{A}$.

%
%

\subsection{The $H$-invariant central elements}
\label{subsec:HinvCenterInvariant}

We consider unital $^*$-algebras in this subsection.  Clearly, any
$H$-equivariant $^*$-homomorphism $\Phi: \mathcal{A} \longrightarrow
\mathcal{B}$ restricts to a $^*$-homomorphism
$\mathcal{Z}(\mathcal{A})^H \longrightarrow
\mathcal{Z}(\mathcal{B})^H$ of the $H$-invariant central elements. In
particular, this gives a groupoid action of the isomorphism groupoid
\begin{equation}
    \label{eq:IsoActsZAH}
    \starIsoH(\mathcal{B}, \mathcal{A}) \times
    \mathcal{Z}(\mathcal{A})^H 
    \longrightarrow \mathcal{Z}(\mathcal{B})^H.
\end{equation}

We shall now extend this to an action of $\starPicH$ in the following
way: First we recall some standard results from Morita theory, see
e.g.~\cite{ara:1999a, bursztyn.waldmann:2004a,
  bursztyn.waldmann:2001a}. If $\BEA$ is a $^*$-equivalence bimodule
then for any central element $a \in \mathcal{Z}(\mathcal{A})$ there
exists a unique central element $h_{\mathcal{E}}(a) \in
\mathcal{Z}(\mathcal{B})$ such that
\begin{equation}
    \label{eq:hEaxxa}
    h_{\mathcal{E}}(a) \cdot x = x \cdot a
\end{equation}
for all $x \in \BEA$ and the map $h_{\mathcal{E}}:
\mathcal{Z}(\mathcal{A}) \longrightarrow \mathcal{Z}(\mathcal{B})$ is
a $^*$-isomorphism. Moreover, $h_{\mathcal{E}} = h_{\mathcal{E}'}$ if
$\BEA$ and $\BEpA$ are isomorphic $^*$-Morita equivalence bimodules
and we have $h_{\mathcal{F}} \circ h_{\mathcal{E}} = h_{\mathcal{F}
  \tensB \mathcal{E}}$ as well as $h_{\mathcal{A}} =
\id_{\mathcal{Z}(\mathcal{A})}$, see
\cite[Section 2.3]{bursztyn.waldmann:2004a} and
\cite[Prop.~7.6]{bursztyn.waldmann:2001a}. This can be rephrased as an
action of the $^*$-Picard groupoid on centers
\begin{equation}
    \label{eq:starPicActsOnCenters}
    h: \starPic(\mathcal{B}, \mathcal{A}) \times
    \mathcal{Z}(\mathcal{A})
    \ni ([\mathcal{E}], a) \mapsto h_{\mathcal{E}}(a) \in
    \mathcal{Z}(\mathcal{B})
\end{equation}
by $^*$-isomorphisms. In particular, centers are invariant as
$^*$-algebras under $^*$-Morita equivalence \cite{ara:1999a}.
\begin{lemma}
    \label{lemma:BEAHcovZAHintoZBH}
    Let $\BEA$ be an $H$-covariant $^*$-Morita equivalence bimodule.
    Then $h_{\mathcal{E}}$ restricts to a $^*$-isomorphism
    $h_{\mathcal{E}}: \mathcal{Z}(A)^H \longrightarrow
    \mathcal{Z}(\mathcal{B})^H$. Thus we have an action of the
    $H$-covariant $^*$-Picard groupoid on $H$-invariant central
    elements
    \begin{equation}
        \label{eq:starPicHActsZH}
        h: \starPicH(\mathcal{B}, \mathcal{A}) \times
        \mathcal{Z}(\mathcal{A})^H 
        \longrightarrow
        \mathcal{Z}(\mathcal{B})^H.
    \end{equation}
\end{lemma}
\begin{proof}
    We only have to check that $h_{\mathcal{E}}$ maps $H$-invariant
    elements to $H$-invariant ones. For $x \in \BEA$ we have
    \[
    (g \acts h_{\mathcal{E}}(a)) \cdot x
    =
    g_\sweedler{1} \acts 
    \left(h_{\mathcal{E}}(a) \cdot (S(g_\sweedler{2}) \acts x)\right)
    =
    g_\sweedler{1} \acts 
    \left((S(g_\sweedler{2}) \acts x) \cdot a \right)
    =
    \epsilon(g) x \cdot a
    =
    \epsilon(g) h_{\mathcal{E}}(a) \cdot x,
    \]
    since $a$ is invariant. By Remark~\ref{remark:InjectiveActions} we
    get $g \acts h_{\mathcal{E}}(a) = \epsilon(g) h_{\mathcal{E}}(a)$.
    Then the action properties for \eqref{eq:starPicHActsZH} follow
    immediately from those of \eqref{eq:starPicActsOnCenters} and
    \eqref{eq:NiceTrianglePicDiagram}.
\end{proof}
\begin{corollary}
    \label{corollary:ZAHInvariant}
    The $H$-covariant $^*$-Picard group $\starPicH(\mathcal{A})$ acts
    on $\mathcal{Z}(\mathcal{A})^H$ by $^*$-isomorphisms whence
    $\mathcal{Z}(\mathcal{A})^H$ as $\starPicH(\mathcal{A})$-space is
    invariant under $H$-covariant $^*$-Morita equivalence.
\end{corollary}
Moreover, we have a compatibility between the canonical groupoid
action \eqref{eq:IsoActsZAH} and the action $h$, adapting
\cite[Prop.~2.4]{bursztyn.waldmann:2004a} to this situation:
\begin{lemma}
    \label{lemma:IsoStarPicActionOnZH}
    The actions \eqref{eq:IsoActsZAH} and \eqref{eq:starPicHActsZH}
    are compatible in the sense that the diagram
    \begin{equation}
        \label{eq:IsoPicActsOnZAH}
        \bfig
        \Dtriangle/>`>`<-/<1200,200>%
        [\starIsoH(\mathcal{B},\mathcal{A})\times\mathcal{Z}(\mathcal{A})^H%
        `\mathcal{Z}(\mathcal{B})^H%
        `\starPicH(\mathcal{B},\mathcal{A})\times\mathcal{Z}(\mathcal{A})^H%
        ;\ell\times\id``h]
        \efig
    \end{equation}
    commutes.
\end{lemma}

In general, the center $\mathcal{Z}(\mathcal{A})$ needs not to be
preserved by the action of $H$. However, if $H$ is
\emph{cocommutative} then this is the case whence
$\mathcal{Z}(\mathcal{A})$ inherits a $^*$-action of $H$. In this
case, the action $h$ of $\starPic$ on centers
\eqref{eq:starPicActsOnCenters} restricts to an action, also denoted
by $h$, of $\starPicH$ on the centers
\begin{equation}
    \label{eq:starPicHActsCentersCocom}
    h: \starPicH(\mathcal{B}, \mathcal{A}) \times
    \mathcal{Z}(\mathcal{A}) \longrightarrow \mathcal{Z}(\mathcal{B})
\end{equation}
by \emph{$H$-equivariant} $^*$-isomorphisms as a simple computation
shows. Hence we have
\begin{lemma}
    \label{lemma:HcocomZAHspaceInv}
    Let $H$ be cocommutative. Then $\starPicH$ acts on the centers by
    $H$-equivariant $^*$-iso\-mor\-phisms whence
    $\mathcal{Z}(\mathcal{A})$ as a $\starPicH(\mathcal{A})$-space is
    invariant under $H$-covariant $^*$-Morita equivalence.
\end{lemma}

%
%

\subsection{Equivariant $K$-theory}
\label{subsec:Ktheory}

Again we shall restrict ourselves to unital $^*$-algebras in this
subsection for simplicity. There are many notions of equivariant
$K$-theory, we shall use a rather naive definition taking care of the
inner products as well.

We consider $H$-covariant pre-Hilbert (right) $\mathcal{A}$-modules
$\PA$ with the following additional properties: The inner product
$\SPPA{\cdot,\cdot}$ is \textdef{strongly non-degenerate}, i.e.~the
map $x \mapsto \SPPA{x, \cdot} \in \Hom_{\mathcal{A}}(\PA,
\mathcal{A})$ is bijective. Moreover, we want $\PA$ to be
\textdef{finitely generated} and \textdef{projective}. The subcategory
of all $H$-covariant pre-Hilbert $\mathcal{A}$-modules with these two
additional properties is denoted by $\StrProjH(\mathcal{A})$, where
the morphisms are adjointable module morphisms as before. By
$\StrProjHClass(\mathcal{A})$ we denote the set of isometric
isomorphism classes of $\StrProjH(\mathcal{A})$. Then
$\StrProjHClass(\mathcal{A})$ becomes an abelian semigroup where the
addition $\oplus$ is induced by the direct orthogonal sum of elements
in $\StrProjH(\mathcal{A})$. The \textdef{$H$-equivariant strong
  $K_0$-group} $\StrKH(\mathcal{A})$ of $\mathcal{A}$ is then by
definition the Grothendieck group associated to
$\StrProjHClass(\mathcal{A})$. Similarly, dropping the complete
positivity of the inner product (but keeping the strong
non-degeneracy) we obtain $^*$-versions $\starProjH(\mathcal{A})$,
$\starProjHClass(\mathcal{A})$ and $\starKH(\mathcal{A})$,
respectively.

A $H$-covariant pre-Hilbert module $\PA$ is in
$\StrProjH(\mathcal{A})$ if and only if there exist $x_i, y_i \in \PA$
with $i = 1, \ldots, n$ such that
\begin{equation}
    \label{eq:HermitianDualBasis}
    x = \sum\nolimits_i x_i \cdot \SPPA{y_i, x}
\end{equation}
for all $x \in \PA$.  This is an easy adaption of the dual basis lemma
for projective modules, see e.g.~\cite[Lem.~2.9]{lam:1999a}. We shall
call such vectors $x_i, y_i$ a \textdef{Hermitian dual basis}.
Then we have the following lemma:
\begin{lemma}
    \label{lemma:ProjectiveHermitianDualBasis}
    Let $\PB \in \StrProjH(\mathcal{B})$ and let $\BEA$ be a
    $H$-covariant strong Morita equivalence bimodule. Then $\BEA$ as
    right $\mathcal{A}$-module is in $\StrProjH(\mathcal{A})$ and
    $\PB \tensM[\mathcal{B}] \BEA \in \StrProjH(\mathcal{A})$, too.
\end{lemma}
\begin{proof}
    The first statement is well-known and follows directly from the
    fullness of $\BSPE{\cdot,\cdot}$ and the compatibility
    \eqref{eq:SPCompatible}. For the second statement, let $\{x_i,
    y_i\}_{i=1, \ldots, n}$ be a Hermitian dual basis for $\PB$ and
    let $\{\xi_\alpha, \eta_\alpha\}_{\alpha = 1, \ldots, m}$ be a
    Hermitian dual basis for $\BEA$ viewed as right
    $\mathcal{A}$-module. Then $\{x_i \tensor[\mathcal{B}] \xi_\alpha,
    y_i \tensor[\mathcal{B}] \eta_\alpha\}_{i, \alpha}$ is easily
    shown to be a Hermitian dual basis for $\PB \tensM[\mathcal{B}]
    \BEA$.  In particular, the inner product on $\PB
    \otimes_{\mathcal{B}} \BEA$ is already non-degenerate whence the
    usual quotient procedure for $\tensM$ is not needed here.
\end{proof}

From this and the associativity properties of $\tensM$ and $\tensB$ as
in Proposition~\ref{proposition:TensorTildeNett} we immediately obtain
the following result:
\begin{proposition}
    \label{proposition:PicActsOnK}
    The $H$-covariant strong Picard groupoid acts on $\StrProjHClass$
    by semi-group isomorphisms from the right, i.e.
    \begin{equation}
        \label{eq:SAction}
        \mathsf{S}:
        \StrProjHClass(\mathcal{B}) 
        \times \StrPicH(\mathcal{B}, \mathcal{A})
        \ni ([\PB], [\BEA]) \mapsto
        [\mathsf{S}_{\mathcal{E}}(\PB)] 
        = [\PA \tensM[\mathcal{B}] \BEA] \in 
        \StrProjHClass(\mathcal{A}),
    \end{equation}
    and hence it also acts on the $H$-equivariant strong $K_0$-groups
    by group isomorphisms
    \begin{equation}
        \label{eq:SActionK}
        \mathsf{S}:
        \StrKH(\mathcal{B}) \times 
        \StrPicH(\mathcal{B}, \mathcal{A})
        \longrightarrow
        \StrKH(\mathcal{A}).
    \end{equation}
    The analogous result holds for $\starPicH$, $\starProjHClass$ and
    $\starKH$.
\end{proposition}
\begin{corollary}
    \label{corollary:KNullInvariant}
    The $H$-covariant strong Picard group $\StrPicH(\mathcal{A})$ acts
    on $\StrKH(\mathcal{A})$ by group automorphisms and
    $\StrKH(\mathcal{A})$ is invariant as
    $\StrPicH(\mathcal{A})$-space under $H$-covariant strong Morita
    equivalence.
\end{corollary}

Note that this result corresponds to the `action' by Rieffel induction
$\mathsf{R}$ on representation theories, where we have replaced the
action from the left via $\mathsf{R}$ by an action from the right via
the change of base ring functors $\mathsf{S}$.

Again the $H$-equivariant isomorphisms $\starIsoH$ act on
$\StrProjHClass$ and hence on $\StrKH$ as well and the above actions
\eqref{eq:SAction} and \eqref{eq:SActionK} restrict to this via the
groupoid morphism $\ell$ from
Proposition~\ref{proposition:IsoToPicMorphism}.

%
%

\subsection{The lattice $\Lattice[\mathcal{D}, H](\mathcal{A})$}
\label{subsec:IdealMoritaInvariant}

Let $\mathcal{D}$ be admissible and all other $^*$-algebras are
idempotent and non-degenerate as before. Then we can act with
$\StrPicH$ on the lattices of $(\mathcal{D}, H)$-closed ideals by the
following construction. Let $\BEA$ be an $H$-covariant strong Morita
equivalence bimodule and let $\mathcal{J} \subseteq \mathcal{A}$ be a
subset. Then we define
\begin{equation}
    \label{eq:PhiELatticeDef}
    \Phi_{\mathcal{E}}(\mathcal{J}) = 
    \left\{
        b \in \mathcal{B} 
        \;\left| \;
            \SPEA{x, b \cdot y} \in \mathcal{J}
            \; \textrm{for all} \; x, y \in \BEA
        \right.
    \right\}.
\end{equation}
We have the following properties of the map $\Phi_{\mathcal{E}}$
generalizing the results of \cite{bursztyn.waldmann:2001b}:
\begin{lemma}
    \label{lemma:PhiEActsLattice}
    Let $\BEA$ be an $H$-covariant strong Morita equivalence
    bimodule and let $\mathcal{D}$ be admissible.
    \begin{compactenum}
    \item If $\mathcal{J} = \ker\pi$ for $(\HD, \pi) \in
        \Rep[\mathcal{D}, H](\mathcal{A})$ then
        $\Phi_{\mathcal{E}}(\mathcal{J}) = \ker
        \mathsf{R}_{\mathcal{E}}\pi$ whence in particular
        $\Phi_{\mathcal{E}}(\mathcal{J}) \in \Lattice[\mathcal{D},
        H](\mathcal{B})$ for any $\mathcal{J} \in
        \Lattice[\mathcal{D}, H](\mathcal{A})$.
    \item If $\BEpA$ is another $H$-covariant strong Morita
        equivalence bimodule isomorphic to $\BEA$ then
        $\Phi_{\mathcal{E}} = \Phi_{\mathcal{E}'}$.
    \item If $\CFB$ is another $H$-covariant strong Morita equivalence
        bimodule then $\Phi_{\mathcal{F}} \circ \Phi_{\mathcal{E}} =
        \Phi_{\mathcal{F} \tensB \mathcal{E}}$ and $\Phi_{\mathcal{A}}
        = \id_{\Lattice[\mathcal{D}, H](\mathcal{A})}$.
    \end{compactenum}
\end{lemma}
\begin{proof}
    The first part is analogous to
    \cite[Prop.~5.1]{bursztyn.waldmann:2001b}. The second part follows
    as $\mathsf{R}_{\mathcal{E}}(\pi)$ and
    $\mathsf{R}_{\mathcal{E}'}(\pi)$ are unitarily equivalent
    $^*$-representations which therefor have the same kernel. The same
    Rieffel induction argument can be used for the third part since we
    can restrict to strongly non-degenerate $^*$-representations by
    Lemma~\ref{lemma:Lattice}.
\end{proof}

From this lemma we easily conclude the following statement
generalizing Rieffel's correspondence from the theory of
$C^*$-algebras, see e.g.~\cite[Thm.~3.24]{raeburn.williams:1998a}, as
well as \cite[Thm.~5.4]{bursztyn.waldmann:2001b}:
\begin{theorem}
    \label{theorem:LatticeAction}
    Let $\mathcal{D}$ be admissible. Then the map
    \begin{equation}
        \label{eq:PhiActionPicLattice}
        \Phi: \StrPicH(\mathcal{B}, \mathcal{A}) 
        \times \Lattice[\mathcal{D}, H](\mathcal{A}) 
        \longrightarrow
        \Lattice[\mathcal{D}, H](\mathcal{B})
    \end{equation}
    defines an action of the $H$-covariant strong Picard groupoid on
    the lattices of $(\mathcal{D}, H)$-closed ideals by lattice
    isomorphisms.
\end{theorem}
\begin{proof}
    The only thing to be checked is that $\Phi_{\mathcal{E}}$ is a
    lattice homomorphism as the well-definedness and the action
    properties follow from Lemma~\ref{lemma:PhiEActsLattice}. Clearly
    we have $\Phi_{\mathcal{E}}(\mathcal{I}) \le
    \Phi_{\mathcal{E}}(\mathcal{J})$ for $\mathcal{I} \le \mathcal{J}$
    and $\Phi_{\mathcal{E}}(\mathcal{I} \wedge \mathcal{J}) =
    \Phi_{\mathcal{E}}(\mathcal{I}) \wedge
    \Phi_{\mathcal{E}}(\mathcal{J})$. From these properties and the
    bijectivity of $\Phi_{\mathcal{E}}$ it follows that
    $\Phi_{\mathcal{E}}(\mathcal{I} \vee \mathcal{J}) =
    \Phi_{\mathcal{E}}(\mathcal{I}) \vee
    \Phi_{\mathcal{E}}(\mathcal{J})$.
\end{proof}
\begin{corollary}
    \label{corollary:LatticeInv}
    Let $\mathcal{D}$ be admissible. Then $\StrPicH(\mathcal{A})$ acts
    on the lattice $\Lattice[\mathcal{D}, H](\mathcal{A})$ by lattice
    automorphisms and $\Lattice[\mathcal{D}, H](\mathcal{A})$ as
    $\StrPicH(\mathcal{A})$-space is invariant under $H$-covariant
    strong Morita equivalence.
\end{corollary}

%
%

\subsection{The groups $\group{U}(H, \mathcal{A})$ and $\group{U}_0(H,
  \mathcal{A})$}
\label{subsec:UHAUNullHAMoritaInvariant}

Also in this subsection the $^*$-algebras are required to be unital.
In the characterization of the kernel of the canonical groupoid
morphism $\StrPicH \longrightarrow \StrPic$ as well as for $\starPicH
\longrightarrow \starPic$ the groups $\group{U}(H, \mathcal{A})$ and
$\group{U}_0(H, \mathcal{A})$ play the dominant role which already
suggests that they are a Morita invariant.

As we have outlined in the appendix, the $H$-equivariant
$^*$-isomorphisms $\starIsoH$ act not only on $\group{U}(H,
\mathcal{A})$ and $\group{U}_0(H, \mathcal{A})$ in a canonical way but
also on the whole exact sequence \eqref{eq:TheUExactSequenz}. We shall
now extend this to an action of $\starPicH$ extending thereby the
action $h$ of $\starPicH$ on the centers.
\begin{lemma}
    \label{lemma:hPicUHA}
    Let $\BEA$ be an $H$-covariant $^*$-Morita equivalence bimodule and
    let $\twist{a} \in \Hom_{\ring{C}}(H, \mathcal{A})$.
    \begin{compactenum}
    \item The definition
        \begin{equation}
            \label{eq:actsSubaDef}
            g \actsS[\twist{a}] x = 
            (g_\sweedler{1} \acts x) \cdot \twist{a}(g_\sweedler{2})
        \end{equation}
        gives another compatible $H$-action on $\BEA$ such that
        $(\BEA, \actsS[\twist{a}])$ is an $H$-covariant $^*$-Morita
        equivalence bimodule if and only if $\twist{a} \in
        \group{U}(H, \mathcal{A})$ and any such action is of this form
        for a uniquely determined $\twist{a} \in \group{U}(H,
        \mathcal{A})$.
    \item The group $\group{U}(H, \mathcal{A})$ acts freely and
        transitively from the right on the set of all compatible
        $H$-actions on $\BEA$ by $(\acts, \twist{a}) \mapsto
        \actsS[\twist{a}]$.
    \item For $\twist{b} \in \group{U}(H, \mathcal{B})$ and $\twist{a}
        \in \group{U}(H, \mathcal{A})$ we have
        $(\actsS[\twist{a}])^\twist{b} =
        (\actsT[\twist{b}])_\twist{a}$ and there exists a unique
        $h_{\mathcal{E}} (\twist{a}) \in \group{U}(H, \mathcal{B})$
        such that
        \begin{equation}
            \label{eq:actsSaactsThEa}
            \actsS[\twist{a}] = \actsT[h_{\mathcal{E}}(\twist{a})].
        \end{equation}
    \item The map
        \begin{equation}
            \label{eq:hEUHAtoUHB}
            h_{\mathcal{E}}:
            \group{U}(H, \mathcal{A}) 
            \ni \twist{a} \mapsto h_{\mathcal{E}}(\twist{a}) \in 
            \group{U}(H, \mathcal{B})
        \end{equation}
        is a group isomorphism.
    \end{compactenum}
\end{lemma}
\begin{proof}
    The first part is lengthy computation but completely analogous to
    Proposition~\ref{proposition:twistbIStwist}. The second part is in
    the same spirit as Proposition~\ref{proposition:twistbIStwist} as
    well. For the third part we have
    \[
    g (\actsS[\twist{a}])^\twist{b} x
    =
    \twist{b}(g_\sweedler{1}) \cdot 
    (g_\sweedler{2} \actsS[\twist{a}] x)
    =
    \twist{b}(g_\sweedler{1}) \cdot 
    (g_\sweedler{2} \acts x) \cdot 
    \twist{a}(g_\sweedler{3})
    =
    (g_\sweedler{1} \actsT[\twist{b}] x) \cdot
    \twist{a}(g_\sweedler{2})
    =
    g (\actsT[\twist{b}])_\twist{a} x,
    \]
    since $\BEA$ is a bimodule. Then the remaining statements are
    general facts on commuting free and transitive group actions.
\end{proof}

The next lemma investigates the dependence of the isomorphism
$h_{\mathcal{E}}$ on the bimodule $\mathcal{E}$:
\begin{lemma}
    \label{lemma:hEClassE}
    Let $\BEA$ and $\BEpA$ be isomorphic $H$-covariant $^*$-Morita
    equivalence bimodules. Then $h_{\mathcal{E}} = h_{\mathcal{E}'}$.
\end{lemma}
\begin{proof}
    Let $U: \mathcal{E} \longrightarrow \mathcal{E}'$ be an
    isomorphism. Then on one hand
    \[
    U\left(g \actsT[h_{\mathcal{E}}(\twist{a})] x\right)
    =
    U\left(
        h_{\mathcal{E}}(\twist{a})(g_\sweedler{1}) \cdot
        (g_\sweedler{2} \acts x)
    \right)
    =
    h_{\mathcal{E}}(\twist{a})(g_\sweedler{1}) \cdot 
    (g_\sweedler{2} \actsp U(x))
    =
    g (\actsp)^{h_{\mathcal{E}}(\twist{a})} U(x)
    \]
    and on the other hand
    \[
    U\left(g \actsT[h_{\mathcal{E}}(\twist{a})] x\right)
    =
    U(g \actsS[\twist{a}] x)
    =
    U\left(
        (g_\sweedler{1} \acts x)\cdot \twist{a}(g_\sweedler{2})
    \right)
    =
    (g_\sweedler{1} \actsp U(x)) \cdot \twist{a}(g_\sweedler{2})
    =
    g (\actsp)^{h_{\mathcal{E}'}(\twist{a})} U(x).
    \]
    This implies $h_{\mathcal{E}}(\twist{a}) =
    h_{\mathcal{E}'}(\twist{a})$ by the uniqueness from
    Lemma~\ref{lemma:hPicUHA}.
\end{proof}
\begin{lemma}
    \label{lemma:hEhFhEFHellPhi}
    Let $\BEA$ and $\CFB$ be $H$-covariant $^*$-Morita equivalence
    bimodules. Then
    \begin{equation}
        \label{eq:hFhEhFE}
        h_{\mathcal{F}} \circ h_{\mathcal{E}} 
        = h_{\mathcal{F} \tensB \mathcal{E}}
        \quad
        \textrm{and}
        \quad
        h_{\mathcal{A}} = \id_{\group{U}(H, \mathcal{A})}.
    \end{equation}
    For $\Phi \in \starIsoH(\mathcal{A})$ we have
    \begin{equation}
        \label{eq:hellPhiPhi}
        h_{\ell(\Phi)} = \Phi_*,
    \end{equation}
    where $\Phi_*: \group{U}(H, \mathcal{A}) \longrightarrow
    \group{U}(H, \mathcal{B})$ as in
    Proposition~\ref{proposition:GLUFunktoriell}.
\end{lemma}
\begin{proof}
    Let $x \in \mathcal{F}$, $\phi \in \mathcal{E}$ and $\twist{a} \in
    \group{U}(H, \mathcal{A})$. Then
    \[
    \begin{split}
      g \actsT[h_{\mathcal{F} \tensB \mathcal{E}}(\twist{a})] 
      (x \otimes \phi)
      &=
      (g_\sweedler{1} \acts (x \otimes \phi)) \cdot
      \twist{a}(g_\sweedler{2}) \\
      &=
      (g_\sweedler{1} \acts x) \otimes 
      (g_\sweedler{2} \acts \phi \cdot \twist{a}(g_\sweedler{3})) \\
      &=
      (g_\sweedler{1} \acts x) \otimes 
      \left(
          h_{\mathcal{E}}(\twist{a})(g_\sweedler{2}) 
          \cdot 
          (g_\sweedler{3} \acts \phi) 
      \right) \\
      &=
      \left(
          (g_\sweedler{1} \acts x) \cdot
          h_{\mathcal{E}}(\twist{a})(g_\sweedler{2})
      \right)
      \otimes 
      (g_\sweedler{3} \acts \phi) \\
      &=
      \left(
          h_{\mathcal{F}}(h_{\mathcal{E}}(\twist{a}))(g_\sweedler{1})
          \cdot (g_\sweedler{2} \acts x)
      \right)
      \otimes
      (g_\sweedler{3} \acts \phi) \\
      &=
      h_{\mathcal{F}}(h_{\mathcal{E}}(\twist{a}))(g_\sweedler{1})
      \cdot
      (g_\sweedler{2} \acts (x \otimes \phi)) \\
      &=
      g \actsT[h_{\mathcal{F}}(h_{\mathcal{E}}(\twist{a}))]
      (x \otimes \phi)
    \end{split}
    \]
    proves the first part. The second statement in \eqref{eq:hFhEhFE}
    is trivial using the `module condition' for $\twist{a}$. The last
    statement \eqref{eq:hellPhiPhi} is also a straightforward
    computation.
\end{proof}

Collecting these results, we get a generalization of the action $h$ of
the Picard groupoid on centers and a generalization of the action of
$\starIsoH$ on the exact sequence \eqref{eq:UZAHAHANullCommutes}.
\begin{theorem}
    \label{theorem:PicActsOnUHA}
    The map
    \begin{equation}
        \label{eq:hActionPicUHA}
        h: \starPicH(\mathcal{B}, \mathcal{A}) 
        \times \group{U}(H, \mathcal{A})
        \ni ([\mathcal{E}], \twist{a}) 
        \mapsto
        h_{\mathcal{E}}(\twist{a}) \in
        \group{U}(H, \mathcal{B})
    \end{equation}
    determines an action of $\starPicH$ on the exact sequence
    \eqref{eq:TheUExactSequenz}, i.e.
    \begin{equation}
        \label{eq:CoolesDiagramPicUHA}
        \bfig
        \hSquares|rrrrrrr|/>`>``>`>`>`>/%
        [1`\group{U}(\zentrum(\mathcal{A}))^H%
        `\group{U}(\zentrum(\mathcal{A}))%
        `1%
        `\group{U}(\zentrum(\mathcal{B}))^H%
        `\group{U}(\zentrum(\mathcal{B}))%
        ;```h_{\mathcal{E}}`h_{\mathcal{E}}``]
        \morphism(1550,0)<300,0>[`;]
        \morphism(1550,500)<300,0>[`;]
        \hSquares(2050,0)|rrrrrrr|/>`>`>`>``>`>/%
        [\group{U}(H, \mathcal{A})%
        `\group{U}_0(H,\mathcal{A})%
        `1%
        `\group{U}(H,\mathcal{B})%
        `\group{U}_0(H,\mathcal{B})%
        `1%
        ;``h_{\mathcal{E}}`h_{\mathcal{E}}```]
        \efig.
    \end{equation}
    commutes and all $h_{\mathcal{E}}$ are group isomorphisms.
    Moreover, this groupoid action is compatible with the groupoid
    morphism $\ell: \starIsoH \longrightarrow \starPicH$ and the
    canonical action of $\starIsoH$ on the exact sequence as in
    Corollary~\ref{corollary:IsoActsOfSequence}. 
\end{theorem}
\begin{proof}
    The only thing left to show is the commutativity of the box in the
    middle of \eqref{eq:CoolesDiagramPicUHA} since then the last
    vertical arrow is defined in such a way, that
    \eqref{eq:CoolesDiagramPicUHA} commutes.  Thus let $c \in
    \group{U}(\mathcal{Z}(\mathcal{A}))$ be given. Then for $x \in
    \BEA$ we have
    \[
    \begin{split}
        g \actsT[h_{\mathcal{E}}(\hat{c})] x
        &=
        (g_\sweedler{1} \acts x) \cdot \hat{c}(g_\sweedler{2}) \\
        &=
        g_\sweedler{1} \acts x \cdot (c g_\sweedler{2} \acts c^{-1}) \\
        &=
        h_{\mathcal{E}}(c) \cdot 
        (g_\sweedler{1} \acts (x \cdot c^{-1})) \\
        &=
        h_{\mathcal{E}}(c) \cdot (g_\sweedler{1} \acts
        (h_{\mathcal{E}}(c)^{-1} \cdot x)) \\
        &=
        \widehat{h_{\mathcal{E}}(c)}(g_\sweedler{1}) 
        \cdot g_\sweedler{2} \acts x \\
        &=
        g \actsT[\widehat{h_{\mathcal{E}}(c)}] x,
    \end{split}
    \]
    whence $h_{\mathcal{E}}(\hat{c}) = \widehat{h_{\mathcal{E}}(c)}$.
\end{proof}
We leave it to the reader to draw the appropriate big commutative
diagram expressing all compatibilities relating $\ell$ and $h$ stated
in this theorem.
\begin{corollary}
    \label{corollary:UHAInvariant}
    The $H$-covariant $^*$-Picard group $\starPicH(\mathcal{A})$ acts
    on the exact sequence \eqref{eq:TheUExactSequenz} by isomorphisms
    whence \eqref{eq:TheUExactSequenz} as a
    $\starPicH(\mathcal{A})$-space is invariant under $H$-covariant
    $^*$-Morita equivalence. In particular, each of the groups
    $\group{U}(\mathcal{Z}(\mathcal{A}))^H$,
    $\group{U}(\mathcal{Z}(\mathcal{A}))$, $\group{U}(H,
    \mathcal{A})$, and $\group{U}_0(H, \mathcal{A})$ carries a
    canonical $\starPicH(\mathcal{A})$-action by group automorphisms.
    They are invariant under $H$-covariant $^*$-Morita equivalence.
\end{corollary}

We can interpret the result of Lemma~\ref{lemma:hEhFhEFHellPhi} also
in another way. According to Theorem~\ref{theorem:KernelPicHPic} the
group $\group{U}_0(H, \mathcal{B})$ acts on $\StrPicH(\mathcal{B},
\mathcal{A})$ freely by twisting the $H$-action
\begin{equation}
    \label{eq:UnullBHTwistsBEA}
    [\twist{b}] \cdot \left[\BEA, \acts\right] 
    = 
    \left[\BEA, \actsT[\twist{b}]\right].
\end{equation}
Similarly, $\group{U}_0(H, \mathcal{A})$ acts from the right by
\begin{equation}
    \label{eq:UNullHAACtsBEA}
    \left[\BEA, \acts\right] \cdot [\twist{a}]
    =
    \left[\BEA, \actsS[\twist{a}]\right].
\end{equation}
Then from the proof of Lemma~\ref{lemma:hEhFhEFHellPhi} we see that
we have the following compatibilities between these two actions and
the tensor product of bimodules, namely
\begin{equation}
    \label{eq:bFEbFE}
    [\twist{c}] \cdot 
    \left(\left[\CFB\right] \tensB \left[\BEA\right]\right)
    =
    \left([\twist{c}] \cdot \left[\CFB\right]\right)
    \tensB \left[\BEA\right],
\end{equation}
\begin{equation}
    \label{eq:bEEhEinvb}
    [\twist{b}] \cdot \left[\BEA\right]
    =
    \left[\BEA\right] \cdot 
    \left[h_{\mathcal{E}}^{-1}(\twist{b})\right],
\end{equation}
\begin{equation}
    \label{eq:FbEFbE}
    \left(\left[\CFB\right] \cdot[\twist{b}]\right)
    \tensB \left[\BEA\right]
    =
    \left[\CFB\right] \tensB
    \left([\twist{b}] \cdot \left[\BEA\right]\right)
\end{equation}
and
\begin{equation}
    \label{eq:FEaFEa}
    \left(\left[\CFB\right] \tensB \left[\BEA\right]\right)
    \cdot [\twist{a}]
    =
    \left[\CFB\right] \tensB
    \left(\left[\BEA\right] \cdot [\twist{a}] \right).
\end{equation}
for $\twist{c} \in \group{U}(H,\mathcal{C})$, $\twist{b} \in
\group{U}(H, \mathcal{B})$ and $\twist{a} \in \group{U}(H,
\mathcal{A})$. From this we conclude the following statement:
\begin{proposition}
    \label{proposition:UnullInStrPicH}
    The map
    \begin{equation}
        \label{eq:UnullHAToStrPicH}
        \group{U}_0(H, \mathcal{A}) \ni [\twist{a}]
        \mapsto
        [\twist{a}] \cdot \left[\AAA\right] 
        = \left[\AAA, \actsT[\twist{a}]\right]
        \in 
        \StrPicH(\mathcal{A})
    \end{equation}
    is an injective group homomorphism such that
    \begin{equation}
        \label{eq:UnullPicHPicExact}
        1 \longrightarrow
        \group{U}_0(H, \mathcal{A})
        \longrightarrow
        \StrPicH(\mathcal{A})
        \longrightarrow
        \StrPic(\mathcal{A})
    \end{equation}
    is exact.
\end{proposition}
\begin{proof}
    It follows from a straightforward computation using
    \eqref{eq:bFEbFE}, \eqref{eq:bEEhEinvb}, \eqref{eq:FbEFbE}, and
    \eqref{eq:FEaFEa} that \eqref{eq:UnullHAToStrPicH} is a
    well-defined group homomorphism. The exactness of
    \eqref{eq:UnullPicHPicExact} is then a consequence of
    Theorem~\ref{theorem:KernelPicHPic}.
\end{proof}

Though this observation helps to understand the $H$-covariant strong
Picard group one should not overestimate its importance as the group
morphism $\StrPicH(\mathcal{A}) \longrightarrow \StrPic(\mathcal{A})$
is only in the very simplest cases surjective.

%
%

\section{Crossed products}
\label{sec:CrossedProducts}

In this section we shall investigate the crossed product algebras
$\mathcal{A} \rtimes H$ and relate their Picard groupoids with the
$H$-covariant Picard groupoids of the underlying algebras
$\mathcal{A}$.

%
%

\subsection{Definitions and preliminary results}
\label{subsec:DefPrelimResultCrossed}

Let $\mathcal{A}$ be a $^*$-algebra over $\ring{C}$ with a $^*$-action
of a Hopf $^*$-algebra $H$. Recall that the crossed product
$^*$-algebra $\mathcal{A} \rtimes H$ is $\mathcal{A} \otimes H$ as a
$\ring{C}$-module with multiplication defined by
\begin{equation}
    \label{eq:CrossedProductMultiplication}
    (a \otimes g) (b \otimes h) 
    = (a (g_\sweedler{1} \acts b)) \otimes g_\sweedler{2} h
\end{equation}
and $^*$-involution
\begin{equation}
    \label{eq:CrossedProductInvolution}
    (a \otimes g)^* 
    = g_\sweedler{1}^* \acts a^* \otimes g_\sweedler{2}^*.
\end{equation}
Then it is well-known that $\mathcal{A} \rtimes H$ is a $^*$-algebra
over $\ring{C}$, sometimes also called the smash product of
$\mathcal{A}$ and $H$, see e.g.~\cite{majid:1995a, kassel:1995a} for
this and more general crossed product constructions and
e.g.~\cite{schmuedgen.wagner:2003a} for their representation theory.

For later use we note the following simple and well-known fact
expressing the functoriality of the crossed product construction:
\begin{lemma}
    \label{lemma:CrossedFunctorial}
    If $\Phi: \mathcal{A} \longrightarrow \mathcal{B}$ is a
    $H$-equivariant $^*$-homomorphism then
    \begin{equation}
        \label{eq:PhiotimesIdCrossed}
        \Phi \otimes \id: \mathcal{A} \rtimes H 
        \longrightarrow \mathcal{B} \rtimes H
    \end{equation}
    is a $^*$-homomorphism. In particular, this induces a groupoid
    morphism
    \begin{equation}
        \label{eq:IsoStarHtoIsoStarCrossed}
        \cdot \rtimes\!H: \starIsoH 
        \longrightarrow
        \starIso,
    \end{equation}
    such that the identities $\mathcal{A}$ in $\starIsoH$ are mapped
    to their crossed products $\mathcal{A} \rtimes H$ with $H$ and
    arrows $\Phi$ are mapped to $\Phi \otimes \id$.
\end{lemma}

On $\mathcal{A} \rtimes H$ one has a canonical $^*$-action of $H$
defined by
\begin{equation}
    \label{eq:HactsAtimesH}
    g \acts (a \otimes h)
    = (g_\sweedler{1} \acts a) 
    \otimes g_\sweedler{2} h S(g_\sweedler{3})
\end{equation}
and there is a canonical $^*$-homomorphism
\begin{equation}
    \label{eq:AintoAH}
    \iota: \mathcal{A} 
    \ni a \mapsto \iota(a) = 
    a \otimes \Unit_H \in \mathcal{A} \rtimes H,
\end{equation}
which, up to possible torsion-effects due to $\tensor[\ring{C}]$, is
injective. Furthermore, $\iota$ is $H$-equivariant, i.e.~$g \acts
\iota(a) = \iota(g \acts a)$.

If $\mathcal{A}$ is unital then $\mathcal{A} \rtimes H$ is
unital with unit $\Unit_{\mathcal{A}} \otimes \Unit_H$ and we
have a canonical $^*$-homomorphism
\begin{equation}
    \label{eq:HintoAH}
    \jmath: H \ni g \mapsto \Unit_{\mathcal{A}} \otimes g \in
    \mathcal{A} \rtimes H,
\end{equation}
such that under this inclusion the action \eqref{eq:HactsAtimesH}
becomes `inner' in the sense that
\begin{equation}
    \label{eq:HactsAHInner}
    g \acts (a \otimes h)
    =
    \jmath(g_\sweedler{1}) (a \otimes h) \jmath(S(g_\sweedler{2})),
\end{equation}
see the adjoint action \eqref{eq:AdjointAction} of $H$ on itself.
Finally, in the unital case the crossed product is \emph{universal}
with respect to these properties, i.e.~if $\mathcal{B}$ is another
unital $^*$-algebra with two unital $^*$-homomorphisms
$\iota_{\mathcal{B}}: \mathcal{A} \longrightarrow \mathcal{B}$ and
$\jmath_{\mathcal{B}}: H \longrightarrow \mathcal{B}$ such that
$\iota_{\mathcal{B}}(g \acts a) = \jmath_{\mathcal{B}}(g_\sweedler{1})
\iota_{\mathcal{B}} (a) \jmath_{\mathcal{B}} (S (g_\sweedler{2}))$
then there exists a unique unital $^*$-homomorphism $\phi: \mathcal{A}
\rtimes H \longrightarrow \mathcal{B}$ such that $\iota_{\mathcal{B}}
= \phi \circ \iota$ and $\jmath_{\mathcal{B}} = \phi \circ \jmath$. In
fact, $\phi(a \otimes g) =
\iota_{\mathcal{B}}(a)\jmath_{\mathcal{B}}(g)$.

This observation immediately implies the following crucial property of
$\mathcal{A} \rtimes H$ which is one of the motivations to study
crossed products. This statement is well-known in various contexts.
\begin{lemma}
    \label{lemma:ReEpARepAHEquiv}
    The categories $\smod[H](\mathcal{A})$ and $\smod(\mathcal{A}
    \rtimes H)$ are equivalent, where the equivalence on objects is
    given by
    \begin{equation}
        \label{eq:RepHAtoRepAH}
        \smod[H](\mathcal{A}) \ni (\mathcal{H}, \pi)
        \mapsto
        (\hat{\mathcal{H}}, \hat{\pi}) \in
        \smod(\mathcal{A} \rtimes H),
    \end{equation}
    where $\hat{\mathcal{H}} = \mathcal{H}$ as pre-Hilbert spaces and
    $\hat{\pi}(a \otimes g)\phi = \pi(a) g \acts \phi$, and on
    morphisms $T:(\mathcal{H}_1, \pi_1) \longrightarrow
    (\mathcal{H}_2, \pi_2)$ it is the identity. The same statement
    holds for $\sMod$, $\rep$ and $\Rep$ instead of $\smod$, too.
\end{lemma}

The following proposition should be well-known and allows to construct
positive functionals for $\mathcal{A} \rtimes H$ and hence
$^*$-representations via the GNS construction.
\begin{proposition}
    \label{proposition:PosFunForCrossedProd}
    Let $\omega: \mathcal{A} \longrightarrow \ring{C}$ be a
    $H$-invariant positive linear functional and let $\mu: H
    \longrightarrow \ring{C}$ be a positive linear functional. Then
    $\omega \otimes \mu: \mathcal{A} \rtimes H \longrightarrow
    \ring{C}$ is again positive.
\end{proposition}
\begin{proof}
    Let $\sum_i a_i \otimes g_i$ be given. Then a
    straightforward computation using the invariance of $\omega$ gives
    \[
    (\omega \otimes \mu)\left(
        \left(\sum\nolimits_i a_i \otimes g_i\right)^*
        \left(\sum\nolimits_j a_j \otimes g_j\right)
    \right)
    =
    \sum\nolimits_{i, j} \omega(a_i^* a_j) \mu(g_i^* g_j)
    \ge 0,
    \]
    since both $\omega$ and $\mu$ are positive linear functionals and
    hence completely positive, see
    e.g.~\cite[Lem.~4.3]{bursztyn.waldmann:2001a}.
\end{proof}

In particular, $\omega \otimes \epsilon$ is always a positive linear
functional on $\mathcal{A} \rtimes H$ whence we can embed the
$H$-invariant positive functionals of $\mathcal{A}$ into the positive
linear functionals of $\mathcal{A} \rtimes H$. More generally, if
$\chi: H \longrightarrow \ring{C}$ is a \textdef{unitary character},
i.e.~a unital $^*$-homomorphism, then $\omega \otimes \chi$ is
positive, the counit is an example. We have the following invariance
with respect to the $^*$-action \eqref{eq:HactsAHInner}
\begin{equation}
    \label{eq:omegaotimeschiInv}
    (\omega \otimes \chi) (g \acts (a \otimes h))
    = \epsilon(g) (\omega \otimes \chi)(a \otimes h).
\end{equation}
\begin{remark}
    \label{remark:GNSForomegaotimesepsilon}
    If $\omega: \mathcal{A} \longrightarrow \ring{C}$ is $H$-invariant
    then the $H$-covariant GNS representation $(\mathcal{H}_\omega,
    \pi_\omega)$ of $\mathcal{A}$
    corresponds to the GNS representation $(\mathcal{H}_{\omega
      \otimes \epsilon}, \pi_{\omega\otimes \epsilon})$ of
    $\mathcal{A} \rtimes H$ under the functor
    \eqref{eq:RepHAtoRepAH}. In fact, the map
    \begin{equation}
        \label{eq:IntertwinerGNSGNS}
        U: (\hat{\mathcal{H}}_\omega, \hat{\pi}_\omega) \ni 
        \psi_a \mapsto \psi_{a \otimes \Unit_H} \in 
        (\mathcal{H}_{\omega\otimes\epsilon},
        \pi_{\omega\otimes\epsilon})
    \end{equation}
    is a unitary intertwiner which can be verified easily. Here
    $\psi_a$ and $\psi_{a \otimes \Unit_H}$ denote the equivalence
    classes of $a$ and $a \otimes \Unit_H$, respectively.
\end{remark}

%
%

\subsection{Crossed products of $^*$-representations}
\label{subsec:CrossedRep}

We shall now extend the crossed product construction to modules and
$^*$-representations.
\begin{lemma}
    \label{lemma:CrossedBimodule}
    Let $\BEA \in \smod[\mathcal{A}, H](\mathcal{B})$. Then on
    $\mathcal{E} \otimes H$ we have a $(\mathcal{B} \rtimes H,
    \mathcal{A} \rtimes H)$-bimodule structure defined by
    \begin{equation}
        \label{eq:CrossedBimoduleLeft}
        (b \otimes g) \cdot (x \otimes h)
        =
        (b \cdot (g_\sweedler{1} \acts x)) 
        \otimes g_\sweedler{2}h
    \end{equation}
    and
    \begin{equation}
        \label{eq:CrossedBimoduleRight}
        (x \otimes g) \cdot (a \otimes h)
        = (x \cdot (g_\sweedler{1} \acts a)) 
        \otimes g_\sweedler{2} h.
    \end{equation}
    Moreover,
    \begin{equation}
        \label{eq:CrossedInnerProductRight}
        \SPEoHAH{x \otimes g, y \otimes h}
        =
        \left(g_\sweedler{1}^* \acts \SPEA{x, y}\right)
        \otimes g_\sweedler{2}^* h
    \end{equation}
    defines a $(\mathcal{A} \rtimes H)$-valued inner product on
    $\mathcal{E} \otimes H$ such that
    \begin{equation}
        \label{eq:CrossedSPCompB}
        \SPEoHAH{(b \otimes g) \cdot (x \otimes h), y \otimes k}
        =
        \SPEoHAH{x \otimes h, (b \otimes g)^* \cdot (y \otimes k)}.
    \end{equation}
\end{lemma}
\begin{proof}
    It is a well-known straightforward computation to show that the
    definitions \eqref{eq:CrossedBimoduleLeft} and
    \eqref{eq:CrossedBimoduleRight} indeed give the described bimodule
    structure. Thus we have to prove
    that~\eqref{eq:CrossedInnerProductRight} is a $(\mathcal{A}
    \rtimes H)$-valued inner product. Clearly, it extends
    $\ring{C}$-sesquilinearily to $\mathcal{E} \otimes H$. We compute
    \[
    \begin{split}
        \left(\SPEoHAH{x \otimes g, y \otimes h}\right)^*
        &=
        \left(
            g_\sweedler{1}^* \acts \SPEA{x, y} \otimes
            g_\sweedler{2}^* h
        \right)^* \\
        &=
        \left(
            h_\sweedler{1}^* g_\sweedler{2} S(g_\sweedler{1}^*)^*
        \right) \acts
        \left(\SPEA{x, y}\right)^* 
        \otimes h_\sweedler{2}^* g_\sweedler{3} \\
        &=
        h_\sweedler{1}^* \acts \SPEA{y, x} 
        \otimes h_\sweedler{2}^* g \\
        &=
        \SPEoHAH{y \otimes h, x \otimes g}.
    \end{split}
    \]
    Moreover,
    \[
    \begin{split}
        \SPEoHAH{x \otimes g, (y \otimes h) \cdot (a \otimes k)}
        &=
        g_\sweedler{1}^* \acts 
        \SPEA{x, y \cdot(h_\sweedler{1} \acts a)}
        \otimes
        g_\sweedler{2}^* h_\sweedler{2} k \\
        &=
        \left(g_\sweedler{1}^* \acts \SPEA{x, y}\right)
        \left((g_\sweedler{2}^* h_\sweedler{1})\acts a \right)
        \otimes
        g_\sweedler{3}^* h_\sweedler{2} k \\
        &=
        \left( g_\sweedler{1}^* \acts \SPEA{x, y} \otimes
            g_\sweedler{2}^* h
        \right)
        (a \otimes k) \\
        &=
        \SPEoHAH{x \otimes g, y \otimes h}
        (a \otimes k),
    \end{split}
    \]
    whence $\SPEoHAH{\cdot,\cdot}$ is indeed a $(\mathcal{A} \rtimes
    H)$-valued inner product. Finally, we compute
    \[
    \begin{split}
        \SPEoHAH{(b \otimes g) \cdot (x \otimes h), y \otimes k}
        &=
        \SPEoHAH{(b \cdot g_\sweedler{1} \acts x) \otimes
          g_\sweedler{2}h, y \otimes k} \\
        &=
        \left(
            (h_\sweedler{1}^* g_\sweedler{2}^*) 
            \acts \SPEA{b \cdot g_\sweedler{1} \acts x, y}
        \right)
        \otimes h_\sweedler{2}^* g_\sweedler{3}^* k \\
        &=
        \left(
            (h_\sweedler{1}^* g_\sweedler{2}^*) 
            \acts \SPEA{g_\sweedler{1} \acts x, b^* \cdot y}
        \right)
        \otimes h_\sweedler{2}^* g_\sweedler{3}^* k \\
        &=
        \left(
            h_\sweedler{1}^* \acts
            \SPEA{
              (S(g_\sweedler{2}^*)^*g_\sweedler{1})
              \acts x,
              g_\sweedler{3}^* \acts (b^* \cdot y)}
        \right)
        \otimes 
        h_\sweedler{2}^* g_\sweedler{4}^* k \\
        &=
        \left(
            h_\sweedler{1}^* \acts
            \SPEA{x, g_\sweedler{1}^* \acts (b^* \cdot y)}
        \right)
        \otimes
        h_\sweedler{2}^* g_\sweedler{2}^* k \\
        &=
        \SPEoHAH{
          x \otimes h,
          g_\sweedler{1}^* \acts (b^* \cdot y) \otimes
          g_\sweedler{2}^* k}
        \\
        &=
        \SPEoHAH{
          x \otimes h,
          (b \otimes g)^* \cdot (y \otimes k)},
    \end{split}
    \]
    whence \eqref{eq:CrossedSPCompB} follows.
\end{proof}

It may happen that the inner product $\SPEoHAH{\cdot,\cdot}$ on
$\mathcal{E} \otimes H$ is degenerate. Thus we can pass in the usual
way to the quotient by the degeneracy space which is compatible with
the $(\mathcal{B} \rtimes H, \mathcal{A} \rtimes H)$-bimodule
structure as usual. We end up with an object in $\smod[\mathcal{A}
\rtimes H](\mathcal{B} \rtimes H)$ which we shall denote by
\begin{equation}
    \label{eq:EHDef}
    \mathcal{E} \rtimes H = 
    (\mathcal{E} \otimes H) \big/ (\mathcal{E} \otimes H)^\bot,
\end{equation}
always understood to be endowed with the $(\mathcal{B} \rtimes H,
\mathcal{A} \rtimes H)$-bimodule structure and the induced
$(\mathcal{A} \rtimes H)$-valued inner product which we shall denote
by $\SPEHAH{\cdot,\cdot}$. The next lemma shows that complete
positivity as well as strong non-degeneracy is always preserved:
\begin{lemma}
    \label{lemma:CrossedPositive}
    Let $\BEA \in \rep[\mathcal{A}, H](\mathcal{B})$. Then the inner
    product $\SPEoHAH{\cdot,\cdot}$ is completely positive, whence
    $\BEAH \in \rep[\mathcal{A} \rtimes H](\mathcal{B} \rtimes H)$ is
    a $^*$-representation on a pre-Hilbert module.  Moreover, if $\BEA
    \in \sMod[\mathcal{A}, H](\mathcal{B})$ then $\BEAH \in
    \sMod[\mathcal{A} \rtimes H](\mathcal{B} \rtimes H)$ and hence
    $\BEA \in \Rep[\mathcal{A}, H](\mathcal{B})$ implies $\BEAH \in
    \Rep[\mathcal{A}\rtimes H](\mathcal{B} \rtimes H)$.
\end{lemma}
\begin{proof}
    Let $\Phi^{(1)}, \ldots \Phi^{(n)} \in \mathcal{E} \otimes H$ be
    given and let $\Phi^{(\alpha)} = \sum_{i=1}^N x^{(\alpha)}_i
    \otimes g^{(\alpha)}_i$ with some $x^{(\alpha)}_i \in \mathcal{E}$
    and $g^{(\alpha)}_i \in H$, where without restriction $N$ is the
    same for all $\alpha = 1, \ldots, n$. Then
    \[
    \SPEoHAH{\Phi^{(\alpha)}, \Phi^{(\beta)}}
    =
    \sum\nolimits_{i,j=1}^N 
    \left(
        \left(g^{(\alpha)}_i\right)^*_\sweedler{1} \acts 
        \SPEA{x^{(\alpha)}_i, x^{(\beta)}_j}
    \right)
    \otimes
    \left(g^{(\alpha)}_i\right)^*_\sweedler{2} g^{(\beta)}_j.
    \]
    Now the map $f: M_{nN}(\mathcal{A}) \longrightarrow
    M_{nN}(\mathcal{A} \rtimes H)$ defined by
    \[
    f: A = \left(A^{\alpha\beta}_{ij}\right)
    \mapsto
    \left(
        \left(g^{(\alpha)}_i\right)^*_\sweedler{1} \acts
        A^{\alpha\beta}_{ij}
        \otimes
        \left(g^{(\alpha)}_i\right)^*_\sweedler{2} g^{(\beta)}_j
    \right)
    \]
    is a positive map. Indeed, we have
    \[
    \begin{split}
        f(A^*A) 
        &=
        \sum\nolimits_{\gamma, k}
        \left(
            \left(g^{(\alpha)}_i\right)^*_\sweedler{1} \acts
            \left(
                (A^{\gamma\alpha}_{ki})^* A^{\gamma\beta}_{kj}
            \right)
            \otimes 
            \left(g^{(\alpha)}_i\right)^*_\sweedler{2} g^{(\beta)}_j
        \right)\\
        &=
        \sum\nolimits_{\gamma, k}
        \left(
            \left(
                \left(g^{(\alpha)}_i\right)^*_\sweedler{1} \acts
                (A^{\gamma\alpha}_{ki})^*
                \otimes
                \left(g^{(\alpha)}_i\right)^*_\sweedler{2}
            \right)
            \left(
                A^{\gamma\beta}_{kj} \otimes g^{(\beta)}_j
            \right)
        \right)\\
        &=
        \sum\nolimits_{\gamma, k}
        \left(A^{\gamma\alpha}_{ki} \otimes g^{(\alpha)}_i\right)^*
        \left(A^{\gamma\beta}_{kj} \otimes g^{(\beta)}_j\right) \\
        &=
        (A \otimes g)^* (A \otimes g),
    \end{split}
    \]
    where $A \otimes g \in M_{nN}(\mathcal{A} \rtimes H)$ is given by
    its matrix coefficients $(A \otimes g)^{\alpha\beta}_{ij} =
    A^{\alpha\beta}_{ij} \otimes g^{(\beta)}_j$. Thus $f(A^*A) \in
    M_{nN}(\mathcal{A} \rtimes H)^{++}$ whence $f$ is positive.  Since
    the matrix $\left(\SPEA{x^{(\alpha)}_i, x^{(\beta)}_j}\right)$ is
    a positive matrix in $M_{nN}(\mathcal{A})$, by complete positivity
    of $\SPEA{\cdot,\cdot}$ we conclude that the matrix
    $f\left(\left(\SPEA{x^{(\alpha)}_i, x^{(\beta)}_j}\right)\right)$
    is positive as well. Then the summation over $i,j$ is the positive
    map $\tau$ from \cite[Ex.~2.1]{bursztyn.waldmann:2003a:pre} whence
    the result is positive again. This is precisely the matrix
    $\left(\SPEoHAH{\Phi^{(\alpha)}, \Phi^{(\beta)}}\right)$. Thus the
    complete positivity of the inner product $\SPEoHAH{\cdot,\cdot}$
    is shown and the complete positivity of $\SPEHAH{\cdot,\cdot}$
    follows. The statement on the strong non-degeneracy is trivial.
\end{proof}

In the unital case one can simplify the above argument by observing
that
\begin{equation}
    \label{eq:CompletePositiveCrossed}
    \SPEoHAH{x \otimes g, y \otimes h}
    =
    (\Unit_{\mathcal{A}} \otimes g)^* 
    \left(\SPEA{x,y} \otimes \Unit_H\right)
    (\Unit_{\mathcal{A}} \otimes h).
\end{equation}
From this the complete positivity of $\SPEoHAH{\cdot,\cdot}$ can be
deduced more easily.
\begin{remark}
    \label{remark:LeftLinearCrossed}
    For a left $\mathcal{B}$-linear $H$-covariant $\mathcal{B}$-valued
    inner product the corresponding definition of the left
    $(\mathcal{B} \rtimes H)$-linear $(\mathcal{B} \rtimes H)$-valued
    inner product on $\mathcal{E} \otimes H$ is
    \begin{equation}
        \label{eq:LeftCrossedInnerProduct}
        \begin{split}
            \BHSPEoH{x \otimes g, y \otimes h}
            &=
            \left(g_\sweedler{2} \acts 
                \BSPE{S^{-1}(g_\sweedler{1}) \acts x,
                  S^{-1}(h_\sweedler{1}) \acts y}
            \right)
            \otimes
            g_\sweedler{3} h_\sweedler{2}^*\\
            &=
            \BSPE{x, S(g_\sweedler{1})^* S^{-1}(h_\sweedler{1}) \acts y}
            \otimes g_\sweedler{2} h_\sweedler{2}^*.
        \end{split}
    \end{equation}
    The motivation for this formula comes from the isomorphism $I_2$
    in Proposition~\ref{proposition:CrossedTensorCompatible} below
    which identifies the complex conjugated bimodule $\cc{\mathcal{E}
      \rtimes H}$ canonically with $\cc{\mathcal{E}} \rtimes H$.  One
    can prove by an analogous computation that $\BHSPEoH{\cdot,\cdot}$
    is indeed $(\mathcal{B} \rtimes H)$-left linear and enjoys the
    correct symmetry properties. Moreover, it is compatible with the
    right $(\mathcal{A} \rtimes H)$-module structure, whence it gives
    a non-degenerate inner product $\BHSPEH{\cdot,\cdot}$ on the
    corresponding quotient $\mathcal{E} \rtimes H$. Finally,
    $\BHSPEoH{\cdot,\cdot}$ is completely positive if
    $\BSPE{\cdot,\cdot}$ is completely positive, as one can check
    directly in the same spirit as for $\SPEoHAH{\cdot,\cdot}$.
    Alternatively, we shall see an argument in
    Remark~\ref{remark:PositivityLeftCrossedIP}.
\end{remark}

The next lemma ensures the functoriality of the construction of
$\mathcal{E} \rtimes H$:
\begin{lemma}
    \label{lemma:CrossedBimoduleFunktor}
    Let $T: \BEA \longrightarrow \BEpA$ be an intertwiner between
    $\BEA, \BEpA \in \smod[\mathcal{A}, H](\mathcal{B})$. Then the map
    $T \otimes \id_H: \mathcal{E} \otimes H \longrightarrow
    \mathcal{E}' \otimes H$ descends to an intertwiner between
    $\mathcal{E} \rtimes H$ and $\mathcal{E}' \rtimes H \in
    \smod[\mathcal{A} \rtimes H](\mathcal{B} \rtimes H)$. The adjoint
    of $T \otimes \id$ is given by $T^* \otimes \id$.
\end{lemma}
\begin{proof}
    This is an easy verification using the $H$-equivariance of $T$ as
    well as the existence of $T^*$. In fact, everything is already
    true on the level of $\mathcal{E} \otimes H$ and $\mathcal{E}'
    \otimes H$.
\end{proof}

Collecting the results of the preceding lemmas we finally arrive at
the following statement:
\begin{proposition}
    \label{proposition:CrossedRepresentation}
    The assignment $\mathcal{E} \mapsto \mathcal{E} \rtimes H$ on
    objects and $T \mapsto T \otimes \id$ on morphisms gives a functor
    \begin{equation}
        \label{eq:CrossedmoduleFunctor}
        \cdot \rtimes\!H:
        \smod[\mathcal{A}, H](\mathcal{B}) \longrightarrow
        \smod[\mathcal{A} \rtimes H](\mathcal{B} \rtimes H)
    \end{equation}
    which restricts to functors
    \begin{equation}
        \label{eq:CrossedModuleFunctor}
        \cdot \rtimes\!H:
        \sMod[\mathcal{A}, H](\mathcal{B}) \longrightarrow
        \sMod[\mathcal{A} \rtimes H](\mathcal{B} \rtimes H)
    \end{equation}
    \begin{equation}
        \label{eq:CrossedrepFunctor}
        \cdot \rtimes\!H:
        \rep[\mathcal{A}, H](\mathcal{B}) \longrightarrow
        \rep[\mathcal{A} \rtimes H](\mathcal{B} \rtimes H)
    \end{equation}
    \begin{equation}
        \label{eq:CrossedRepFunctor}
        \cdot \rtimes\!H:
        \Rep[\mathcal{A}, H](\mathcal{B}) \longrightarrow
        \Rep[\mathcal{A} \rtimes H](\mathcal{B} \rtimes H).
    \end{equation}
\end{proposition}

In a next step we want to discuss the compatibility of the
\textdef{crossed product functors} \eqref{eq:CrossedmoduleFunctor},
\eqref{eq:CrossedModuleFunctor}, \eqref{eq:CrossedrepFunctor} and
\eqref{eq:CrossedRepFunctor}, respectively, with the tensor product
functors from~\eqref{eq:CovariantInternalTensorProduct} and
\eqref{eq:repCovInternalTensor}, respectively. Again, we only have to
investigate the case of $\smod$, the other cases will follow
analogously. 
\begin{proposition}
    \label{proposition:CrossedTensorCompatible}
    Let $\CFB \in \smod[\mathcal{B}, H](\mathcal{C})$ and $\BEA \in
    \smod[\mathcal{A}, H](\mathcal{B})$. Then we have:
    \begin{enumerate}
    \item The map
        \begin{equation}
            \label{eq:Ieins}
            I_1:
            (\mathcal{F} \rtimes H) 
            \tensM[\mathcal{B} \rtimes H]
            (\mathcal{E} \rtimes H)
            \ni 
            (x \otimes g) 
            \tensM[\mathcal{B} \rtimes H] 
            (y \otimes h)
            \mapsto
            (x \tensM[\mathcal{B}] (g_\sweedler{1} \acts y))
            \otimes g_\sweedler{2} h
            \in 
            (\mathcal{F} \tensM[\mathcal{B}] \mathcal{E}) \otimes H
        \end{equation}
        is a canonical isomorphism of $^*$-representations of
        $\mathcal{C} \rtimes H$ on $(\mathcal{A} \rtimes H)$-inner
        product modules.
    \item The map
        \begin{equation}
            \label{eq:Izwei}
            I_2:
            \cc{\mathcal{E} \rtimes H}
            \ni \cc{x \otimes g}
            \mapsto
            g_\sweedler{1}^* \ccacts \cc{x} \otimes g_\sweedler{2}^*
            =
            \cc{S^{-1}(g_\sweedler{1}) \acts x} 
            \otimes g_\sweedler{2}^*
            \in \cc{\mathcal{E}} \rtimes H
        \end{equation}
        is a canonical isomorphism of (right) $\mathcal{B} \rtimes
        H$-representations on (left) $(\mathcal{A} \rtimes H)$-inner
        product modules with inverse given explicitly by
        \begin{equation}
            \label{eq:IzweiInv}
            I_2^{-1} (\cc{x} \otimes g) = 
            \cc{g_\sweedler{1}^* \acts x \otimes g_\sweedler{2}^*}.
        \end{equation}
    \end{enumerate}
\end{proposition}
\begin{proof}
    For the first part one checks easily that $I_1$ is well-defined
    over $\otimes_{\mathcal{B}\rtimes H}$. Moreover, it is a
    straightforward computation that $I_1$ is a bimodule map as
    specified. For the isometry we compute
    \[
    \begin{split}
        &\SPFEHAH{
          (x \tensor[\mathcal{B}] (g_\sweedler{1} \acts y)) 
          \otimes g_\sweedler{2}h,
          (x' \tensor[\mathcal{B}] (g'_\sweedler{1} \acts y'))
          \otimes g'_\sweedler{2}h'
        }\\
        &=
        \left(
            (h_\sweedler{1}^* g_\sweedler{2}^*) 
            \acts 
            \SPFEA{
              x \tensor[\mathcal{B}] (g_\sweedler{1} \acts y),
              x' \tensor[\mathcal{B}] (g'_\sweedler{1} \acts y)'
            }
        \right)
        \otimes
        h_\sweedler{2}^* g_\sweedler{3}^* g'_\sweedler{2} h' \\
        &=
        \left(
            (h_\sweedler{1}^* g_\sweedler{2}^*) 
            \acts 
            \SPEA{
              g_\sweedler{1} \acts y,
              \SPFB{x, x'} \cdot (g'_\sweedler{1} \acts y')
            }
        \right)
        \otimes
        h_\sweedler{2}^* g_\sweedler{3}^* g'_\sweedler{2} h' \\
        &=
        \left(
            h_\sweedler{1}^* 
            \acts 
            \SPEA{
              (S(g_\sweedler{2}^*)^*g_\sweedler{1}) \acts y,
              g_\sweedler{3}^* \acts \left(
                  \SPFB{x, x'} \cdot (g'_\sweedler{1} \acts y')
              \right)
            }
        \right)
        \otimes
        h_\sweedler{2}^* g_\sweedler{4}^* g'_\sweedler{2} h' \\
        &=
        \left(
            h_\sweedler{1}^* 
            \acts 
            \SPEA{
              y,
              (g_\sweedler{1}^* \acts \SPFB{x, x'}) \cdot 
              (g_\sweedler{2}^*g'_\sweedler{1}) \acts y'
            }
        \right)
        \otimes
        h_\sweedler{2}^* g_\sweedler{3}^* g'_\sweedler{2} h' \\
        &=
        \SPEHAH{
          y \otimes h,
          \left(
              (g_\sweedler{1}^* \acts \SPFB{x, x'}) \cdot 
              (g_\sweedler{2}^*g')_\sweedler{1} \acts y'
          \right)
          \otimes
          (g_\sweedler{2}^* g')_\sweedler{2} h'
        }\\
        &=
        \SPEHAH{
          y \otimes h,
          \left(
              (g_\sweedler{1}^* \acts \SPFB{x, x'})
              \otimes g_\sweedler{2}^*g'
          \right)
          \cdot
          (y' \otimes h')
        }\\
        &=
        \SPEHAH{
          y \otimes h,
          \SPFHBH{
            x \otimes g, 
            x' \otimes g
          }
          \cdot
          (y' \otimes h')
        } \\
        &=
        \SPFHEHAH{
          (x \otimes g) \tensor[\mathcal{B}] (y \otimes h),
          (x' \otimes g') \tensor[\mathcal{B}] (y' \otimes h')
        },
    \end{split}
    \]
    whence $I_1$ is isometric already on the level of
    $\tensor[\mathcal{B}]$ instead of $\tensM[\mathcal{B}]$. Finally,
    surjectivity is clear since $(x \otimes \Unit_H)
    \tensor[\mathcal{B}] (y \otimes g)$ is mapped to $(x
    \tensor[\mathcal{B}] y) \otimes g$. The injectivity follows as on
    the quotients both inner products are, by definition,
    non-degenerate whence an isometric map is injective. This shows
    that $I_1$ is an isomorphism indeed.  Moreover, it is canonical in
    the following sense: Let $S: \mathcal{F} \longrightarrow
    \mathcal{F}'$ and $T: \mathcal{E} \longrightarrow \mathcal{E}'$ be
    morphisms in $\smod[\mathcal{B}, H](\mathcal{C})$ and
    $\smod[\mathcal{A}, H](\mathcal{B})$, respectively. Then $S \tensM
    T$ is a morphism in $\smod[\mathcal{A}, H](\mathcal{C})$ and $S
    \otimes \id$ and $T \otimes \id$ are the corresponding morphisms
    in $\smod[\mathcal{B} \rtimes H](\mathcal{C} \rtimes H)$ and
    $\smod[\mathcal{A} \rtimes H](\mathcal{B} \rtimes H)$,
    respectively, according to
    Lemma~\ref{lemma:CrossedBimoduleFunktor}. Then $I_1$ is compatible
    with morphisms as it is easy to check that
    \[
    I_1 \circ \left((S \otimes \id)\tensM (T \otimes \id)\right)
    = \left((S \tensM T) \otimes \id\right) \circ I_1.
    \]
    This proves the first part. For the second we first observe that
    $I_2$ certainly has the correct $\ring{C}$-linearity properties. A
    lengthy but straightforward computation shows by successively
    unwinding the definitions that $I_2$ is a bimodule map. Thus we
    compute using \eqref{eq:ccacts} and
    \eqref{eq:LeftCrossedInnerProduct}
    \[
    \begin{split}
        &\AHSPccEH{I_2(\cc{x \otimes g}), I_2(\cc{y \otimes h})}
        \\
        &=
        \AHSPccEH{
          (g_\sweedler{1}^* \ccacts \cc{x}) \otimes g_\sweedler{2}^*,
          (h_\sweedler{1}^* \ccacts \cc{y}) \otimes h_\sweedler{2}^*
        }\\
        &=
        \AHSPccEH{
          \cc{S(g_\sweedler{1}^*)^* \acts x} \otimes g_\sweedler{2}^*,
          \cc{S(h_\sweedler{1}^*)^* \acts y} \otimes h_\sweedler{2}^*
        }\\
        &=
        \left(
            g_\sweedler{3}^* \acts
            \ASPccE{
              S^{-1}(g_\sweedler{2}^*) \ccacts 
              \left(\cc{S^{-1}(g_\sweedler{1}) \acts x}\right),
              S^{-1}(h_\sweedler{2}^*) \ccacts
              \left(\cc{S^{-1}(h_\sweedler{1}) \acts y}\right)
            }
        \right)
        \otimes
        g_\sweedler{4}^* h_\sweedler{3} \\
        &=
        \left(
            g_\sweedler{3}^* \acts
            \ASPccE{
              \cc{(g_\sweedler{2} S^{-1}(g_\sweedler{1})) \acts x},
              \cc{(h_\sweedler{2} S^{-1}(h_\sweedler{1})) \acts y}
            }
        \right)
        \otimes 
        g_\sweedler{4}^* h_\sweedler{3} \\
        &=
        \left(
            g_\sweedler{1}^* \acts \SPEA{x, y}
        \right)
        \otimes 
        g_\sweedler{2}^* h \\
        &=
        \AHSPCCEH{\cc{x \otimes g}, \cc{y \otimes h}},
    \end{split}
    \]
    whence $I_2$ is isometric. It is a simple computation that
    \eqref{eq:IzweiInv} provides an inverse for $I_2$. The
    compatibility with intertwiners is shown analogously as for $I_1$.
\end{proof}
\begin{remark}
    \label{remark:PositivityLeftCrossedIP}
    Using the second part of the proposition we also obtain an easy
    proof for the complete positivity of the inner product on
    $\mathcal{E} \rtimes H$ if we had a \emph{left}
    $\mathcal{B}$-linear $\mathcal{B}$-valued inner product on
    $\mathcal{E}$. In this case we can pass to $\cc{\mathcal{E}}$
    instead, making the inner product right $\mathcal{B}$-linear and
    use $\cc{\mathcal{E}}\rtimes H$, which has, by
    Lemma~\ref{lemma:CrossedPositive}, a completely positive right
    $(\mathcal{B} \rtimes H)$-linear $(\mathcal{B} \rtimes H)$-valued
    inner product. This is isometric to the corresponding inner
    product on $\cc{\mathcal{E} \rtimes H}$ by
    Proposition~\ref{proposition:CrossedTensorCompatible} and by
    Remark~\ref{remark:LeftRightAndCC} the complete positivity of the
    inner product on $\mathcal{E} \rtimes H$ follows.
\end{remark}

Rephrasing the statement of the proposition in terms of the functors
$\tensM$ and $\cdot \rtimes\!H$ we have the following result:
\begin{corollary}
    \label{corollary:TensRtimesHCommute}
    The diagram
    \begin{equation}
        \label{eq:TensCrossCommute}
        \bfig
        \Square/>`>`>`>/%
        [{\smod[\mathcal{B}, H](\mathcal{C})} \times%
        {\smod[\mathcal{A}, H](\mathcal{B})}%
        `{\smod[\mathcal{A},H](\mathcal{C})}%
        `{\smod[\mathcal{B} \rtimes H](\mathcal{C}\!\rtimes\!H)}%
        \times {\smod[\mathcal{A}\rtimes H](\mathcal{B}\!\rtimes\!H)}%
        `{\smod[\mathcal{A}\rtimes H](\mathcal{C}\!\rtimes\!H)}%
        ;\tensM%
        `(\cdot \rtimes H) \times (\cdot \rtimes H)%
        `\cdot\rtimes H%
        `\tensM]
        \efig
    \end{equation}
    commutes in the sense of functors, i.e.~up to the natural
    transformation $I_1$.
\end{corollary}
Analogous statements hold for the complex conjugation exchanging left
and right linear inner products.  Let us also remark that
Proposition~\ref{proposition:CrossedTensorCompatible} still holds if
we restrict ourselves to $^*$-representations in $\sMod$, $\rep$ or
$\Rep$, respectively.

%
%

\subsection{The Picard groupoid of crossed products}
\label{subsec:PicardForCrossedProducts}

After our discussion of crossed product constructions for general
$^*$-representations we turn now to the equivalence bimodules. The
first lemma ensures that the functor $\cdot\rtimes\!H$ applied to
equivalence bimodules gives again equivalence bimodules:
\begin{lemma}
    \label{lemma:CrossedEquivalenceBimodule}
    Let $\BEA$ be an $H$-covariant $^*$-Morita equivalence bimodule.
    Then $\BEAH$, endowed with the induced $(\mathcal{B} \rtimes
    H)$-left linear $(\mathcal{B} \rtimes H)$-valued inner product
    $\BHSPEH{\cdot,\cdot}$ and the induced right $(\mathcal{A} \rtimes
    H)$-linear $(\mathcal{A} \rtimes H)$-valued inner product
    $\SPEHAH{\cdot,\cdot}$, is a $^*$-Morita equivalence bimodule for
    $\mathcal{B} \rtimes H$ and $\mathcal{A} \rtimes H$. Moreover, if
    $\BEA$ is even a strong equivalence bimodule then $\BEAH$ is a
    strong equivalence bimodule as well.
\end{lemma}
\begin{proof}
    We already know that on $\mathcal{E} \otimes H$ we have two inner
    products $\BHSPEoH{\cdot,\cdot}$ and $\SPEoHAH{\cdot,\cdot}$ which
    have the correct linearity and compatibility with respect to the
    $(\mathcal{B}\rtimes H, \mathcal{A} \rtimes H)$-bimodule
    structure. To show the compatibility of the inner products we
    compute
    \[
    \begin{split}
        \BHSPEoH{x \otimes g, y \otimes h} \cdot (z \otimes k)
        &= 
        \left(
            \BSPE{x, 
              (S(g_\sweedler{1})^*S^{-1}(h_\sweedler{1})) \acts y}
        \right)
        \cdot (z \otimes k) \\
        &= 
        \BSPE{x, 
          (S(g_\sweedler{1})^*S^{-1}(h_\sweedler{1})) \acts y}
        \cdot 
        \left( (g_\sweedler{2} h_\sweedler{2}^*)\acts z\right)
        \otimes
        g_\sweedler{3} h_\sweedler{3}^* k \\
        &=
        x \cdot
        \SPEA{(S(g_\sweedler{1})^*S^{-1}(h_\sweedler{1})) \acts y,
          (g_\sweedler{2} h_\sweedler{2}^*)\acts z}
        \otimes 
        g_\sweedler{3} h_\sweedler{3}^* k \\
        &=
        (x \otimes g) \cdot
        \left(
            \SPEA{
              S^{-1}(h_\sweedler{1}) \acts y,
              h_\sweedler{2}^* \acts z}
            \otimes
            h_\sweedler{3}^* k
        \right) \\
        &=
        (x \otimes g) \cdot
        \left(
            h_\sweedler{1}^* \acts
            \SPEA{y,z} \otimes h_\sweedler{2}^*k
        \right)\\
        &=
        (x \otimes g)
        \cdot
        \SPEoHAH{y \otimes h, z \otimes k},
    \end{split}
    \]
    whence \eqref{eq:SPCompatible} follows. Thus it follows that their
    degeneracy spaces coincide whence we can non-ambiguously define
    $\mathcal{E} \rtimes H$ and obtain a $(\mathcal{B} \rtimes H,
    \mathcal{A} \rtimes H)$-bimodule with compatible non-degenerate
    inner products $\BHSPEH{\cdot,\cdot}$ and $\SPEHAH{\cdot,\cdot}$.
    Moreover, since $\mathcal{B} \cdot \mathcal{E} = \mathcal{E} =
    \mathcal{E} \cdot \mathcal{A}$ it is easy to check that
    $\mathcal{E} \rtimes H$ is also strongly non-degenerate for both
    module structures (in the unital case this is trivial). It remains
    to check whether the inner products are full. If $a = \sum_i
    \SPEA{x_i, y_i}$ by fullness of $\SPEA{\cdot,\cdot}$ then
    \[
    a \otimes g = \sum\nolimits_i \SPEHAH{x_i \otimes \Unit_H, y_i
      \otimes g}
    \]
    implies fullness of $\SPEHAH{\cdot,\cdot}$ and analogously for
    $\BHSPEH{\cdot,\cdot}$. Thus $\BEAH$ is indeed a $^*$-Morita
    equivalence bimodule. Since complete positivity of inner products
    is preserved under the crossed product construction by
    Lemma~\ref{lemma:CrossedPositive}, the remaining statement follows
    as well.
\end{proof}
\begin{remark}
    \label{remark:CstarVersionCrossed}
    This lemma has a remarkable and well-known counterpart in
    $C^*$-algebra theory: Here it is known that for a locally compact
    group $G$ acting (strongly) continuously on $C^*$-algebras,
    $G$-covariant strong Morita equivalence implies strong Morita
    equivalence of the corresponding $C^*$-algebraic crossed products,
    see \cite{curto.muhly.williams:1984a, combes:1984a}. The above
    lemma reproduces this statements e.g.~for the case of a discrete
    group by using the group algebra $H = \mathbb{C}[G]$. It will be
    left to a future project to investigate topological versions of
    the above lemma in order to recover the statements of
    \cite{curto.muhly.williams:1984a, combes:1984a} fully. We also
    refer to \cite{echterhoff.kaliszewski.quigg.raeburn:2002a:pre} for
    more general crossed product constructions in the $C^*$-algebraic
    framework related to Morita theory and for further references.
\end{remark}

Adapting Lemma~\ref{lemma:CrossedBimoduleFunktor} to equivalence
bimodules we immediately have to following statement:
\begin{lemma}
    \label{lemma:CrossedEquivalenceFunctiorial}
    Let $\BEA$, $\BEpA$ be isomorphic $H$-covariant $^*$-Morita
    equivalence bimodules such that $T: \BEA \longrightarrow \BEpA$ is
    an isomorphism. Then $T\otimes \id: \BEAH \longrightarrow \BEpAH$
    is an isomorphism of $^*$-Morita equivalence bimodules. The
    analogous statement holds for strong Morita equivalence bimodules.
\end{lemma}

To obtain a good Morita theory for crossed products we have to
guarantee that $\mathcal{A} \rtimes H$ is idempotent and
non-degenerate. While it is easy to see that $\mathcal{A} \rtimes H$
is idempotent if $\mathcal{A}$ is idempotent, there could be some
torsion-effects due to $\otimes_{\ring{C}}$ which spoil the
non-degeneracy of $\mathcal{A} \rtimes H$ even if $\mathcal{A}$ was
non-degenerate. Nevertheless it is safe to assume that $\mathcal{A}
\rtimes H$ is non-degenerate as e.g.~in the unital case $\mathcal{A}
\rtimes H$ is unital and thus non-degenerate. Also if $\ring{C}$ is a
field we will have no problems. Thus we shall ignore these subtleties
in the following and always assume that all occurring crossed products
are non-degenerate.
\begin{lemma}
    \label{lemma:AAACrossedHtoAcrossedH}
    The map
    \begin{equation}
        \label{eq:AAAHToACrossedH}
        I_3:
        \AAA \rtimes H \ni x \otimes g
        \mapsto x \otimes g \in \AHAHAH
    \end{equation}
    is an isomorphism of strong Morita equivalence bimodules.
\end{lemma}
\begin{proof}
    Thanks to the assumption that $\mathcal{A} \rtimes H$ is
    non-degenerate, the canonical inner products on $\mathcal{A}
    \rtimes H$ are non-degenerate whence there is no quotient
    procedure necessary. Then an easy check shows that the bimodule
    structures and the inner products on $\mathcal{A} \otimes H$
    coming from both interpretations simply coincide.
\end{proof}

Now we can formulate the main result of this section which relates the
$H$-covariant Picard groupoid to the Picard groupoid of the
corresponding crossed products:
\begin{theorem}
    \label{theorem:PicHtoPicCrossed}
    The crossed product with $H$ induces groupoid morphisms
    \begin{equation}
        \label{eq:starPicHtostarPicCrossed}
        \cdot \rtimes\!H:
        \starPicH \longrightarrow \starPic
    \end{equation}
    and
    \begin{equation}
        \label{eq:StrPicHToStrPicCrossed}
        \cdot \rtimes\!H:
        \StrPicH \longrightarrow \StrPic,
    \end{equation}
    where units $\mathcal{A} \in \starPicH$ (or $\StrPicH$,
    respectively) are mapped to the crossed product algebras
    $\mathcal{A} \rtimes H$ and arrows $[\BEA] \in
    \starPicH(\mathcal{B}, \mathcal{A})$ (or $\StrPicH(\mathcal{B},
    \mathcal{A})$, respectively) are mapped to the arrows $[\BEAH] \in
    \starPic(\mathcal{B} \rtimes H, \mathcal{A} \rtimes H)$ (or
    $\StrPic(\mathcal{B} \rtimes H, \mathcal{A} \rtimes H)$,
    respectively).
\end{theorem}
\begin{proof}
    From Lemma~\ref{lemma:CrossedEquivalenceFunctiorial} we see that
    $\cdot\rtimes\!H$ is well-defined on isomorphism classes.
    Moreover, Lemma~\ref{lemma:AAACrossedHtoAcrossedH} ensures that
    units are mapped to units indeed, in the stated way. Finally,
    Proposition~\ref{proposition:CrossedTensorCompatible} can easily
    be adapted to the case of two compatible inner products whence it
    shows that tensor products are mapped to tensor products: for
    equivalence bimodules one inner product determines the other by
    compatibility \eqref{eq:SPCompatible}. Again,
    Proposition~\ref{proposition:CrossedTensorCompatible} shows that
    complex conjugated bimodules are mapped to complex conjugated
    bimodules whence products and inverses in $\starPicH$ are mapped
    to products and inverses in $\starPic$. The same holds for
    $\StrPicH$ and $\StrPic$, whence
    \eqref{eq:starPicHtostarPicCrossed} and
    \eqref{eq:StrPicHToStrPicCrossed} are groupoid morphisms indeed.
\end{proof}

Let us now collect a few conclusions from this theorem. The first
transfers well-known results from $C^*$-algebra theory
\cite{curto.muhly.williams:1984a, combes:1984a}, where strongly
continuous actions of a locally compact group is considered, to our
algebraic framework.
\begin{corollary}
    \label{corollary:HMEImpliesMECrossed}
    If $\mathcal{A}$ and $\mathcal{B}$ are $H$-covariantly strongly
    Morita equivalent (or $^*$-Morita equivalent, respectively) then
    $\mathcal{A} \rtimes H$ and $\mathcal{B} \rtimes H$ are strongly
    Morita equivalent ($^*$-Morita equivalent, respectively).
\end{corollary}
\begin{corollary}
    \label{corollary:PicHPicGroupMorphisms}
    There are group homomorphisms
    \begin{equation}
        \label{eq:StarPicHstarPicGroupCrossed}
        \starPicH(\mathcal{A}) 
        \longrightarrow \starPic(\mathcal{A} \rtimes H)
    \end{equation}
    and
    \begin{equation}
        \label{eq:StrPicHStrPicGroupCrossed}
        \StrPicH(\mathcal{A}) 
        \longrightarrow \StrPic(\mathcal{A} \rtimes H).
    \end{equation}
\end{corollary}

The following easy proposition shows that the groupoid morphism in
Theorem~\ref{theorem:PicHtoPicCrossed} is compatible with the (much
easier) groupoid morphism from \eqref{eq:IsoStarHtoIsoStarCrossed} via
the canonical groupoid morphism $\ell$ from
Proposition~\ref{proposition:IsoToPicMorphism}.
\begin{proposition}
    \label{proposition:CrossedHCompatibleEll}
    The groupoid morphism $\cdot\rtimes\!H$ is compatible with the
    groupoid morphism $\ell$, i.e.~the diagram
    \begin{equation}
        \label{eq:OnceAgainCommutativeDiagram}
        \bfig
        \Square/>`>`>`>/%
        [\starIsoH`\starIso`\StrPicH`\StrPic%
        ;\cdot\rtimes H`\ell`\ell`\cdot\rtimes H]
        \efig         
    \end{equation}
    commutes.
\end{proposition}
\begin{proof}
    Let $\Phi \in \starIsoH(\mathcal{B}, \mathcal{A})$. Then
    $\ell(\Phi \otimes \id)$ is represented by the bimodule
    $\mathcal{B} \otimes H$ endowed with the usual left $(\mathcal{B}
    \rtimes H)$-module structure and the right $(\mathcal{A} \rtimes
    H)$-module structure
    \[
    (b \otimes g) \TMult{\Phi \otimes \id} (a \otimes h)
    =
    (b \otimes g) \cdot (\Phi(a) \otimes h)
    =
    (b \Phi(g_\sweedler{1} \acts a)) \otimes g_\sweedler{2}h
    =
    (b \dotPhi (g_\sweedler{1} \acts a)) \otimes g_\sweedler{2} h,
    \]
    which coincides with the right $(\mathcal{A} \rtimes H)$-module
    structure on $\BBPhiA \rtimes H$ thanks to the $H$-equivariance of
    $\Phi$. Analogously, for the $(\mathcal{A} \rtimes H)$-valued
    inner products we have
    \[
    \begin{split}
        \SPBHPhiIdAH{b \otimes g, b' \otimes g'}
        &= (\Phi \otimes \id)^{-1} 
        \left(
            \SPBHBH{b \otimes g, b' \otimes g'}
        \right) \\
        &=
        (\Phi^{-1} \otimes \id)
        \left(
            (g_\sweedler{1}^* \acts b^*)(g_\sweedler{2}^* \acts b')
            \otimes 
            g_\sweedler{3}^* g'
        \right) \\
        &=
        (g_\sweedler{1}^* \acts \Phi^{-1}(b^*b')) 
        \otimes g_\sweedler{2}^*g' \\
        &=
        (g_\sweedler{1}^* \acts \SPBPhiA{b, b'})
        \otimes g_\sweedler{2}^*g' \\
        &=
        \SPBPhiHAH{b \otimes g, b' \otimes g'},
    \end{split}
    \]
    whence they coincide, too. The left $(\mathcal{B} \rtimes
    H)$-module structure and the $(\mathcal{B} \rtimes H)$-valued
    inner products are unchanged whence the statement follows.
\end{proof}

%
%

\subsection{An example: $\StrPicH(\ring{C}) \longrightarrow
  \StrPic(H)$}
\label{subsec:ExamplePicHCToPicH}

We illustrate our general techniques by a simple but yet interesting
example where $\mathcal{A} = \ring{C}$, endowed with the trivial
action of $H$. First we need the following lemma:
\begin{lemma}
    \label{lemma:GLHCintoAutH}
    Let $\chi \in \group{GL}(H, \ring{C})$.
    \begin{compactenum}
    \item $\chi^{-1}(g) = \chi(S^{-1}(g)) = \chi(S(g))$.
    \item $\chi \in \group{U}(H, \ring{C})$ if and only if $\chi(g^*)
        = \cc{\chi(g)}$.
    \item $\Phi^\chi(g) = \chi(S(g_\sweedler{1})) g_\sweedler{2}$
        defines an automorphism $\Phi^\chi \in \Aut(H)$ and $\Phi^\chi
        \in \starAut(H)$ if $\chi \in \group{U}(H, \ring{C})$.
    \item The map
        \begin{equation}
            \label{eq:GLHCintoAut}
            \group{GL}(H, \ring{C}) 
            \ni \chi \mapsto \Phi^\chi \in 
            \Aut(H)
        \end{equation}
        is an injective group homomorphism.
    \item $\Phi^\chi$ is an inner automorphism if and only if $\chi =
        \twist{e}$.
    \end{compactenum}
\end{lemma}
\begin{proof}
    Clearly, $\chi^{-1}(g) = \chi(S^{-1}(g))$ by
    \eqref{eq:InverseTwist} since the action is trivial. Moreover, one
    easily checks that $\tilde{\chi}(g) = \chi(S(g))$ defines a
    inverse with respect to the convolution product whence by
    uniqueness $\tilde{\chi} = \chi^{-1}$. For the second part we
    define $\cc{\chi}(g) = \cc{\chi(\cc{g})}$. Then a straightforward
    computation using the first part and the unitarity condition for
    $\chi$ shows that $\cc{\chi}$ is a convolution inverse of
    $\chi^{-1}$ and thus equal to $\chi$. The converse direction if
    trivial. For the third and fourth part we have $\Phi^{\twist{e}} =
    \id$ and $\Phi^\chi(gh) = \Phi^\chi(g) \Phi^\chi(h)$ by a little
    computation. Hence $\Phi^\chi$ is a homomorphism and if $\chi \in
    \group{U}(H, \ring{C})$ one immediately has $\Phi^\chi(g^*) =
    \Phi^\chi(g)^*$. Next we prove $\Phi^\chi(\Phi^{\tilde{\chi}}(g))
    = \Phi^{\chi * \tilde{\chi}}(g)$ again by a simple computation.
    Then it follows that $\Phi^\chi$ is bijective and
    \eqref{eq:GLHCintoAut} is a group homomorphism.  For the
    injectivity assume $\Phi^\chi(g) = g$. Applying $\epsilon$ to this
    equation gives immediately $\chi(S(g)) = \epsilon(g)$ whence $\chi
    = \twist{e}$.
\end{proof}

This lemma is the Hopf-algebraic version of the well-known
construction of automorphisms of the group algebra $\ring{C}[G]$ out
of group characters of a group $G$. It shows that $\group{U}(H,
\ring{C})$ always gives a non-trivial contribution to $\StrPic(H)$: In
fact, from \cite[Eq.~(2.4)]{bursztyn.waldmann:2004a}, i.e.~the
non-covariant version of \eqref{eq:ellIsoToPic}, we have an injective
group homomorphism
\begin{equation}
    \label{eq:ellPhichiintoStrPic}
    \group{U}(H, \ring{C}) \ni \chi
    \mapsto
    \ell(\Phi^\chi) \in \StrPic(H),
\end{equation}
where injectivity follows from the last part of the lemma.  On the
other hand, we know that the crossed product algebra $\ring{C} \rtimes
H$ is just $H$ itself where the canonical identification is simply $z
\otimes g \mapsto zg$. Thus the general groupoid morphism
$\cdot\rtimes\!H$ from \eqref{eq:StrPicHToStrPicCrossed} gives a group
homomorphism
\begin{equation}
    \label{eq:StrPicHCtoStrPicH}
    \StrPicH(\ring{C}) \longrightarrow \StrPic(H),
\end{equation}
which we shall relate to \eqref{eq:ellPhichiintoStrPic}. From
Remark~\ref{remark:TrivialesZentrum} we know that $\group{U}(H,
\ring{C}) = \group{U}_0(H, \ring{C})$ and from
Proposition~\ref{proposition:UnullInStrPicH} we know that we can view
$\group{U}_0(H, \ring{C})$ as a subgroup of $\StrPicH(\ring{C})$.
Putting these group homomorphisms together we obtain the following
statement:
\begin{proposition}
    \label{proposition:StrPicHCinStrPicHUHCAutH}
    The diagram of group homomorphisms
    \begin{equation}
        \label{eq:StrPicHCUHCundsoweiter}
        \bfig
        \Square/^{ (}->`^{ (}->`>`>/%
        [\group{U}(H, \ring{C})%
        `\StrPicH(\ring{C})%
        `\starAut(H)%
        `\StrPic(H)%
        ;``\cdot\rtimes\!H`\ell]
        \efig
    \end{equation}
    commutes whence $\group{U}(H, \ring{C})$ can be viewed as a
    subgroup of $\StrPic(H)$.
\end{proposition}
\begin{proof}
    Let $\chi \in \group{U}(H, \ring{C}) = \group{U}_0(H, \ring{C})$.
    Then the image of $\chi$ in $\StrPicH(\ring{C})$ is given by the
    isomorphism class of the trivial bimodule $\ring{C}$ with
    canonical inner products and $H$-action $g \actsT[\chi] z =
    \chi(g_\sweedler{1}) g_\sweedler{2} \acts z = \chi(g)z$. We denote
    this bimodule by $\ring{C}^\chi$. Then $[\ring{C}^\chi]$ is mapped
    to $[\ring{C}^\chi \rtimes H]$ where $\ring{C}^\chi \rtimes H
    \cong H$ as $\ring{C}$-modules and the left $H$-module structure is
    given by
    \[
    g \cdot h = \chi(g_\sweedler{1}) g_\sweedler{2} h =
    \Phi^{\chi^{-1}}(g) h = g \TMult{\Phi^\chi} h,
    \]
    while the right $H$-module structure is the canonical one. The
    left-linear inner product is easily shown to be $\Phi^\chi (g
    h^*)$ and the right-linear inner product is the canonical one. Thus
    $\ring{C}^\chi \rtimes H$ is isomorphic to
    $\Bimod{\Phi^\chi}{\mathit{H}}{\mathit{H}}{}{\mathit{H}}$ whose
    class in $\StrPic(H)$ is just $\ell(\Phi^\chi)$. This proves the
    commutativity of \eqref{eq:StrPicHCUHCundsoweiter}. The
    injectivity of the inclusion of $\group{U}(H, \ring{C})$ into
    $\StrPic(H)$ was shown in \eqref{eq:ellPhichiintoStrPic}.
\end{proof}

If $\ring{C}$ is even an algebraically closed field (which is the case
if $\ring{R}$ is a real closed field, see
e.g.~\cite[Sect.~5.1]{jacobson:1985a}) then we can make
\eqref{eq:StrPicHCUHCundsoweiter} more precise:
\begin{corollary}
    \label{corollary:CAlgClosedField}
    Let $\ring{R}$ be real closed field whence $\ring{C}$ is
    algebraically closed.
    \begin{compactenum}
    \item $\StrPic(\ring{C}) = \{\id\}$.
    \item $\StrPicH(\ring{C}) = \group{U}(H, \ring{C}) =
        \group{U}_0(H, \ring{C})$.
    \item $\StrPicH(\ring{C}) \longrightarrow \StrPic(H)$ is injective
        and its image is determined by
        \eqref{eq:StrPicHCUHCundsoweiter}.
    \end{compactenum}
\end{corollary}
\begin{proof}
    The first part is clear since the only equivalence bimodule (up to
    isomorphism) is the one-dimensional vector space $\ring{C}$ with
    the (uniquely determined up to isometries) canonical positive
    definite inner product $\SP{z, w} = \cc{z}w$. Note that
    $\starPic(\ring{C}) = \mathbb{Z}_2$ in this case.  Then the second
    part follows from Proposition~\ref{proposition:UnullInStrPicH}
    whence the third part follows from
    Proposition~\ref{proposition:StrPicHCinStrPicHUHCAutH}.
\end{proof}
\begin{remark}
    \label{remark:CrossedNotSurjective}
    Though $\cdot\rtimes\!H: \StrPicH(\ring{C}) \longrightarrow
    \StrPic(H)$ is injective in this example it needs not to be
    surjective: Neither the map $\group{U}(H, \ring{C})
    \longrightarrow \starAut(H)$ nor the map $\ell$ need to be
    surjective. An example can be obtained for a \emph{commutative}
    Hopf algebra $H$. In this case the antipode $S$ is a
    $^*$-automorphism, $S \in \starAut(H)$, and if $\Phi^\chi = S$ for
    some $\chi$ then applying $\epsilon$ gives $\chi(g) = \epsilon(g)$
    whence $S = \id$ follows.  Beside this case, which is rarely of
    interest, we conclude that $\ell(S) \in \StrPic(H)$ gives a
    non-trivial element not in the image of $\StrPicH(\ring{C})$.  In
    general, the interesting elements of $\StrPic(H)$ are those which
    are \emph{not} in the image of $\ell$ anyway, see e.g.~the
    discussion in \cite[Sect.~2]{bursztyn.waldmann:2004a}. According
    to Proposition~\ref{proposition:StrPicHCinStrPicHUHCAutH} they can
    never be obtained from $\StrPicH(\ring{C})$.
\end{remark}
\begin{example}
    \label{example:CharactersOfLieAlgebra}
    To have a more concrete example we consider again $H =
    U_{\ring{C}}(\lie{g})$. Then $\twist{a} \in \group{U}(H,
    \ring{C})$ satisfies $\twist{a}(\xi X) = \twist{a}(\xi)
    \twist{a}(X)$ by use of the action condition, where $\xi \in
    \lie{g}$ and $X \in U_{\ring{C}}(\lie{g})$ since $\epsilon(\xi) =
    0$. Since $\lie{g}$ together with $\Unit$ generates
    $U_{\ring{C}}(\lie{g})$, we obtain $\twist{a}(XY) = \twist{a}(X)
    \twist{a}(Y)$ for all $X, Y \in U_{\ring{C}}(\lie{g})$. From the
    normalization $\twist{a}(\Unit) = 1$ and the unitarity
    $\twist{a}(X^*) = \cc{\twist{a}(X)}$ we finally see that any
    $\twist{a} \in \group{U}(H, \ring{C})$ is given by a
    $^*$-homomorphism $\twist{a}: U_{\ring{C}}(\lie{g})
    \longrightarrow \ring{C}$. Thus the group $\group{U}(H, \ring{C})$
    coincides with the unitary $\ring{C}$-valued characters of
    $\lie{g}$.
\end{example}

%
%

\appendix

%
%

\section{The groups $\group{GL}(H, \mathcal{A})$, $\group{GL}_0(H,
  \mathcal{A})$, $\group{U}(H, \mathcal{A})$ and $\group{U}_0(H,
  \mathcal{A})$}
\label{sec:FunnyGroup}

In this appendix we shall describe several groups naturally associated
to any \emph{unital} $^*$-algebra with a $^*$-action of a Hopf
$^*$-algebra on it.

%
%

\subsection{Definitions and fundamental properties}
\label{subsec:DefinitionFundamentaProperties}

Recall that on $\Hom_{\ring{C}}(H, \mathcal{A})$ one has the
associative \textdef{convolution product}
\begin{equation}
    \label{eq:ConvolutionProduct}
    (\twist{a} * \twist{b})(g) 
    = \twist{a}(g_\sweedler{1}) \twist{b}(g_\sweedler{2})
\end{equation}
where $\twist{a}, \twist{b} \in \Hom_{\ring{C}}(H,
\mathcal{A})$. Since $\mathcal{A}$ is unital and $H$ counital,
$\Hom_{\ring{C}}(H, \mathcal{A})$ is known to be unital with unit
$\twist{e}$ given by
\begin{equation}
    \label{eq:UnitConvolution}
    \twist{e}(g) = \epsilon(g) \Unit_{\mathcal{A}},
\end{equation}
see e.g.~\cite{kassel:1995a,majid:1995a}. We are now interested in
particular subgroups of the group of all invertible elements in the
convolution algebra $\Hom_{\ring{C}}(H, \mathcal{A})$.
\begin{definition}
    \label{definition:GLHAandUHA}
    An element $\twist{a} \in \Hom_{\ring{C}}(H, \mathcal{A})$ belongs
    to $\group{GL}(H, \mathcal{A})$ if for all $g, h \in H$ and $b \in
    \mathcal{A}$ we have
    \begin{compactenum}
    \item $\twist{a}(\Unit_H) = \Unit_\mathcal{A}$ (normalization),
    \item $\twist{a}(gh) = \twist{a}(g_\sweedler{1}) (g_\sweedler{2}
        \acts \twist{a}(h))$ (action condition),
    \item $(g_\sweedler{1} \acts b) \twist{a}(g_\sweedler{2}) =
        \twist{a}(g_\sweedler{1}) (g_\sweedler{2} \acts b)$ (module
        condition),
    \end{compactenum}
    and it belongs to $\group{U}(H, \mathcal{A})$ if in addition
    \begin{compactenum}
        \addtocounter{enumi}{3}
    \item $\twist{a}(g_\sweedler{1})
        \left(\twist{a}(S(g_\sweedler{2})^*)\right)^* = \epsilon(g)
        \Unit_{\mathcal{A}}$ (unitarity condition).
    \end{compactenum}
\end{definition}
The `action condition' can also be interpreted as a cocycle condition
while the `module condition' expresses a certain centrality property
of the values $\twist{a}(g) \in \mathcal{A}$. The subsets
$\group{U}(H, \mathcal{A}) \subseteq \group{GL}(H, \mathcal{A})$ turn
out to be subgroups of the group of invertible elements
$\group{GL}(\Hom_{\ring{C}}(H, \mathcal{A}), *)$:
\begin{proposition}
    \label{proposition:GLUGroups}
    The set $\group{GL}(H, \mathcal{A})$ becomes a group with respect
    to the convolution product $*$ and $\group{U}(H, \mathcal{A})$ is
    a subgroup. The inverse of $\twist{a} \in \group{GL}(H,
    \mathcal{A})$ is explicitly given by
    \begin{equation}
        \label{eq:InverseTwist}
        \twist{a}^{-1}(g) 
        = g_\sweedler{2} \acts \twist{a}(S^{-1}(g_\sweedler{1})).
    \end{equation}
\end{proposition}
\begin{proof}
    Clearly $\twist{e} \in \group{U}(H, \mathcal{A}) \subseteq
    \group{GL}(H, \mathcal{A})$ and $*$ is associative. Now let
    $\twist{a}, \twist{b} \in \group{GL}(H, \mathcal{A})$. Then
    $\twist{a} * \twist{b}$ fulfills the normalization
    condition. Moreover, by \textit{ii.)} and \textit{iii.)}
    \[
    \begin{split}
        (\twist{a} * \twist{b})(gh)
        &= 
        \twist{a}(g_\sweedler{1} h_\sweedler{1})
        \twist{b}(g_\sweedler{2} h_\sweedler{2}) \\
        &= 
        \twist{a}(g_\sweedler{1}) 
        (g_\sweedler{2} \acts \twist{a}(h_\sweedler{1}))
        \twist{b}(g_\sweedler{3})
        (g_\sweedler{4} \acts \twist{b}(h_\sweedler{2})) \\
        &= 
        \twist{a}(g_\sweedler{1}) 
        \twist{b}(g_\sweedler{2})
        (g_\sweedler{3} \acts \twist{a}(h_\sweedler{1}))
        (g_\sweedler{4} \acts \twist{b}(h_\sweedler{2})) \\
        &= 
        \twist{a}(g_\sweedler{1}) 
        \twist{b}(g_\sweedler{2})
        (g_\sweedler{3} \acts
        (\twist{a}(h_\sweedler{1}) \twist{b}(h_\sweedler{2}))) \\
        &= 
        (\twist{a} * \twist{b})(g_\sweedler{1})
        (g_\sweedler{2} \acts (\twist{a} * \twist{b})(h)),
    \end{split}
    \]
    whence $\twist{a} * \twist{b}$ fulfills the action condition. For
    the module condition we compute
    \[
    \begin{split}
        (g_\sweedler{1} \acts c) 
        ((\twist{a} * \twist{b})(g_\sweedler{2}))
        &= 
        (g_\sweedler{1} \acts c)
        \twist{a}(g_\sweedler{2}) \twist{b}(g_\sweedler{3}) \\
        &= 
        \twist{a}(g_\sweedler{1}) 
        (g_\sweedler{2} \acts c) \twist{b}(\sweedler{3}) \\
        &= 
        \twist{a}(g_\sweedler{1}) 
        \twist{b}(\sweedler{2})(g_\sweedler{3} \acts c) \\
        &= 
        ((\twist{a} * \twist{b})(g_\sweedler{1}))
        (g_\sweedler{2} \acts c),
    \end{split}
    \]
    whence $\twist{a} * \twist{b} \in \group{GL}(H, \mathcal{A})$,
    indeed. Now let $\twist{a}^{-1} \in \Hom_{\ring{C}}(H,
    \mathcal{A})$ be defined as in \eqref{eq:InverseTwist} then
    $\twist{a}^{-1}$ satisfies the normalization condition. For the
    action condition we compute using $S \otimes S \circ \Delta^\op =
    \Delta \circ S$
    \[
    \begin{split}
        \twist{a}^{-1}(g_\sweedler{1}) 
        (g_\sweedler{2} \acts \twist{a}^{-1}(h))
        &=
        (g_\sweedler{2} \acts \twist{a}(S^{-1}(g_\sweedler{1})))
        \left(
            (g_\sweedler{3} h_\sweedler{2}) \acts
            \twist{a}(S^{-1}(h_\sweedler{1}))
        \right) \\
        &= 
        (g_\sweedler{4} h_\sweedler{3}) \acts
        \left(
            \left(
                S^{-1}(g_\sweedler{3} h_\sweedler{2})
                \acts
                (g_\sweedler{2} \acts
                \twist{a}(S^{-1}(g_\sweedler{1})))
            \right)
            \twist{a}(S^{-1}(h_\sweedler{1}))
        \right) \\
        &= (g_\sweedler{2} h_\sweedler{3}) \acts
        \left(
            \left(
                S^{-1}(h_\sweedler{2}) \acts 
                \twist{a}(S^{-1}(g_\sweedler{1})) 
            \right)
            \twist{a}(S^{-1}(h_\sweedler{1}))
        \right) \\
        &= (g_\sweedler{2} h_\sweedler{2}) \acts
        \left(
            \left(
                S^{-1}(h_\sweedler{1})_\sweedler{1} \acts
                \twist{a}(S^{-1}(g_\sweedler{1}))
            \right)
            \twist{a}(S^{-1}(h_\sweedler{1})_\sweedler{2})
        \right) \\
        &\stackrel{\textit{iii.)}}{=}
        (g_\sweedler{2} h_\sweedler{2}) \acts
        \left(
            \twist{a}(S^{-1}(h_\sweedler{1})_\sweedler{1})
            \left(
                S^{-1}(h_\sweedler{1})_\sweedler{2} 
                \acts
                \twist{a}(S^{-1}(g_\sweedler{1}))
            \right)
        \right) \\
        &\stackrel{\textit{ii.)}}{=}
        (g_\sweedler{2} h_\sweedler{2}) \acts
        \left(
            \twist{a}
            \left(
                S^{-1}(h_\sweedler{1})
                S^{-1}(g_\sweedler{1})
            \right)
        \right) \\
        &=
        (gh)_\sweedler{2} \acts 
        \twist{a}(S^{-1}((gh)_\sweedler{1})) \\ 
        &= \twist{a}^{-1}(gh).
    \end{split}
    \]
    For the module condition we compute
    \[
    \begin{split}
        (g_\sweedler{1} \acts b) \twist{a}^{-1}(g_\sweedler{2})
        &=
        (g_\sweedler{1} \acts b) 
        \left(
            g_\sweedler{3} \acts \twist{a}(S^{-1}(g_\sweedler{2}))
        \right) \\
        &=
        g_\sweedler{4} \acts
        \left(
            (S^{-1}(g_\sweedler{3}) \acts (g_\sweedler{1} \acts b))
            \twist{a}(S^{-1}(g_\sweedler{2}))
        \right) \\
        &=
        g_\sweedler{3} \acts
        \left(
            (S^{-1}(g_\sweedler{2})_\sweedler{1} \acts 
            (g_\sweedler{1} \acts b))
            \twist{a}(S^{-1}(g_\sweedler{2})_\sweedler{2})
        \right) \\
        &\stackrel{\textit{iii.)}}{=}
        g_\sweedler{3} \acts
        \left(
            \twist{a}(S^{-1}(g_\sweedler{2})_\sweedler{1})
            \left(
                S^{-1}(g_\sweedler{2})_\sweedler{2}
                \acts (g_\sweedler{1} \acts b)
            \right)
        \right) \\
        &=
        g_\sweedler{2} \acts
        \left(
            \twist{a}(S^{-1}(g_\sweedler{1})) b
        \right) \\
        &=
        (g_\sweedler{2} \acts \twist{a}(S^{-1}(g_\sweedler{1})))
        (g_\sweedler{3} \acts b) \\
        &=
        \twist{a}^{-1}(g_\sweedler{1}) (g_\sweedler{2} \acts b),
    \end{split}
    \]
    whence $\twist{a}^{-1} \in \group{GL}(H, \mathcal{A})$ is shown.
    It remains to show that $\twist{a}^{-1}$ is the convolution
    inverse of $\twist{a}$. Indeed,
    \[
    \begin{split}
        \left(\twist{a}^{-1} * \twist{a}\right)(g)
        &= 
        (g_\sweedler{2} \acts \twist{a}(S^{-1}(g_\sweedler{1})))
        \twist{a}(g_\sweedler{3}) \\
        &\stackrel{\textit{iii.)}}{=}
        \twist{a}(g_\sweedler{2}) 
        (g_\sweedler{3} \acts \twist{a}(S^{-1}(g_\sweedler{1}))) \\
        &\stackrel{\textit{ii.)}}{=}
        \twist{a}(g_\sweedler{2} S^{-1}(g_\sweedler{1})) \\
        &\stackrel{\textit{i.)}}{=}
        \epsilon(g) \Unit_{\mathcal{A}},
    \end{split}
    \]
    and similarly for $\twist{a} * \twist{a}^{-1} = \twist{e}$. Thus
    $\group{GL}(H, \mathcal{A})$ is a group and the inverses are given
    by formula~\eqref{eq:InverseTwist}. Thus let $\twist{a}, \twist{b}
    \in \group{U}(H, \mathcal{A})$ be given.  Then
    \[
    \begin{split}
        (\twist{a} * \twist{b})(g_\sweedler{1}) 
        \left(
            (\twist{a} * \twist{b})(S(g_\sweedler{2})^*)
        \right)^*
        &=
        \twist{a}(g_\sweedler{1})
        \twist{b}(g_\sweedler{2})
        \left(\twist{b}(S(g_\sweedler{3})^*)\right)^*
        \left(\twist{a}(S(g_\sweedler{4})^*)\right)^*
        \\
        &\stackrel{\textit{iv.)}}{=}
        \twist{a}(g_\sweedler{1})
        \epsilon(g_\sweedler{2})
        \left(\twist{a}(S(g_\sweedler{3})^*)\right)^*
        \\
        &\stackrel{\textit{iv.)}}{=}
        \epsilon(g) \Unit_{\mathcal{A}}
    \end{split}
    \]
    whence $\twist{a} * \twist{b} \in \group{U}(H, \mathcal{A})$. Finally, 
    \[
    \begin{split}
        \epsilon(g) \Unit_{\mathcal{A}}
        &=
        \epsilon(g_\sweedler{1})
        \cc{\epsilon(S(g_\sweedler{2})^*)} \Unit_{\mathcal{A}} \\
        &=
        (\twist{a}^{-1} * \twist{a})(g_\sweedler{1})
        \left(
            (\twist{a}^{-1} * \twist{a}) (S(g_\sweedler{2})^*)
        \right)^* \\
        &=
        \twist{a}^{-1} (g_\sweedler{1})
        \twist{a}(g_\sweedler{2})
        \left(\twist{a}(S(g_\sweedler{3})^*)\right)^*
        \left(\twist{a}^{-1}(S(g_\sweedler{4})^*)\right)^* \\
        &\stackrel{\textit{iv.)}}{=}
        \twist{a}^{-1} (g_\sweedler{1})
        \epsilon(g_\sweedler{2})
        \left(\twist{a}^{-1}(S(g_\sweedler{3})^*)\right)^* \\
        &=
        \twist{a}^{-1} (g_\sweedler{1})
        \left(\twist{a}^{-1}(S(g_\sweedler{2})^*)\right)^*
    \end{split}
    \]
    shows $\twist{a}^{-1} \in \group{U}(H, \mathcal{A})$ as well. This
    completes the proof.
\end{proof}

Note that the group $\group{GL}(H, \mathcal{A})$ is defined for any
action of a Hopf algebra $H$ on an unital associative algebra as long
as the antipode of $H$ is \emph{invertible}. For $\group{U}(H,
\mathcal{A})$ we need the $^*$-involutions.

The next proposition describes how certain central elements of
$\mathcal{A}$ contribute to $\group{GL}(H, \mathcal{A})$ and
$\group{U}(H, \mathcal{A})$, respectively. We denote by
$\group{GL}(\zentrum(\mathcal{A}))$ the abelian group of
\textdef{invertible central elements} in $\mathcal{A}$ and
$\group{U}(\zentrum(\mathcal{A}))$ denotes the subgroup of
\textdef{unitary central elements}. Moreover,
$\group{GL}(\zentrum(\mathcal{A}))^H$ and
$\group{U}(\zentrum(\mathcal{A}))^H$ denote the $H$-invariant elements
in $\group{GL}(\zentrum(\mathcal{A}))$ and
$\group{U}(\zentrum(\mathcal{A}))$, respectively, which are subgroups.
\begin{proposition}
    \label{proposition:CenterContributionToGLUHA}
    Let $c \in \group{GL}(\zentrum(\mathcal{A}))$ then
    \begin{equation}
        \label{eq:HatcDef}
        \hat{c}(g) = c (g \acts c^{-1})
    \end{equation}
    defines an element $\hat{c} \in \group{GL}(H, \mathcal{A})$ and $c
    \mapsto \hat{c}$ is a group homomorphism such that
    \begin{equation}
        \label{eq:GLZHGLZGLHAExakt}
        1 
        \longrightarrow \group{GL}(\zentrum(\mathcal{A}))^H
        \longrightarrow \group{GL}(\zentrum(\mathcal{A}))
        \stackrel{\hat{\;}}{\longrightarrow}
        \group{GL}(H, \mathcal{A})
    \end{equation}
    is exact. Similarly, for $c \in \group{U}(\zentrum(\mathcal{A}))$
    we have $\hat{c} \in \group{U}(H, \mathcal{A})$ and
    \begin{equation}
        \label{eq:UZHUZUHAExakt}
        1 
        \longrightarrow \group{U}(\zentrum(\mathcal{A}))^H
        \longrightarrow \group{U}(\zentrum(\mathcal{A}))
        \stackrel{\hat{\;}}{\longrightarrow}
        \group{U}(H, \mathcal{A})
    \end{equation}
    is an exact sequence of group homomorphisms. Moreover, the image
    of $\group{GL}(\zentrum(\mathcal{A}))$ under $\hat{\;}$ is in the
    center of $\group{GL}(H, \mathcal{A})$.
\end{proposition}
\begin{proof}
    First we check that $\hat{c} \in \group{GL}(H, \mathcal{A})$. The
    normalization is clear. For the action condition we compute
    \[
    \begin{split}
        \hat{c}(g_\sweedler{1}) (g_\sweedler{2} \acts \hat{c}(h))
        &=
        c (g_\sweedler{1} \acts c^{-1}) 
        (g_\sweedler{2} \acts (c (h \acts c^{-1}))) \\
        &=
        c \left(g_\sweedler{3} \acts 
            \left(
                (S^{-1}(g_\sweedler{2} g_\sweedler{1}) \acts c^{-1})
                c (h \acts c^{-1})
            \right)
        \right) \\
        &=
        c (g \acts (c^{-1} c (h \acts c^{-1}))) \\
        &= 
        c ((gh) \acts c^{-1}) \\
        &= 
        \hat{c}(gh).
    \end{split}
    \]
    The module condition is shown by
    \[
    \begin{split}
        (g_\sweedler{1} \acts b) \hat{c}(g_\sweedler{2})
        &=
        (g_\sweedler{1} \acts b) c (g_\sweedler{2} \acts c^{1}) \\
        &=
        c (g \acts (bc)^{-1}) \\
        &=
        c (g_\sweedler{1} \acts c^{-1}) (g_\sweedler{2} \acts b) \\
        &= 
        \hat{c}(g_\sweedler{1}) (g_\sweedler{2} \acts b),
    \end{split}
    \]
    whence $\hat{c} \in \group{GL}(H, \mathcal{A})$. Clearly
    $\widehat{\Unit_\mathcal{A}} = \twist{e}$ and for $c, d \in
    \group{GL}(\zentrum(\mathcal{A}))$ we have
    \[
    \widehat{cd}(g)
    = cd (g_\sweedler{1} \acts c^{-1}) (g_\sweedler{2} \acts d^{-1})
    = c(g_\sweedler{1} \acts c^{-1}) d (g_\sweedler{2} \acts d^{-1})
    = \hat{c}(g_\sweedler{1}) \hat{d}(g_\sweedler{2})
    = (\hat{c} * \hat{d})(g),
    \]
    whence $\hat{\;}$ is a group morphism. If $c$ is $H$-invariant it
    is easy to see that $\hat{c} = \twist{e}$. Conversely, if $\hat{c}
    = \twist{e}$, then $c (g \acts c^{-1}) = \epsilon(g)$ whence $g
    \acts c^{-1} = \epsilon(g) c^{-1}$. Hence $c^{-1}$ and thus $c$ is
    $H$-invariant. This proves the exactness of
    \eqref{eq:GLZHGLZGLHAExakt}. Now let $\twist{a} \in \group{GL}(H,
    \mathcal{A})$ be arbitrary and $c \in
    \group{GL}(\zentrum(\mathcal{A}))$. Then using the centrality of
    $c$ as well as \textit{iii.)} for $\twist{a}$ we get
    \[
    (\twist{a} * \hat{c})(g)
    = \twist{a}(g_\sweedler{1}) c (g_\sweedler{2} \acts c^{-1})
    = c \twist{a}(g_\sweedler{1})(g_\sweedler{2} \acts c^{-1})
    = c (g_\sweedler{1} \acts c^{-1})\twist{a}(g_\sweedler{2})
    = (\hat{c} * \twist{a})(g),
    \]
    whence $\hat{c}$ is central in $\group{GL}(H, \mathcal{A})$.
    For the last part let $c \in \group{U}(\zentrum(\mathcal{A}))$ be
    unitary. Then
    \[
    \begin{split}
        \hat{c}(g_\sweedler{1})
        \left(\hat{c}(S(g_\sweedler{2})^*)\right)^*
        &=
        c (g_\sweedler{1} \acts c^{-1})
        \left(S(g_\sweedler{2})^* \acts c^{-1}\right)^* c^* \\
        &=
        (g_\sweedler{1} \acts c^{-1})
        (S(S(g_\sweedler{2})^*)^* \acts c) \\
        &=
        g \acts (c^{-1} c) \\
        &= \epsilon(g) \Unit_{\mathcal{A}}
    \end{split}
    \]
    shows $\hat{c} \in \group{U}(H, \mathcal{A})$. The remaining
    statements follow easily.
\end{proof}

Thus we can divide by the image of $\group{GL}(\zentrum(\mathcal{A}))$
under $\hat{\;}$ and obtain the quotient groups
\begin{equation}
    \label{eq:GLNullDef}
    \group{GL}_0(H, \mathcal{A}) 
    = \group{GL}(H, \mathcal{A}) \big/
    \widehat{\group{GL}(\zentrum(\mathcal{A}))}
\end{equation}
and
\begin{equation}
    \label{eq:UNullDef}
    \group{U}_0(H, \mathcal{A}) 
    = \group{U}(H, \mathcal{A}) \big/
    \widehat{\group{U}(\zentrum(\mathcal{A}))}.
\end{equation}
This way, we can complete the sequences \eqref{eq:GLZHGLZGLHAExakt}
and \eqref{eq:UZHUZUHAExakt} to the exact sequences
\begin{equation}
    \label{eq:TheGLExactSequenz}
    1 
    \longrightarrow \group{GL}(\zentrum(\mathcal{A}))^H
    \longrightarrow \group{GL}(\zentrum(\mathcal{A}))
    \stackrel{\hat{\;}}{\longrightarrow} \group{GL}(H, \mathcal{A})
    \longrightarrow \group{GL}_0(H, \mathcal{A})
    \longrightarrow
    1
\end{equation}
and
\begin{equation}
    \label{eq:TheUExactSequenz}
    1 
    \longrightarrow \group{U}(\zentrum(\mathcal{A}))^H
    \longrightarrow \group{U}(\zentrum(\mathcal{A}))
    \stackrel{\hat{\;}}{\longrightarrow}\group{U}(H, \mathcal{A})
    \longrightarrow \group{U}_0(H, \mathcal{A})
    \longrightarrow
    1.
\end{equation}
\begin{remark}
    \label{remark:TrivialesZentrum}
    If the center $\zentrum(\mathcal{A}) = \ring{C}
    \Unit_{\mathcal{A}}$ is \emph{trivial} then
    $\group{GL}(\zentrum(\mathcal{A}))^H =
    \group{GL}(\zentrum(\mathcal{A}))$ as well as
    $\group{U}(\zentrum(\mathcal{A}))^H =
    \group{U}(\zentrum(\mathcal{A}))$. Thus we have in this case
    \begin{equation}
        \label{eq:TrivialesZentrumGLUGruppen}
        \group{GL}(H, \mathcal{A}) = \group{GL}_0(H, \mathcal{A})
        \quad
        \textrm{and}
        \quad
        \group{U}(H, \mathcal{A}) = \group{U}_0(H, \mathcal{A}).
    \end{equation}
\end{remark}

The groups $\group{GL}(H, \mathcal{A})$, $\group{GL}_0(H,
\mathcal{A})$, $\group{U}(H, \mathcal{A})$ and $\group{U}_0(H,
\mathcal{A})$ enjoy nice functorial properties which we shall discuss
now. Under general homomorphisms or $^*$-homomorphisms, respectively,
we can not conclude any good behavior of elements in $\group{GL}(H,
\mathcal{A})$ or $\group{U}(H, \mathcal{A})$, respectively, as the
module condition requires information about commutation relations with
\emph{arbitrary} algebra elements. Nevertheless, for \emph{surjective}
homomorphisms we have the following statement expressing the
functorial properties:
\begin{proposition}
    \label{proposition:GLUFunktoriell}
    Let $\phi: \mathcal{A} \longrightarrow \mathcal{B}$ be a
    $H$-equivariant surjective homomorphism.
    \begin{compactenum}
    \item For any $\twist{a} \in \group{GL}(H, \mathcal{A})$ we have
        $\phi_* \twist{a} = \phi \circ \twist a \in \group{GL}(H,
        \mathcal{B})$.
    \item The map $\phi_*: \group{GL}(H, \mathcal{A}) \longrightarrow
        \group{GL}(H, \mathcal{B})$ is a group homomorphism.
    \item If $\psi: \mathcal{B} \longrightarrow \mathcal{C}$ is
        another $H$-equivariant surjective homomorphism then
        \begin{equation}
            \label{eq:phipsistarpsistarpsistar}
            (\psi \circ \phi)_* = \psi_* \circ \phi_*
            \quad
            \textrm{and}
            \quad
            (\id_{\mathcal{A}})_* = \id_{\group{GL}(H, \mathcal{A})}.
        \end{equation}
    \item The group homomorphism $\phi_*: \group{GL}(H, \mathcal{A})
        \longrightarrow \group{GL}(H, \mathcal{B})$ induces a group
        homomorphism $\group{GL}_0(H, \mathcal{A}) \longrightarrow
        \group{GL}_0(H, \mathcal{B})$, also denoted by $\phi_*$, such
        that the diagram
        \begin{equation}
            \label{eq:GLZAHZAHANullCommutes}
            \bfig
            \hSquares|rrrrrrr|/>`>``>`>`>`>/%
            [1`\group{GL}(\zentrum(\mathcal{A}))^H%
            `\group{GL}(\zentrum(\mathcal{A}))%
            `1%
            `\group{GL}(\zentrum(\mathcal{B}))^H%
            `\group{GL}(\zentrum(\mathcal{B}))%
            ;```\phi`\phi``]
            \morphism(1680,0)<300,0>[`;]
            \morphism(1680,500)<300,0>[`;]
            \hSquares(2230,0)|rrrrrrr|/>`>`>`>``>`>/%
            [\group{GL}(H, \mathcal{A})%
            `\group{GL}_0(H,\mathcal{A})%
            `1%
            `\group{GL}(H,\mathcal{B})%
            `\group{GL}_0(H,\mathcal{B})%
            `1%
            ;``\phi_*`\phi_*```]
            \efig
        \end{equation}
        commutes.
    \item If $\phi$ is in addition a $^*$-homomorphism we can replace
        `\/$\group{GL}$' by `\/$\group{U}$' everywhere. In particular we
        have a commutative diagram
        \begin{equation}
            \label{eq:UZAHAHANullCommutes}
            \bfig
            \hSquares|rrrrrrr|/>`>``>`>`>`>/%
            [1`\group{U}(\zentrum(\mathcal{A}))^H%
            `\group{U}(\zentrum(\mathcal{A}))%
            `1%
            `\group{U}(\zentrum(\mathcal{B}))^H%
            `\group{U}(\zentrum(\mathcal{B}))%
            ;```\phi`\phi``]
            \morphism(1550,0)<300,0>[`;]
            \morphism(1550,500)<300,0>[`;]
            \hSquares(2050,0)|rrrrrrr|/>`>`>`>``>`>/%
            [\group{U}(H, \mathcal{A})%
            `\group{U}_0(H,\mathcal{A})%
            `1%
            `\group{U}(H,\mathcal{B})%
            `\group{U}_0(H,\mathcal{B})%
            `1%
            ;``\phi_*`\phi_*```]
            \efig.
        \end{equation}
    \end{compactenum}
\end{proposition}
\begin{proof}
    The first part is a simple verification of the axioms and the
    second part follows from well-known properties of the convolution
    product. The third part is obvious. For the fourth part we have to
    show that the first three vertical arrows give commuting diagrams
    since the last arrow is induced precisely in the way that
    \eqref{eq:GLZAHZAHANullCommutes} commutes in total. For the first
    box in \eqref{eq:GLZAHZAHANullCommutes} this is obvious. Let $c
    \in \group{GL}(\zentrum(\mathcal{A}))$ then
    \[
    \left(\phi_* \hat{c}\right) (g) 
    = \phi(c (g \acts c^{-1}))
    = \phi(c) (g \acts \phi(c)^{-1})
    = \widehat{\phi(c)}(g)
    \]
    shows the commutativity of the second box in
    \eqref{eq:GLZAHZAHANullCommutes}. Note that since
    $\phi(\Unit_\mathcal{A}) = \Unit_\mathcal{B}$ we indeed have
    $\phi(c^{-1}) = \phi(c)^{-1}$. Then the last part is again a
    simple consequence of the fact that
    $\phi(\group{U}(\zentrum(\mathcal{A}))) \subseteq
    \group{U}(\zentrum(\mathcal{B}))$.
\end{proof}

In a slightly more fancy way we can rephrase the content of the
proposition as follows:
\begin{corollary}
    \label{corollary:IsoActsOfSequence}
    The groupoid of $H$-equivariant isomorphisms $\IsoH$
    acts on the exact sequence \eqref{eq:TheGLExactSequenz} by
    isomorphisms,
    whence in particular the whole exact sequence of groups together
    with its $\AutH(\mathcal{A})$-action on it is an invariant of
    $\mathcal{A}$ as associative algebra with $H$-action. Analogously,
    the groupoid $\starIsoH$ acts on the exact sequence
    \eqref{eq:TheUExactSequenz} by isomorphisms,
    whence the exact sequence \eqref{eq:TheUExactSequenz} together
    with its $\starAutH(\mathcal{A})$-action on it is an invariant of
    $\mathcal{A}$ as $^*$-algebra with $^*$-action of $H$.
\end{corollary}

%
%

\subsection{The cocommutative case}
\label{subsec:CocommutativeCaseGroupsGLUNull}

We specialize now for a cocommutative Hopf $^*$-algebra $H$. In this
case the situation simplifies as follows:
\begin{proposition}
    \label{proposition:HcocomGLGLnullUUnull}
    Let $H$ be cocommutative.
    \begin{compactenum}
    \item $H \acts \zentrum(\mathcal{A}) \subseteq
        \zentrum(\mathcal{A})$ whence $\zentrum(\mathcal{A})$ inherits
        the $H$-action.
    \item For $\twist{a} \in \group{GL}(H, \mathcal{A})$ we have
        $\twist{a}(g) \in \zentrum(\mathcal{A})$ for all $g \in H$
        whence
        \begin{equation}
            \label{eq:GLGLZentrum}
            \group{GL}(H, \mathcal{A}) 
            = \group{GL}(H, \zentrum(\mathcal{A}))
            \quad
            \textrm{and}
            \quad
            \group{GL}_0(H, \mathcal{A}) 
            = \group{GL}_0(H, \zentrum(\mathcal{A}))
        \end{equation}
        and analogously
        \begin{equation}
            \label{eq:UUZentrum}
            \group{U}(H, \mathcal{A}) 
            = \group{U}(H, \zentrum(\mathcal{A}))
            \quad
            \textrm{and}
            \quad
            \group{U}_0(H, \mathcal{A}) 
            = \group{U}_0(H, \zentrum(\mathcal{A})).
        \end{equation}
    \item The groups $\group{GL}(H, \mathcal{A})$, $\group{GL}_0(H,
        \mathcal{A})$, $\group{U}(H, \mathcal{A})$ and $\group{U}_0(H,
        \mathcal{A})$ are abelian.
    \item The space of unital algebra homomorphisms $H \longrightarrow
        \zentrum(\mathcal{A})^H$ is a subgroup of $\group{GL}(H, 
        \mathcal{A})$  and the space of unital $^*$-homomorphisms $H
        \longrightarrow \zentrum(\mathcal{A})^H$ is a subgroup of
        $\group{U}(H, \mathcal{A})$. The inverse of such a
        homomorphism $\twist{a}: H \longrightarrow
        \zentrum(\mathcal{A})$ is explicitly given by
        \begin{equation}
            \label{eq:InverseHomomorphHtoZAH}
            \twist{a}^{-1}(g) 
            = \twist{a}(S(g)) 
            = \twist{a}(S^{-1}(g)).
        \end{equation}
    \end{compactenum}
\end{proposition}
\begin{proof}
    The first statement is well-known. For the second we compute
    \[
    \twist{a}(g) b 
    = 
    \twist{a}(g_\sweedler{1}) \epsilon(g_\sweedler{2}) b
    =
    \twist{a}(g_\sweedler{1}) 
    ((g_\sweedler{2} S(g_\sweedler{3})) \acts b)
    =
    \left(g_\sweedler{1} \acts ( S(g_\sweedler{3}) \acts b)\right)
    \twist{a}(g_\sweedler{2})
    = \epsilon(g_\sweedler{1}) b \twist{a}(g_\sweedler{2})
    = b \twist{a}(g),
    \]
    using the module condition for $\twist{a}$ as well as the
    cocommutativity. The third part is then a simple consequence. For
    the fourth part part we consider a unital homomorphism $\twist{a}:
    H \longrightarrow \zentrum(\mathcal{A})^H$. Then
    $\twist{a}(\Unit_H) = \Unit_{\mathcal{A}}$ by definition and
    \[
    \twist{a}(gh) 
    = \twist{a}(g) \twist{a}(h)
    = \twist{a}(g_\sweedler{1}) \epsilon(g_\sweedler{2}) \twist{a}(h)
    =\twist{a}(g_\sweedler{1}) (g_\sweedler{2} \acts \twist{a}(h)),
    \]
    since $\twist{a}(h)$ is $H$-invariant. Moreover, 
    \[
    (g_\sweedler{1} \acts b) \twist{a}(g_\sweedler{2})
    = \twist{a}(g_\sweedler{2}) (g_\sweedler{1} \acts b)
    = \twist{a}(g_\sweedler{1}) (g_\sweedler{2}) \acts b),
    \]
    since $\twist{a}(g_\sweedler{2})$ is central and $H$ is
    cocommutative. Finally, if $\twist{a}$ is even a $^*$-homomorphism
    then
    \[
    \twist{a}(g_\sweedler{1})
    \left(\twist{a}(S(g_\sweedler{2})^*)\right)^*
    = 
    \twist{a}(g_\sweedler{1})\twist{a}(S(g_\sweedler{2}))
    =
    \twist{a}(g_\sweedler{1}S(g_\sweedler{2}))
    =
    \epsilon(g) \Unit_{\mathcal{A}}.
    \]
    Now let $\twist{a}, \twist{b}: H \longrightarrow
    \zentrum(\mathcal{A})^H$ be unital homomorphisms then
    a simple computation shows that $\twist{a} * \twist{b}$ is again a
    unital homomorphism taking its values in
    $\zentrum(\mathcal{A})^H$. Moreover, the general formula for the
    inverse \eqref{eq:InverseTwist} leads to
    \[
    \twist{a}^{-1}(g) 
    = \epsilon(g_\sweedler{1})\twist{a}(S^{-1}(g_\sweedler{2}))
    = \twist{a}(S^{-1}(g))
    = \twist{a}(S(g))
    \]
    since $\twist{a}(g)$ is invariant and $S^2 = \id$ in the
    cocommutative case. It is easy to see that $\twist{a}^{-1}$ is
    still a homomorphism since $S$ is an antihomomorphism and the
    images all commute. If $\twist{a}$ is even a $^*$-homomorphism
    then $\twist{a}^{-1}$ is a $^*$-homomorphism, too, since $S$
    commutes with the $^*$-involution in the cocommutative case. This
    completes the proof.
\end{proof}

In particular, the unital algebra homomorphisms
\begin{equation}
    \label{eq:Characters}
    \chi: H \longrightarrow \ring{C},
\end{equation}
i.e.~the \textdef{characters} of $H$, always contribute to
$\group{GL}(H, \mathcal{A})$ by setting $\twist{a}^\chi(g) = \chi(g)
\Unit_{\mathcal{A}}$. If $\chi$ is in addition a $^*$-homomorphism we
call $\chi$ a \textdef{unitary character}.  In fact, if the center of
$\mathcal{A}$ is trivial the characters of $H$ give the whole group
$\group{GL}(H, \mathcal{A})$:
\begin{proposition}
    \label{proposition:HCocomCharactersGLU}
    Let $H$ be cocommutative.
    \begin{compactenum}
    \item The characters of $H$ constitute a subgroup of
        $\group{GL}(H, \mathcal{A})$ via $\chi \mapsto \twist{a}^\chi$
        and the unitary characters of $H$ constitute a subgroup of
        $\group{U}(H, \mathcal{A})$.
    \item If the center of $\mathcal{A}$ is trivial,
        $\zentrum(\mathcal{A}) = \ring{C}\Unit_{\mathcal{A}}$, then
        any element of $\group{GL}(H, \mathcal{A}) = \group{GL}_0(H,
        \mathcal{A})$ is a character and any element of $\group{U}(H,
        \mathcal{A}) = \group{U}_0(H, \mathcal{A})$ is a unitary
        character.
    \end{compactenum}
\end{proposition}
\begin{proof}
    First it is clear that $\twist{a}^\chi$ is a unital algebra
    homomorphism $H \longrightarrow \zentrum(\mathcal{A}) =
    \zentrum(\mathcal{A})^H = \ring{C}\Unit_{\mathcal{A}}$ and thus an
    element of $\group{G}(H, \mathcal{A})$. If in addition
    $\cc{\chi(g)} = \chi(g^*)$ then $\twist{a}^\chi$ is a
    $^*$-homomorphism and hence $\twist{a}^\chi \in \group{U}(H,
    \mathcal{A})$. Since $\chi * \chi'$ is clearly a character if
    $\chi$, $\chi'$ are characters and since $\twist{a}^\chi *
    \twist{a}^{\chi'} = \twist{a}^{\chi * \chi'}$ we see that the
    elements of the form $\twist{a}^\chi$ are closed under
    multiplication in $\group{GL}(H, \mathcal{A})$. Moreover, denoting
    the `inverse character' of $\chi$ by $\chi^{-1}(g) = \chi(S(g))$
    we see that $(\twist{a}^\chi)^{-1} = \twist{a}^{\chi^{-1}}$ whence
    the characters are a subgroup indeed. If $\twist{a} \in
    \group{GL}(H, \mathcal{A})$ is of the form $\twist{a}(g) = \chi(g)
    \Unit_{\mathcal{A}}$ for any $g \in H$ with some $\chi(g) \in
    \ring{C}$ then $g \mapsto \chi(g)$ is necessarily a character. The
    unitary case is treated analogously.  Then the second part follows
    from $\group{GL}(H, \mathcal{A}) = \group{GL}_0(H, \mathcal{A})$
    and $\group{U}(H, \mathcal{A}) = \group{U}_0(H, \mathcal{A})$ by
    Remark~\ref{remark:TrivialesZentrum}.
\end{proof}

This statement allows to construct easily $^*$-algebras $\mathcal{A}$
with $^*$-actions of a cocommutative Hopf $^*$-algebra $H$ such that
the groups $\group{GL}_0(H, \mathcal{A})$ and $\group{U}_0(H,
\mathcal{A})$ are non-trivial. Note however, that if the center
$\zentrum(\mathcal{A})$ is non-trivial then it may well happen that
characters of $H$, viewed as non-trivial elements of $\group{GL}(H,
\mathcal{A})$, are killed when passing to the quotient group
$\group{GL}_0(H, \mathcal{A})$. Hence in general $\group{GL}(H,
\mathcal{A}) \ne \group{GL}_0(H, \mathcal{A})$ as well as
$\group{U}(H, \mathcal{A}) \ne \group{U}_0(H, \mathcal{A})$.

%
%

\begin{footnotesize}

\begin{thebibliography}{10}

\bibitem {ara:1999a}
{\sc Ara, P.: }\newblock {\em {M}orita equivalence for rings with involution}.
\newblock Alg. Rep. Theo.  {\bf 2} (1999), 227--247.

\bibitem {ara:1999b}
{\sc Ara, P.: }\newblock {\em Morita equivalence and Pedersen Ideals}.
\newblock Proc. AMS  {\bf 129}.4 (2000), 1041--1049.

\bibitem {arnal.cortet.molin.pinczon:1983a}
{\sc Arnal, D., Cortet, J.~C., Molin, P., Pinczon, G.: }\newblock {\em
  Covariance and Geometrical Invariance in {$*$-Quantization}}.
\newblock J. Math. Phys.  {\bf 24}.2 (1983), 276--283.

\bibitem {bass:1968a}
{\sc Bass, H.: }\newblock {\em Algebraic ${K}$-theory}.
\newblock W. A. Benjamin, Inc., New York, Amsterdam, 1968.

\bibitem {bayen.et.al:1978a}
{\sc Bayen, F., Flato, M., Fr{{\o}}nsdal, C., Lichnerowicz, A., Sternheimer,
  D.: }\newblock {\em Deformation Theory and Quantization}.
\newblock Ann. Phys.  {\bf 111} (1978), 61--151.

\bibitem {benabou:1967a}
{\sc B{\'e}nabou, J.: }\newblock {\em Introduction to Bicategories}.
\newblock In: {\em Reports of the Midwest Category Seminar},   1--77.
  Springer-Verlag, 1967.

\bibitem {bursztyn:2002a}
{\sc Bursztyn, H.: }\newblock {\em Semiclassical geometry of quantum line
  bundles and {M}orita equivalence of star products}.
\newblock Int. Math. Res. Not.  {\bf 2002}.16 (2002), 821--846.

\bibitem {bursztyn.waldmann:2001b}
{\sc Bursztyn, H., Waldmann, S.: }\newblock {\em {$^*$}-Ideals and Formal
  Morita Equivalence of {$^*$}-Algebras}.
\newblock Int. J. Math.  {\bf 12}.5 (2001), 555--577.

\bibitem {bursztyn.waldmann:2001a}
{\sc Bursztyn, H., Waldmann, S.: }\newblock {\em Algebraic Rieffel Induction,
  Formal Morita Equivalence and Applications to Deformation Quantization}.
\newblock J. Geom. Phys.  {\bf 37} (2001), 307--364.

\bibitem {bursztyn.waldmann:2002a}
{\sc Bursztyn, H., Waldmann, S.: }\newblock {\em The characteristic classes of
  {M}orita equivalent star products on symplectic manifolds}.
\newblock Commun. Math. Phys.  {\bf 228} (2002), 103--121.

\bibitem {bursztyn.waldmann:2003a:pre}
{\sc Bursztyn, H., Waldmann, S.: }\newblock {\em Completely positive inner
  products and strong {M}orita equivalence}.
\newblock Preprint (FR-THEP 2003/12)  {\bf math.QA/0309402} (September 2003),
  36 pages.
\newblock To appear in Pacific J. Math.

\bibitem {bursztyn.waldmann:2004a}
{\sc Bursztyn, H., Waldmann, S.: }\newblock {\em Bimodule deformations,
  {P}icard groups and contravariant connections}.
\newblock K-Theory  {\bf 31} (2004), 1--37.

\bibitem {bursztyn.weinstein:2004a}
{\sc Bursztyn, H., Weinstein, A.: }\newblock {\em Picard groups in {P}oisson
  geometry}.
\newblock Moscow Math. J.  {\bf 4} (2004), 39--66.

\bibitem {bursztyn.weinstein:2004a:pre}
{\sc Bursztyn, H., Weinstein, A.: }\newblock {\em Poisson geometry and Morita
  equivalence}.
\newblock Preprint  {\bf math.SG/0402347} (2004), 52 pages.
\newblock To appear in the London Math. Society Lecture Notes series.

\bibitem {cannasdasilva.weinstein:1999a}
{\sc Cannas~da Silva, A., Weinstein, A.: }\newblock {\em Geometric Models for
  Noncommutative Algebras}.
\newblock {\em Berkeley Mathematics Lecture Notes}.
\newblock AMS, 1999.

\bibitem {combes:1984a}
{\sc Combes, F.: }\newblock {\em Crossed products and {M}orita equivalence}.
\newblock Proc. London Math. Soc. (3)  {\bf 49}.2 (1984), 289--306.

\bibitem {connes:1994a}
{\sc Connes, A.: }\newblock {\em Noncommutative Geometry}.
\newblock Academic Press, San Diego, New York, London, 1994.

\bibitem {curto.muhly.williams:1984a}
{\sc Curto, R.~E., Muhly, P.~S., Williams, D.~P.: }\newblock {\em Cross
  Products of Strongly Morita Equivalent {$C^*$}-Algebras}.
\newblock Proc. Am. Math. Soc.  {\bf 90}.4 (1984), 528--530.

\bibitem {dito.sternheimer:2002a}
{\sc Dito, G., Sternheimer, D.: }\newblock {\em Deformation quantization:
  genesis, developments and metamorphoses}.
\newblock In: {\sc Halbout, G. (eds.): }\newblock {\em Deformation
  quantization}, vol.~1 in {\em IRMA Lectures in Mathematics and Theoretical
  Physics},   9--54. Walter de Gruyter, Berlin, New York, 2002.

\bibitem {echterhoff.kaliszewski.quigg.raeburn:2002a:pre}
{\sc Echterhoff, S., Kaliszewski, S., Quigg, J., Raeburn, I.: }\newblock {\em A
  Categorical Approach to Imprimitivity Theorems for {$C^*$}-Dynamical
  Systems}.
\newblock Preprint  {\bf math.OA/0205322} (2002), 152 pages.

\bibitem {gerstenhaber.schack:1988a}
{\sc Gerstenhaber, M., Schack, S.~D.: }\newblock {\em Algebraic Cohomology and
  Deformation Theory}.
\newblock In: {\sc Hazewinkel, M., Gerstenhaber, M. (eds.): }\newblock {\em
  Deformation Theory of Algebras and Structures and Applications},   13--264.
  Kluwer Academic Press, Dordrecht, 1988.

\bibitem {gutt:2000a}
{\sc Gutt, S.: }\newblock {\em Variations on deformation quantization}.
\newblock In: {\sc Dito, G., Sternheimer, D. (eds.): }\newblock {\em
  Conf{\'e}rence Mosh{\'e} Flato 1999. Quantization, Deformations, and
  Symmetries}, {\em Mathematical Physics Studies} no. {\bf 21},   217--254.
  Kluwer Academic Publishers, Dordrecht, Boston, London, 2000.

\bibitem {gutt.rawnsley:2003a}
{\sc Gutt, S., Rawnsley, J.: }\newblock {\em Natural Star Products on
  Symplectic Manifolds and Quantum Moment Maps}.
\newblock Lett. Math. Phys.  {\bf 66} (2003), 123--139.

\bibitem {jacobson:1985a}
{\sc Jacobson, N.: }\newblock {\em Basic Algebra I}.
\newblock Freeman and Company, New York, 2. edition, 1985.

\bibitem {jansen.neumaier.waldmann:2005a:pre}
{\sc Jansen, S., Neumaier, N., Waldmann, S.: }\newblock {\em Covariant {M}orita
  Equivalence of Star Products}.
\newblock In preparation.

\bibitem {kadison.ringrose:1997a}
{\sc Kadison, R.~V., Ringrose, J.~R.: }\newblock {\em Fundamentals of the
  Theory of Operator Algebras. Volume I: Elementary Theory}, vol.~15 in {\em
  Graduate Studies in Mathematics}.
\newblock American Mathematical Society, Providence, 1997.

\bibitem {kassel:1995a}
{\sc Kassel, C.: }\newblock {\em Quantum Groups}.
\newblock {\em Graduate Texts in Mathematics} no. {\bf 155}.
\newblock Springer-Verlag, New York, Berlin, Heidelberg, 1995.

\bibitem {klimyk.schmuedgen:1997a}
{\sc Klimyk, A., Schm{\"u}dgen, K.: }\newblock {\em Quantum Groups and Their
  Representations}.
\newblock {\em Texts and Monographs in Physics}.
\newblock Springer-Verlag, Heidelberg, Berlin, New York, 1997.

\bibitem {kustermans:2002a}
{\sc Kustermans, J.: }\newblock {\em Induced corepresentations of locally
  compact quantum groups}.
\newblock J. Funct. Anal.  {\bf 194} (2002), 410--459.

\bibitem {lam:1999a}
{\sc Lam, T.~Y.: }\newblock {\em Lectures on Modules and Rings}, vol. 189 in
  {\em Graduate Texts in Mathematics}.
\newblock Springer-Verlag, Berlin, Heidelberg, New York, 1999.

\bibitem {lance:1995a}
{\sc Lance, E.~C.: }\newblock {\em {H}ilbert {$C^*$}-modules. A toolkit for
  operator algebraists}, vol. 210 in {\em London Mathematical Society Lecture
  Note Series}.
\newblock Cambridge University Press, Cambridge, 1995.

\bibitem {landsman:1998a}
{\sc Landsman, N.~P.: }\newblock {\em Mathematical Topics between Classical and
  Quantum Mechanics}.
\newblock {\em Springer Monographs in Mathematics}.
\newblock Springer-Verlag, Berlin, Heidelberg, New York, 1998.

\bibitem {landsman:2001c}
{\sc Landsman, N.~P.: }\newblock {\em Bicategories of operator algebras and
  {P}oisson manifolds}.
\newblock In: {\sc Longo, R. (eds.): }\newblock {\em Mathematical physics in
  mathematics and physics}, vol.~30 in {\em Fields Inst. Commun.},   271--286.
  Amer. Math. Soc., Providence, RI, 2001.
\newblock Proceedings of the conference, dedicated to Sergio Doplicher and John
  E. Roberts on the occasion of their 60th birthday, held in Siena, June
  20--24, 2000.

\bibitem {landsman:2001b}
{\sc Landsman, N.~P.: }\newblock {\em Quantized reduction as a tensor product}.
\newblock In: {\sc Landsman, N.~P., Pflaum, M., Schlichenmaier, M. (eds.):
  }\newblock {\em Quantization of Singular Symplectic Quotients},   137--180.
  Birkh{\"a}user, Basel, Boston, Berlin, 2001.

\bibitem {majid:1995a}
{\sc Majid, S.: }\newblock {\em Foundations of quantum group theory}.
\newblock Cambridge University Press, 1995.

\bibitem {moerdijk.mrcun:2003a}
{\sc Moerdijk, I., Mr{\v{c}}un, J.: }\newblock {\em Introduction to Foliations
  and {L}ie Groupoids}.
\newblock {\em Cambridge studies in advanced mathematics} no. {\bf 91}.
\newblock Cambridge University Press, Cambridge, UK, 2003.

\bibitem {morita:1958a}
{\sc Morita, K.: }\newblock {\em Duality for modules and its applications to
  the theory of rings with minimum condition}.
\newblock Sci. Rep. Tokyo Kyoiku Daigaku Sect. A  {\bf 6} (1958), 83--142.

\bibitem {mueller-bahns.neumaier:2004a}
{\sc M{\"u}ller-Bahns, M.~F., Neumaier, N.: }\newblock {\em Some remarks on
  {$\mathfrak{g}$}-invariant Fedosov star products and quantum momentum
  mappings}.
\newblock J. Geom. Phys.  {\bf 50} (2004), 257--272.

\bibitem {raeburn.williams:1998a}
{\sc Raeburn, I., Williams, D.~P.: }\newblock {\em Morita equivalence and
  continuous-trace {$C^*$}-algebras}, vol.~60 in {\em Mathematical Surveys and
  Monographs}.
\newblock American Mathematical Society, Providence, RI, 1998.

\bibitem {rieffel:1974a}
{\sc Rieffel, M.~A.: }\newblock {\em Induced representations of
  {$C^*$}-algebras}.
\newblock Adv. Math.  {\bf 13} (1974), 176--257.

\bibitem {rieffel:1974b}
{\sc Rieffel, M.~A.: }\newblock {\em Morita equivalence for {$C^*$}-algebras
  and {$W^*$}-algebras}.
\newblock J. Pure. Appl. Math.  {\bf 5} (1974), 51--96.

\bibitem {schmuedgen:1990a}
{\sc Schm{\"{u}}dgen, K.: }\newblock {\em Unbounded Operator Algebras and
  Representation Theory}, vol.~37 in {\em Operator Theory: Advances and
  Applications}.
\newblock Birkh{\"{a}}user Verlag, Basel, Boston, Berlin, 1990.

\bibitem {schmuedgen.wagner:2003a}
{\sc Schm{\"u}dgen, K., Wagner, E.: }\newblock {\em Hilbert space
  representations of cross product algebras}.
\newblock J. Funct. Anal.  {\bf 200} (2003), 451--493.

\bibitem {vaes:2004a:pre}
{\sc Vaes, S.: }\newblock {\em A new approach to induction and imprimitivity
  results}.
\newblock Preprint  {\bf math.QA/0407525} (2004), 39 pages.

\bibitem {waldmann:2004a}
{\sc Waldmann, S.: }\newblock {\em The Picard Groupoid in Deformation
  Quantization}.
\newblock Lett. Math. Phys.  {\bf 69} (2004), 223--235.

\bibitem {waldmann:2005a}
{\sc Waldmann, S.: }\newblock {\em States and Representation Theory in
  Deformation Quantization}.
\newblock Rev. Math. Phys. \textbf{17} (2005), 15--75.

\bibitem {xu:1991a}
{\sc Xu, P.: }\newblock {\em Morita Equivalence of Poisson Manifolds}.
\newblock Commun. Math. Phys.  {\bf 142} (1991), 493--509.

\end{thebibliography}

\end{footnotesize}

\end{document}